\documentclass[10pt]{siamltex}

\usepackage{amsfonts}
\usepackage{}

\usepackage{amsfonts}
\usepackage[latin1]{inputenc}
\usepackage{graphicx,pgfplots,tikz,multirow}
\usepackage{cite}
\usepackage[T1]{fontenc}
\usepackage{pgflibraryarrows}
\usepackage{pgflibrarysnakes}
\usepackage{amsmath}
\usepackage{epsfig,graphics, color,bm}
\usepackage{cases}
\usepackage{subeqnarray}
\usepackage{amsfonts}
\usepackage{mathrsfs,bm}
\usepackage{amssymb}
\usepackage[toc,page]{appendix}
\usepackage{graphicx,epsfig,amsmath,amsfonts,amscd}
\usepackage{float}
\usepackage{color}
\usepackage{latexsym,amssymb}
\usepackage{dsfont}
\usepackage{fancyhdr}
\usepackage{stmaryrd}
\usepackage{subfig}
\usepackage{float}
\usepackage{epstopdf}
\usepackage{comment}
\usepackage{array}

\usepackage{amsmath}
\renewcommand{\theequation}{\arabic{section}.\arabic{equation}}

\makeatletter
\newcommand\figcaption{\def\@captype{figure}\caption}
\newcommand\tabcaption{\def\@captype{table}\caption}
\makeatother
\newtheorem{Property}{Property}[section]
\newtheorem{example}{Example}[section]

\newcommand{\beq}{\begin{equation}}
\newcommand{\eeq}{\end{equation}}
\newcommand{\beqn}{\begin{equation*}}
\newcommand{\eeqn}{\end{equation*}}
\newcommand{\be}{\begin{eqnarray}}
\newcommand{\ee}{\end{eqnarray}}
\newcommand{\bex}{\begin{eqnarray*}}
\newcommand{\eex}{\end{eqnarray*}}
\newcommand{\ba}{\begin{array}}
\newcommand{\ea}{\end{array}}

\newcommand{\0}{\boldsymbol{0}}
\definecolor{Red}{rgb}{1,0,0}
\definecolor{Blue}{rgb}{0,0,1}
\definecolor{Green}{rgb}{0,1,0}
\definecolor{magenta}{rgb}{1,0,0.6}
\definecolor{lightblue}{rgb}{0,0.5,1}
\definecolor{lightpurple}{rgb}{0.6,0.4,1}
\definecolor{gold}{rgb}{0.6,0.5,0}
\definecolor{orange}{rgb}{1,0.4,0}
\definecolor{hotpink}{rgb}{1,0,0.5}
\definecolor{newcolor2}{rgb}{0.5,0.3,0.5}
\definecolor{newcolor}{rgb}{0,0.3,1}
\definecolor{newcolor3}{rgb}{1,0,0.35}
\definecolor{darkgreen1}{rgb}{0,0.35, 0}
\definecolor{darkgreen}{rgb}{0,0.6, 0}
\definecolor{darkred}{rgb}{0.75,0,0}
\font\tenbi=cmmib10   at 11 pt
\font\sevenbi=cmmib10 at 9pt
\font\fivebi=cmmib7 at 6pt
\newfam\bifam
\textfont\bifam=\tenbi \scriptfont\bifam=\sevenbi
\scriptscriptfont\bifam=\fivebi

\font\tendb=msbm10 at 12 pt
\font\sevendb=msbm7
\newfam\dbfam
\textfont\dbfam=\tendb \scriptfont\dbfam=\sevendb

\renewcommand{\arraystretch}{1}

\def\0{{\bm0}}

 \makeatletter
\def\ps@pprintTitle{%
 \let\@oddhead\@empty
 \let\@evenhead\@empty
 \def\@oddfoot{\centerline{\thepage}}%
 \let\@evenfoot\@oddfoot}
\makeatother
\usepackage[mathlines]{lineno}
\title{Fractional Gray-Scott Model: Well-posedness, Discretization, and Simulations\footnotemark[1]}
\author{Tingting Wang\footnotemark[2] \footnotemark[3]
\and Fangying Song\footnotemark[3]
\and Hong Wang\footnotemark[4]
\and George Em Karniadakis\footnotemark[3]
}

\begin{document}
%\widowpenalty-10
\maketitle

 \footnotetext[1]{This work was supported by the OSD/ARO/MURI on ``Fractional PDEs for Conservation Laws and Beyond: Theory, Numerics and Applications (W911NF-15-1-0562)" and the National Science Foundation under Grant DMS-1620194. The first author was supported by the China Scholarship Council under 201706220157.}
\footnotetext[2]{School of Mathematics, Shandong University, Jinan 250100, Shandong, China.}
 
\footnotetext[3]{Division of Applied Mathematics, Brown University, Providence, RI 02912, USA
(george$\_$karniadakis@brown.edu, fangying$\_$song@brown.edu,  tingting$\_$wang@brown.edu).}
 
\footnotetext[4]{Department of Mathematics, University of South Carolina, Columbia, SC 29208, USA
(hwang@math.sc.edu).}

%%%

\begin{abstract}
The Gray-Scott (GS) model represents the dynamics and steady state pattern formation in reaction-diffusion systems and has been extensively studied in the past. In this paper, we consider the effects of anomalous diffusion on pattern formation by introducing the fractional Laplacian into the GS model. First, we prove that the continuous solutions of the fractional GS model are unique. We then introduce the Crank-Nicolson (C-N) scheme for time discretization and weighted shifted Gr\"unwald difference operator for spatial discretization. We perform stability analysis for the time semi-discrete numerical scheme, and furthermore, we analyze numerically the errors with benchmark solutions that show second-order convergence both in time and space. We also employ the spectral collocation method in space and C-N scheme in time to solve the GS model in order to verify the accuracy of our numerical solutions. We observe the formation of different patterns at different values of the fractional order, which are quite different than the patterns of the corresponding integer-order GS model, and quantify them by using the radial distribution function (RDF). Finally, we discover the scaling law for steady patterns of the RDFs in terms of the fractional order $1<\alpha \leq 2 $.
\end{abstract}


\section{Introduction}
\setcounter{equation}{0}
In the past several decades, the formation of spatial and temporal patterns has become a very active area of research. There are many diverse patterns formed in physical, biological and chemical systems \cite{kerr2008fast,nicolis1977self,vastano1987chemical}. Among various systems, the reaction and diffusion systems attract much attention, since they create a variety of patterns that could be found in nature, for instance, spots, spot replication, stripes, and travelling waves (e.g. \cite{chen2011stability,pearson1993complex}). A representative reaction and diffusion model is the GS model, which is a variant of the autocatalytic Selkov model of glycolysis  \cite{gray1983autocatalytic,sel1968self}. This model includes the following two reactions
\begin{equation}\label{si:e1}
\begin{aligned}
U+2&V \to 3V,\\
&V\to P,
\end{aligned}
\end{equation}
where $U,V$ and $P$ represent the chemical species. The two reactions take place in an open flow reactor, in which $U$ is continuously supplied and the final product $P$ is removed {\cite{mazin1996pattern}}. The first reaction shows a process of autocatalysis and the second reaction describes the decay of $V$ into $P$. In addition, note that there is a non-equilibrium constraint on $U$ by constantly feeding it and removing $P$, which leads to a variety of unstable phenomena. 

The classical (integer-order) GS model is expressed in the form 
\begin{equation}\label{si:e2}\begin{aligned}
&\frac{\partial u}{\partial t}=\mu_{u}\Delta u-uv^2+F(1-u),\\
&\frac{\partial v}{\partial t}=\mu_{v}\Delta v+uv^2-(F+\kappa)v,
\end{aligned}\end{equation}
where $u$ and $v$ are the concentrations of the two chemical components. $F$ is the feed rate and $\kappa$ is the decay rate of the second reaction; $\mu_{u}$ and $\mu_{v}$ are the diffusion coefficients. 

The response of one dimensional GS model was studied by Vastano and co-works previously \cite{vastano1987chemical}. Doelman et al. \cite{doelman1998stability, doelman1997pattern} investigated the asymptotic scaling of parameters and variables necessary for the analysis of the patterns in one dimension. Pearson \cite{pearson1993complex} studied this system in two dimensions and presented very complicated spatio-temporal patterns. Pearson also did a thorough numerical study for this system. Many complex structures were observed in the numerical solutions. Hale studied the exact homoclinic and heteroclinic solutions of the GS model for autocatalysis \cite{hale2000exact}. McGough and Riley produced the bifurcation analysis to support the existing numerical evidence for patterns and derived the bifurcation results for nonuniform steady states \cite{mcgough2004pattern}.  

\begin{figure}[htp]
\centering
\subfloat[$\alpha=2.0.$]{
\begin{minipage}{.30\textwidth}\centering
\includegraphics[width=1.0\textwidth]{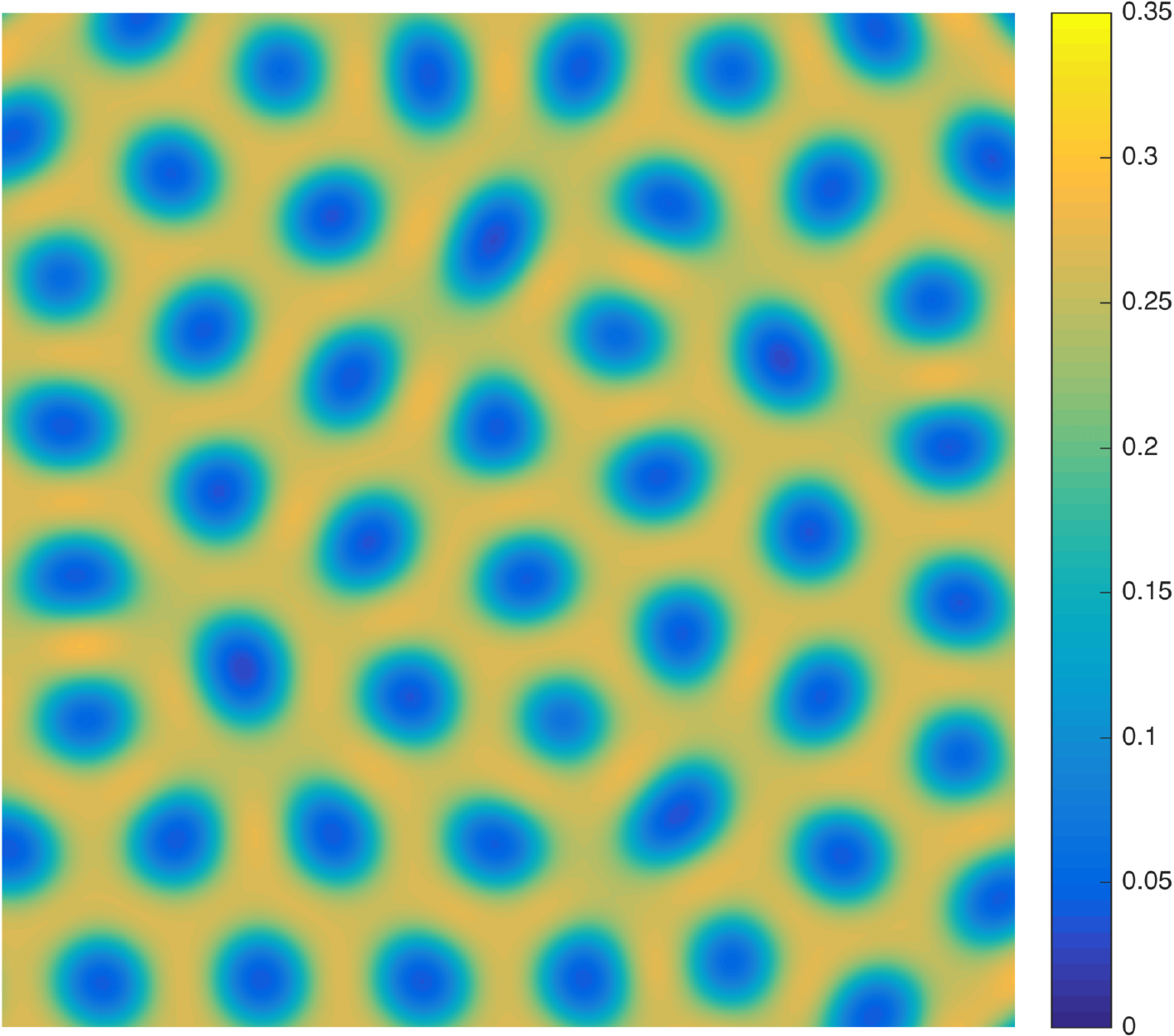}\\
\includegraphics[width=1.0\textwidth]{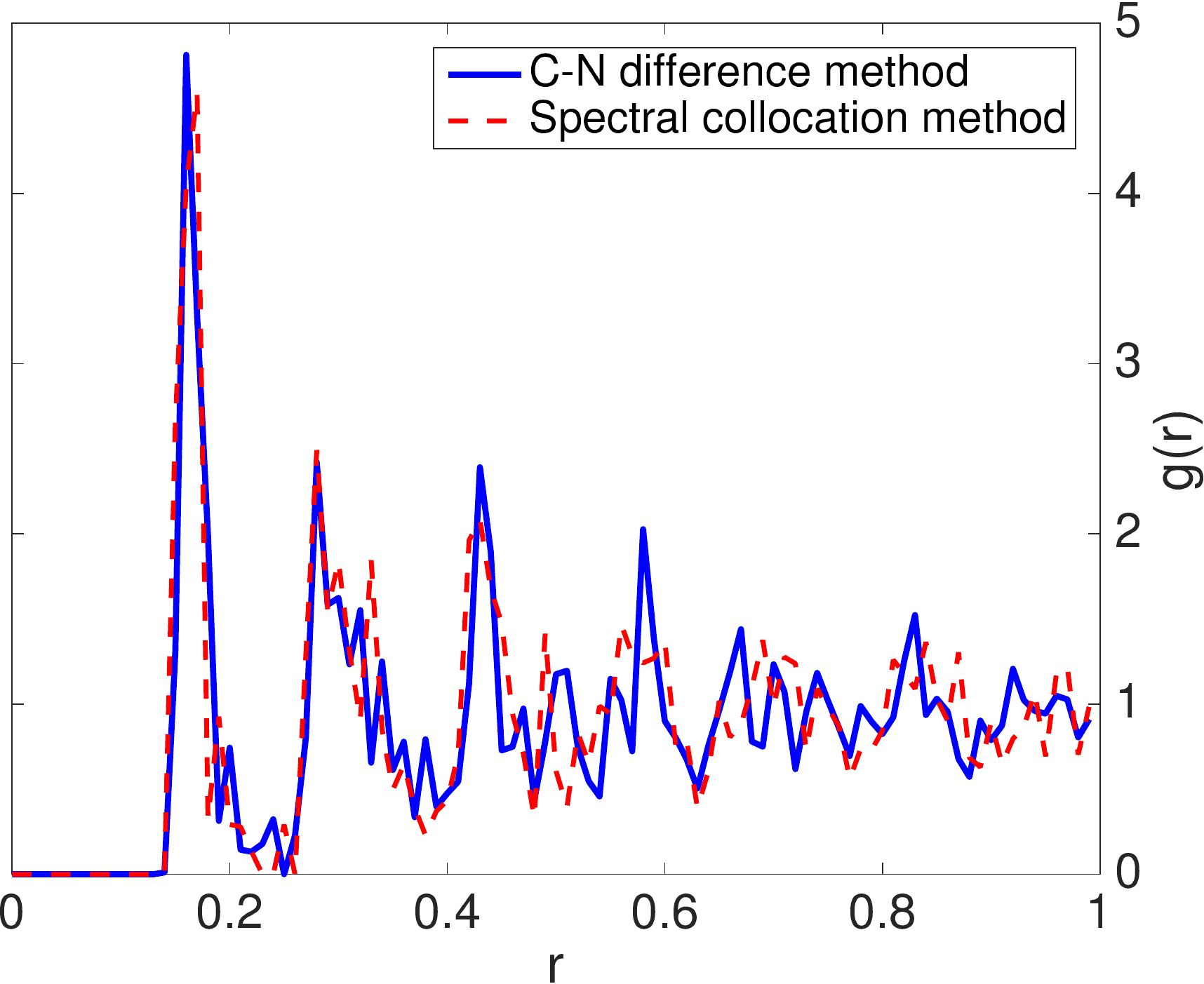}
\end{minipage}}
\subfloat[$\alpha=1.5.$]{
\begin{minipage}{.30\textwidth}\centering
\includegraphics[width=1.0\textwidth]{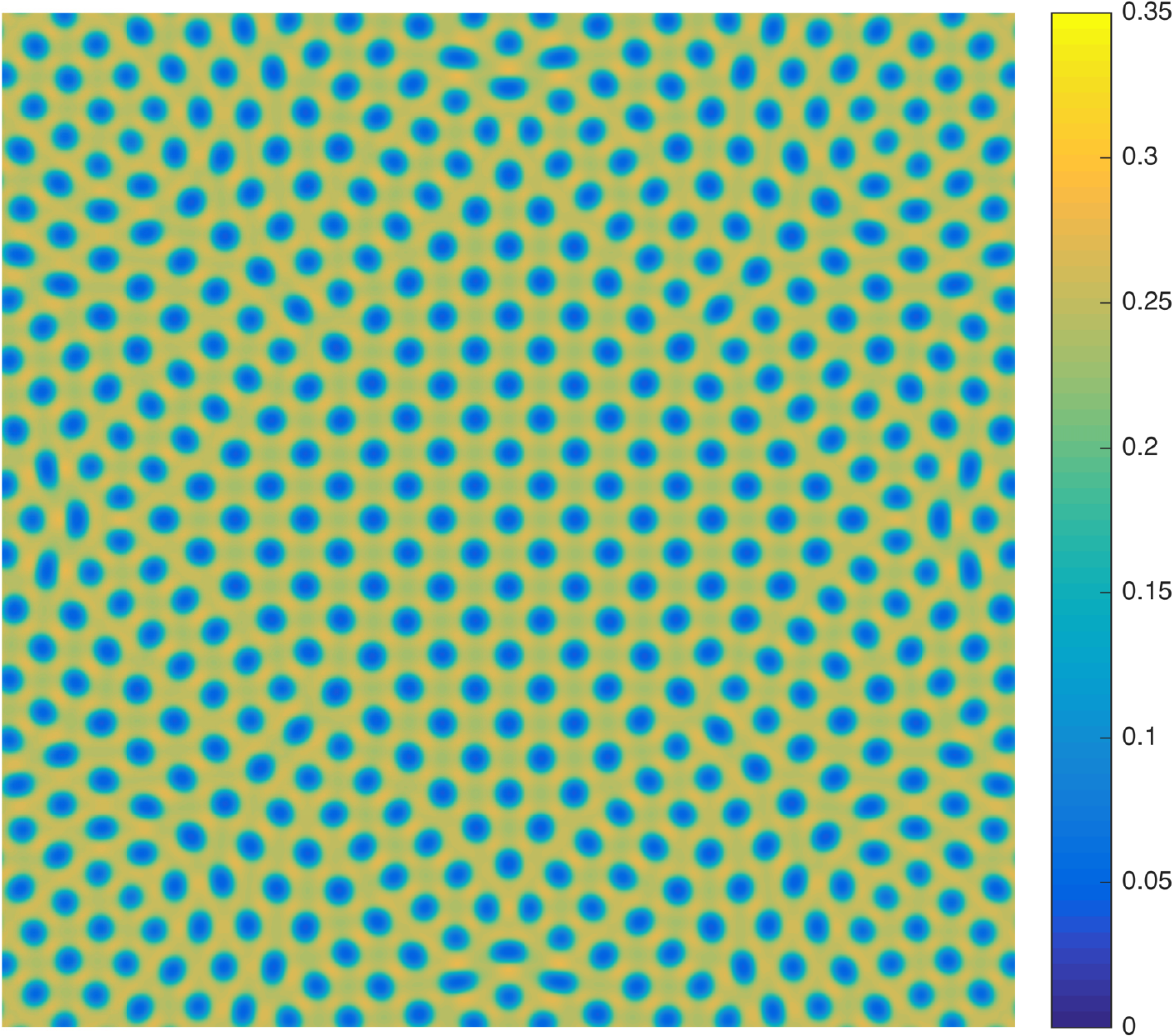}\\
\includegraphics[width=1.0\textwidth]{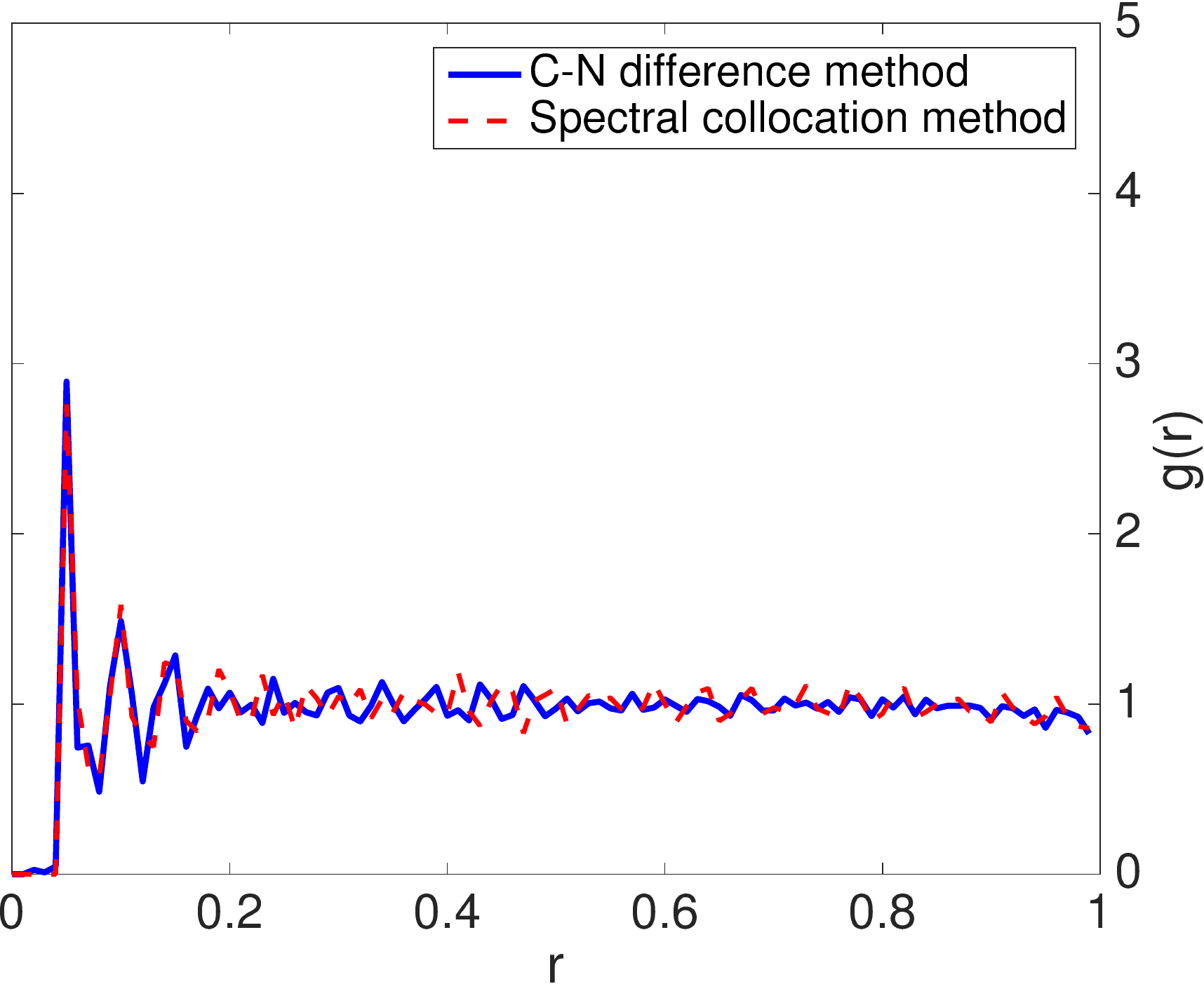}
\end{minipage}}
\caption{The difference in pattern formation and the corresponding RDFs between normal diffusion and anomalous diffusion described by the GS model for $\alpha=2.0, 1.5$.\label{Fig. 0}}
\end{figure}
Note that, in the aforementioned literature on pattern dynamics of the GS model, the models are all with standard diffusion, i.e. the diffusion operator is the normal Laplacian. In addition, different pattern formation was observed by changing the values of parameters $F$ and $\kappa$. Hence, here we aim to study the effects of the super-diffusion (with fractional Laplacian $(-\Delta)^{\frac{\alpha}{2}}$ for $1<\alpha \leq 2$) on pattern formation of this model and perform some theoretical analysis for the GS model. Fig. \ref{Fig. 0} shows the different pattern formation and the corresponding radial distribution functions (RDFs) between normal diffusion and anomalous diffusion for the GS model. The standard diffusion systems correspond to Brownian motion, while the fractional diffusion systems correspond to L\'{e}vy process. The relationship between L\'{e}vy process and certain types of space fractional models has been established in \cite{metzler2000random}. In the paper by Wu et al. \cite{liu2017turing}, the Turing instability and pattern formation of the Lengyel-Epstein model with super-diffusion were studied. In that paper, it was emphasized that more complex dynamics will appear under the super-diffusion. Bueno-Orovio and co-workers used the Fourier spectral method to solve several types of fractional reaction diffusion equations in their paper \cite{bueno2014fourier}. Lee \cite{lee2017second} also introduced a second-order operator splitting Fourier spectral method to approximate the fractional-in-space reaction diffusion equations. All these papers described investigations of the fractional reaction-diffusion models. So in this paper, we will consider the fractional GS model with space fractional Laplacian (super-diffusion effects) in two-dimensions. 

Since the exact solutions for the fractional GS model are not known, we will develop second-order numerical methods to solve the GS equations. The standard approach for solving space fractional diffusion equations is to use finite difference method, finite element method to discrete fractional derivatives and then use Euler formulation for the evolution of time \cite{deng2008finite, ilic2005numerical, liu2007stability, meerschaert2006finite}. However, these methods require the solution of a linear system at each time step, which correspond to a large, dense matrix due to the nonlocal nature of the fractional operator. Roop et al. \cite{ervin2007numerical} analyzed a fully discrete finite element approximation to a time dependent fractional order diffusion equation, which contains a nonlocal quadratic nonlinearity. In the paper by Meerschaert \cite{meerschaert2006finite, tadjeran2006second}, the approach based on the C-N method combined with spatial extrapolation was used to derive temporally and spatially second-order accurate numerical estimates. Deng et al. \cite{tian2015class} introduced the weighted shifted Gr\"{u}nwald difference operator to approximate the Riemann-Liouville fractional derivative and obtained a second-order accuracy in space. To overcome the bottleneck of expensive computations, Wang et al. introduced a fast numerical algorithm, which tackled the problems for two or three dimensions effectively \cite{wang2014fast,wang2010direct}. In addition to these classical numerical methods, spectral methods also have been used for space fractional equations \cite{bueno2014fourier, lee2017second}.  

In this paper we shall use the finite difference method to derive a second-order numerical scheme for the fractional GS model. In this fractional model, the classical Laplacian operator is replaced by the fractional Laplacian operator $(-\Delta)^{\alpha/2}$ with $1< \alpha \leq 2$. We apply the weighted shifted Gr\"unwald difference discretization method, which leads to well structured, relatively sparse and positive definite coefficient matrix \cite{tian2015class}. For time direction, the C-N scheme is employed to obtain a temporally second-order estimate. Since this fractional GS model is a nonlinear system, we use the second-order implicit-explicit methods \cite{ascher1995implicit} to handle this problem, i.e. an implicit scheme is used for the linear terms and an explicit scheme is used for the nonlinear terms. Moreover, we carry out the linear stability analysis for the steady states of the GS model. We also derive the well-posedness of the fractional GS model. In addition, we provide the stability analysis for the time semi-discrete scheme. Several numerical experiments have been conducted to verify the accuracy in time and space of this numerical scheme. In the simulations of the fractional GS model, a small perturbation has been added to the initial states. We observe the formation of patterns under the condition of different parameter values. Moreover, the spectral collocation method is used in space discretization to simulate this fractional model. We compare the steady patterns obtained by the two different numerical methods and calculate the corresponding RDFs. Finally, we estimate the scaling law between the fractional orders and the distance between all spot pairs in the steady spot patterns. 

This paper is organized as follows. In section 2, we introduce the fractional GS model and perform stability analysis for steady states. The well-posedness of this model is presented in section 3. In section 4, we propose a second-order accurate both in time and space numerical scheme for the discretization of the fractional GS model. We also provide the stability analysis for the time semi-discrete scheme. In section 5, we present numerical simulations of the fractional GS model, including convergence results and the scaling law for steady patterns of the RDFs. We conclude in section 6 with a summary.

%%%%%%%%%%%%%%%%%%%%
\section{The fractional GS model}
\setcounter{equation}{0}
The fractional GS model that describes an autocatalytic reaction-diffusion process between two chemical species with concentrations $u$ and $v$ is written as:
\begin{equation}\label{s1:e1}\begin{cases}
&\frac{\partial u}{\partial t}=-\mu_{u}(-\Delta )^{\frac{\alpha}{2}}u-uv^{2}+F(1-u), \\
&\frac{\partial v}{\partial t}=-\mu_{v}(-\Delta )^{\frac{\alpha}{2}}v+uv^{2}-(F+\kappa)v,~~~~ \qquad \qquad(x,y,t)\in \Omega \times [0,T],\\
&u(x,y,0)=u_{0}(x,y),\quad v(x,y,0)=v_{0}(x,y),\quad   \qquad (x,y)\in \Omega, \\
&u(x,y,t)=0, \quad v(x,y,t)=0, \qquad \qquad \qquad \qquad   ~~(x,y,t)\in \partial \Omega \times [0,T],
\end{cases}\end{equation}
where $\Omega=(a,b)\times(c,d), 1 < \alpha \le 2$. The diffusion coefficients satisfy $\mu_{u} \ge 0, \mu_{v} \ge 0$. The parameters $F, \kappa$ are positive constants representing feed rate and decay rate, respectively. In this paper we define the fractional Laplacian operator by Riesz fractional derivatives as follows
\begin{equation}\label{s2:ee1}
\begin{aligned}
-(-\Delta)^{\frac{\alpha}{2}}u= \frac{\partial^{\alpha}u}{\partial|x|^{\alpha}}
+\frac{\partial^{\alpha}u}{\partial|y|^{\alpha}} 
=&-\frac{1}{2\cos(\frac{\pi \alpha}{2})} \Bigl({}_{a}D_{x}^{\alpha} u+{}_{x}D_{b}^{\alpha}u \Bigr) \qquad \qquad \qquad\\
&-\frac{1}{2\cos(\frac{\pi \alpha}{2})} \Bigl({}_{c}D_{y}^{\alpha} u+{}_{y}D_{d}^{\alpha}u \Bigr),
\end{aligned}
\end{equation}
with ${}_{a}D_{x}^{\alpha}, {}_{x}D_{b}^{\alpha}$ and ${}_{c}D_{y}^{\alpha}, {}_{y}D_{d}^{\alpha}$ being the Riemann-Liouville fractional operators.

The fractional GS system is an activation-substrate depletion system \cite{mazin1996pattern}. The chemical specie $V$ grows auto-catalytically on the the specie $U$ i.e., the continuously fed substrate. For instance, the concentrations $u, v$ vary opposite to each other. In other words, the existence of $V$ in the second equation of (\ref{s1:e1}) will prompt the production of $V$ and reduce the concentration of $U$. In practical experiments, the formation of patterns is affected by the parameter values and anomalous diffusion.

\subsection{Steady states and linear stability analysis}\label{2.1}

In this subsection, we consider the spatially uniform steady states of the fractional GS model with $(-\Delta)^{\frac{\alpha}{2}}u\equiv 0$:
\begin{equation}\label{ss:e1}
\begin{aligned}
0&= -uv^2+F(1-u),\\
0 &= uv^2-(F+\kappa)v.
\end{aligned}
\end{equation}
This system has an unique trivial steady state $(u_{*}, v_{*})=(1, 0)$ for all the values of $F$ and $\kappa$. In addition, there also exists two  steady states $(u_{+}, v_{-})$ and $(u_{-}, v_{+})$ when $F\geq 4(F+\kappa)^2 $. Namely, we have
\begin{equation}\label{ss:e2}
\begin{aligned}
u_{\pm}&=\frac{1}{2}(1\pm \sqrt{1-4\gamma^2 F}),\\
v_{\mp}&=\frac{1}{2\gamma}(1\mp \sqrt{1-4\gamma^2 F}), \quad \gamma=\frac{F+\kappa}{F}.
\end{aligned}
\end{equation}
As a result, we can get the saddle-node bifurcation
\begin{equation}\label{ss:e3}
\kappa_{c}=-F+\frac{1}{2}\sqrt{F}, \quad 0\leq F \leq \frac{1}{4}.
\end{equation}

In order to analyze the stability of these spatially uniform steady states, we introduce perturbations to the fractional system by $\epsilon \tilde{u}$ and $\epsilon \tilde{v}$ and obtain 
\begin{equation}\label{ss:e4}
\begin{aligned}
\frac{\partial \tilde{u}}{\partial t}&=-\mu_{u}(-\Delta)^{\frac{\alpha}{2}}\tilde{u}-\tilde{u}(v^2+F)-2uv\tilde{v}+O(\epsilon),\\
\frac{\partial \tilde{v}}{\partial t}&=-\mu_{v}(-\Delta)^{\frac{\alpha}{2}}\tilde{v}-(F+\kappa-2uv)\tilde{v}+v^2\tilde{u}+O(\epsilon).
\end{aligned}
\end{equation}
By neglecting terms $O(\epsilon)$, we can get the normal mode solution 
\begin{equation}\label{ss:e5}
\tilde{u}=u_{o}e^{\lambda t-i(k_{1}x+k_{2}y)},\quad \tilde{v}=v_{o}e^{\lambda t-i(k_{1}x+k_{2}y)},
\end{equation}
with amplitudes $u_{o}, v_{o}$ and wave number $k_{1}, k_{2}$, $|\boldsymbol{k}|^{\alpha}=|k_{1}|^{\alpha}+|k_{2}|^{\alpha}$.

We first investigate the eigenvalues of their characteristic equations to examine the stability of the steady states. By substituting $(\ref{ss:e5})$ into $(\ref{ss:e4})$, the corresponding characteristic equation is derived as follows
\begin{equation}\label{ss:e6}
\left |\begin{array}{cc}
\lambda+\mu_{u}|\boldsymbol{k}|^{\alpha}+v^2+F &2uv\\
-v^2 &\lambda+\mu_{v}|\boldsymbol{k}|^{\alpha}+F+\kappa-2uv
\end{array}\right |=0.
\end{equation}
Then we obtain the following dispersion relation
\begin{equation}\label{ss:ee6}
\lambda^2+\lambda T_{k}+D_{k}=0,
\end{equation}
where
\begin{equation*}\begin{aligned}
T_{k}=&(\mu_{u}+\mu_{v})|\boldsymbol{k}|^{\alpha}+2F+\kappa+v^2-2uv,\\
D_{k}=&\mu_{u}\mu_{v}(|\boldsymbol{k}|^{\alpha})^2+\Bigl[(v^2+F)\mu_{v}+(F+\kappa-2uv)\mu_{u}\Bigr]|\boldsymbol{k}|^{\alpha}\\
&+(F+\kappa)(v^2+F)-2Fuv.
\end{aligned}\end{equation*}
From $(\ref{ss:ee6})$, we get the eigenvalues corresponding to the trivial steady state,
\begin{equation}\label{ss:e7}
\lambda_{1}=-\mu_{u}|\boldsymbol{k}|^{\alpha}-F,\quad \quad \lambda_{2}=-\mu_{v}|\boldsymbol{k}|^{\alpha}-F-\kappa,
\end{equation}
which depend on the fractional order $\alpha$. These are extensive models of the integral case \cite{mcgough2004pattern}. For $\lambda_{1}, \lambda_{2}<0$, it is obvious that this trivial steady state is stable for all values of $F$ and $\kappa$. In terms of the other two alternative steady states, first we insert the point $(u_{+}, v_{-})$ given by (\ref{ss:e2}) to (\ref{ss:ee6}). The computation shows that $v_{-}^{2}\leq F$ for $k_{1},k_{2}=0$. Since the steady states satisfy the equations in (\ref{ss:e1}), making use of $uv=F+\kappa$, we find that $T_{k}<0$ and $D_{k}<0$ for $k_{1},k_{2}=0$. Therefore, the steady state $(u_{+}, v_{-})$ is always unstable. However, the steady state $(u_{-}, v_{+})$ may have stable and unstable structure. 

We consider the case that the diffusion constants are negligibly small with ratio of order 1. In this case, when $\lambda$ in (\ref{ss:ee6}) is purely imaginary, the system undergoes a Hopf bifurcation which closes to the lower branch of the saddle-node curve (see \cite{mazin1996pattern}). When $T_{k}^2-4D_{k}<0$, the eigenvalues of equation (\ref{ss:ee6}) are complex conjugate. Then the Hopf bifurcation can be obtained under the condition that $T_{k}=0, D_{k}>0$ for $k_{1},k_{2}=0$. Inserting the steady state $(u_{-}, v_{+})$, the critical feed rate is given by
\begin{equation}\label{ss:e9}
F_{c}=\frac{\sqrt{\kappa}-2\kappa-\sqrt{(2\kappa-\sqrt{\kappa})^2-4\kappa^{2}}}{2}, \quad 0\leq \kappa \leq \kappa_{c}.
\end{equation}
\begin{figure}[htpb]
\begin{center}
\includegraphics[width=3.8in,height=2.75in]{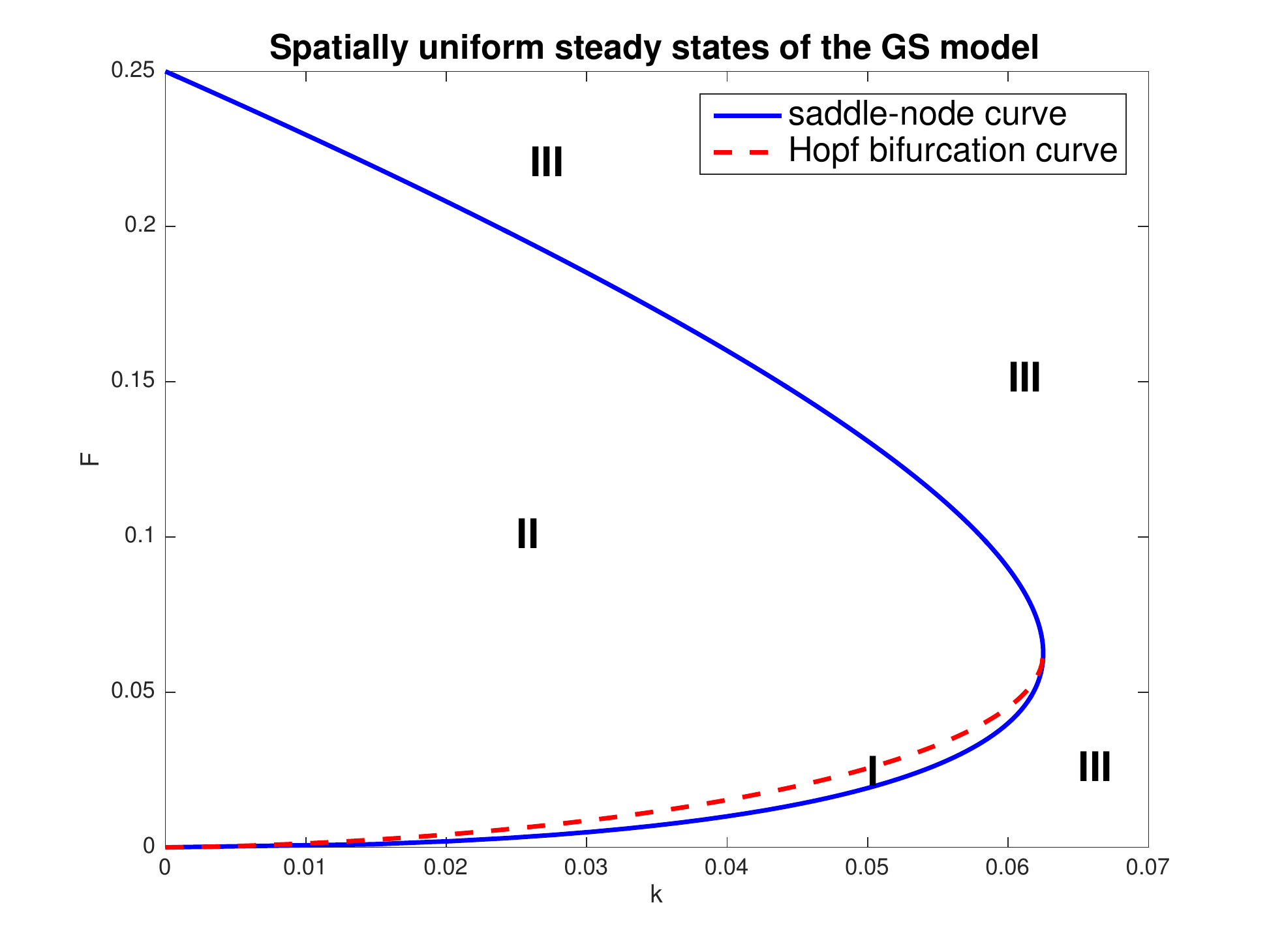} 
\caption{Phase diagram of the inviscid dynamics. In region III, there is only one spatially uniform steady state $(u_{*}=1, v_{*}=0)$, which is stable for all $F$ and $\kappa$. In region I and II , there are three steady states. In region II, only two steady states $(1, 0)$ and $(u_{-}, v_{+})$ are stable. However, the steady state $(1,0)$ is stable in region I, while the uniform state $(u_{-}, v_{+})$ loses stability when $F$ is decreased through the Hopf bifurcation curve (the red dashed-line). The nontrivial fixed point $(u_{+}, v_{-})$ is always unstable. The bifurcating periodic orbit is stable for $\kappa <0.035$ and unstable for $\kappa >0.035$.\label{Fig. 1}}
%\label{default}
\end{center}
\end{figure}

In the $(\kappa, F)$ plane, we plot the phase diagram for the homogeneous fractional GS system. In Fig. $\ref{Fig. 1}$, outside the region bounded by the saddle-node curve (region III), there only exists a single spatially uniform state (the trivial state $(u_{*}=1, v_{*}=0)$) that is stable for all $(\kappa, F)$, while in the remainder region (region I and II) of this plane, there are three spatially uniform steady states. In region II, the system is bistable with the steady states $(1, 0)$ and $(u_{-}, v_{+})$. Since the coefficient $T_{k}$ is less than zero in region I, the steady state $(u_{-}, v_{+})$ loses stability when $F$ is decreased through the Hopf bifurcation curve (the red dashed-line). The third steady state $(u_{+}, v_{-})$ is always unstable for arbitrary parameter. In addition, the intersection point of the two curves is $(\kappa_{c}, F_{c})=(\frac{1}{16}, \frac{1}{16})$.

%%%%%%%%%%%%%%%%%%%%%%%%%%%%%%%%%

\section{Well-posedness}
\setcounter{equation}{0}

The well-posedness of the GS model with classical diffusion has  been analyzed in \cite{you2008global}. In this section, we consider the well-posedness of the fractional GS model. Firstly, we give some useful property and lemmas (see \cite{ervin2006variational,li2009space,podlubny1998fractional}).
\begin{Property}\label{s3:p1}
If $0<p<1,~0<q<1,~v(0)=0,~t>0$, then
\begin{equation*}
{}_{0}D^{p+q}_{t}v(t)={}_{0}D^{p}_{t}{}_{0}D^{q}_{t}v(t)={}_{0}D^{q}_{t}{}_{0}D^{p}_{t}v(t).
\end{equation*} 
\end{Property}
\begin{lemma}\label{s3:l1}
For real $s$, $0<s<1$, if $w(t)\in H^{s}(I), ~v(t)\in C^{\infty}_{0}(I)$, then
\begin{equation*}
({}_{0}D^{s}_{t}w(t),~v(t))_{I}=(w(t),~{}_{t}D^{s}_{T}v(t))_{I}.
\end{equation*}
\end{lemma}
\begin{lemma}\label{s3:l2}
For real $s>0, ~v\in C_{0}^{\infty}(\mathbb{R})$, then
\begin{equation*}\begin{aligned}
&({}_{-\infty}D^{s}_{t}v(t),~{}_{t}D^{s}_{\infty}v(t))=\cos(\pi s)\|{}_{-\infty}D^{s}_{t}v(t)\|_{L^{2}(\mathbb{R})}^{2},\\
&({}_{-\infty}D^{s}_{t}v(t),~{}_{t}D^{s}_{\infty}v(t))=\cos(\pi s)\|{}_{t}D^{s}_{\infty}v(t)\|_{L^{2}(\mathbb{R})}^{2}.
\end{aligned}
\end{equation*}
\end{lemma}

We can extend the fractional derivatives ${}_{a}D^{s}_{x}v, {}_{x}D^{s}_{b}v$ to ${}_{-\infty}D^{s}_{x}v, {}_{x}D^{s}_{\infty}v$ by zero outside $(a,b)$ easily. Using the above property, lemmas and the definition of the fractional Laplacian (\ref{s2:ee1}), we derive the following inner product formula i.e., for $1<\alpha \leq 2,~u\in H^{\frac{\alpha}{2}}_{0}(\Omega)$,
\begin{equation}\label{s3:ee2}\begin{aligned}
 ((-\Delta)^{\frac{\alpha}{2}}u,~u)&=\frac{1}{2\cos(\frac{\pi \alpha}{2})}\Bigl(({}_{a}D_{x}^{\alpha}u,~u)+({}_{x}D_{b}^{\alpha}u,~u)+({}_{c}D_{y}^{\alpha}u,~u)+({}_{y}D_{d}^{\alpha}u,~u)\Bigr)\\
&= \frac{1}{2\cos(\frac{\pi \alpha}{2})}\Bigl(2({}_{a}D_{x}^{\frac{\alpha}{2}}u,~{}_{x}D_{b}^{\frac{\alpha}{2}}u)+2({}_{c}D_{y}^{\frac{\alpha}{2}}u,~{}_{y}D_{d}^{\frac{\alpha}{2}}u)\Bigr)\\
&=\|{}_{a}D^{\frac{\alpha}{2}}_{x}u\|^{2}+\|{}_{c}D^{\frac{\alpha}{2}}_{y}u\|^{2}=:\|D^{\frac{\alpha}{2}}u\|^2,
\end{aligned}\end{equation}
where $u \in H^{\frac{\alpha}{2}}_{0}(\Omega)$ is the completely space of $C^{\infty}_{0}(\Omega)$ in $H^{\frac{\alpha}{2}}(\Omega)$. We omit the subscripts of the $L^2$ norm in this paper whenever it is clear from the context. Then the solution of the fractional GS model can be bounded as in the following theorem.

\begin{theorem}\label{s3:t1}
Suppose the solutions $u(t),v(t)$ of the fractional GS model belong to $H^{\frac{\alpha}{2}}_{0}(\Omega)$ for any initial data $u_{0}, v_{0}\in H^{\frac{\alpha}{2}}_{0}(\Omega)$. Then the following estimates hold
\begin{equation}\label{tt:1}
\|u(t)\|^{2} \leq e^{-Ft}\|u_{0}\|^{2}+|\Omega|,~~~~ t \in [0,T],  
\end{equation}
\begin{equation}\label{tt:2}
\|w(t)\|^{2} \leq \|u_{0}+v_{0}\|^{2}+\Bigl(\frac{\kappa}{F}+\frac{|\mu_{v}-\mu_{u}|^{2}}{2\mu_{u}\mu_{v}|\cos(\frac{\pi \alpha}{2})|^2}\Bigr)\|u_{0}\|^{2}
+\Bigl(F+\kappa+\frac{|\mu_{v}-\mu_{u}|^{2}F}{2\mu_{u}\mu_{v}|\cos(\frac{\pi \alpha}{2})|^2}\Bigr)|\Omega|t, \\[0.05in]
\end{equation}
where $w(t)=u(t)+v(t)$ is the solution of the following equation
\begin{equation}\label{tt:3}
w_{t}=-\mu_{v}(-\Delta)^{\frac{\alpha}{2}}w-(F+\kappa)w+[(\mu_{v}-\mu_{u})(-\Delta)^{\frac{\alpha}{2}}u+\kappa u+F].
\end{equation}
\end{theorem}

The proof is presented in Appendix \ref{adx1}.

%%%%%%%%%%%%%%%%%%%%%%%%%%%%%%%%%
\section{Numerical discretization and stability analysis}\label{s1}
\setcounter{equation}{0}
In this section, we use the C-N difference scheme for time discretization and use the weighted shifted Gr\"unwald difference operator introduced in \cite{tian2015class} to approximate the spatial fractional Laplacian operator. In addition, we apply the second-order implicit-explicit method to handle the nonlinear terms.
\subsection{Numerical discretization}
Let $N_{x},N_{y}$ and $M$ be positive integers. We define the space steps and time step as $h_{x}=(b-a)/N_{x}, h_{y}=(d-c)/N_{y}$ and $\tau=T/M$. Then we partition the space domain $\Omega$ and time interval $[0,T]$ into the uniform mesh with $x_{i}=a+ih_{x}, y_{j}=c+jh_{y}$ and $t_{n}=n\tau$ for $1\leq i\leq N_{x}-1, 1\leq j\leq N_{y}-1$ and $0\leq n\leq M$. We let $t_{n+1/2}=(t_{n}+t_{n+1})/2$, for $0 \leq n \leq M-1$, and introduce the following notations: 
\begin{flalign}\label{s1:ee2}\begin{split}
&u_{i,j}^{n}=u(x_{i}, y_{j}, t_{n}), \quad v_{i,j}^{n}=v(x_{i}, y_{j}, t_{n}),\quad \quad \quad \quad \\
&u_{i,j}^{n+\frac{1}{2}}=(u_{i,j}^{n+1} +u_{i,j}^{n})/2, \quad v_{i,j}^{*, n+\frac{1}{2}}=(3v_{i,j}^{n} -v_{i,j}^{n-1})/2,\quad \quad \quad \quad\\
&\delta_{t}u_{i,j}^{n+\frac{1}{2}}=(u_{i,j}^{n+1} -u_{i,j}^{n})/\tau, \quad \delta_{t}v_{i,j}^{n+\frac{1}{2}}=(v_{i,j}^{n+1} -v_{i,j}^{n})/\tau.\quad \quad \quad 
\end{split}\end{flalign}
We consider the first equation of $u$ in problem $(\ref{s1:e1})$. In time direction, we use the Taylor expansion as the first step $n=0$,
\begin{equation}\label{ss:ee1}\begin{aligned}
\delta_{t}u_{i,j}^{\frac{1}{2}}=&-K_{u} \Bigl(({}_{a}D_{x}^{\alpha} u)_{i,j}^{1}+({}_{x}D_{b}^{\alpha}u)_{i,j}^{1}+({}_{c}D_{y}^{\alpha} u)_{i,j}^{1} +({}_{y}D_{d}^{\alpha} u)_{i,j}^{1}\Bigr) \\
&-K_{u} \Bigl(({}_{a}D_{x}^{\alpha} u)_{i,j}^{0}+({}_{x}D_{b}^{\alpha}u)_{i,j}^{0}+({}_{c}D_{y}^{\alpha} u)_{i,j}^{0} +({}_{y}D_{d}^{\alpha} u)_{i,j}^{0}\Bigr) \\
&-u_{i,j}^{\frac{1}{2}}(v_{i,j}^{0})^2+F(1-u_{i,j}^{\frac{1}{2}})+O(\tau),
\end{aligned}
\end{equation}
and for $1 \leq n \leq M-1$,
\begin{equation}\label{s1:e3}
\begin{aligned}
\delta_{t}u_{i,j}^{n+\frac{1}{2}}=&-K_{u} \Bigl(({}_{a}D_{x}^{\alpha} u)_{i,j}^{n+1}+({}_{x}D_{b}^{\alpha}u)_{i,j}^{n+1}+({}_{c}D_{y}^{\alpha} u)_{i,j}^{n+1} +({}_{y}D_{d}^{\alpha} u)_{i,j}^{n+1}\Bigr) \\
&-K_{u} \Bigl(({}_{a}D_{x}^{\alpha} u)_{i,j}^{n}+({}_{x}D_{b}^{\alpha}u)_{i,j}^{n}+({}_{c}D_{y}^{\alpha} u)_{i,j}^{n} +({}_{y}D_{d}^{\alpha} u)_{i,j}^{n}\Bigr) \\
&-u_{i,j}^{n+\frac{1}{2}}(v_{i,j}^{*, n+\frac{1}{2}})^2+F(1-u_{i,j}^{n+\frac{1}{2}})+O(\tau^{2}),
\end{aligned}
\end{equation}
where $K_{u}=\frac{\mu_{u}}{4\cos(\frac{\pi \alpha}{2})}$.

In space, we use the second-order Gr\"unwald difference operators  ${}_{L}D_{h_{x}}^{\alpha}u, {}_{R}D_{h_{x}}^{\alpha}u$ and ${}_{L}D_{h_{y}}^{\alpha}u, {}_{R}D_{h_{y}}^{\alpha}u$ to approximate the fractional diffusion operators ${}_{a}D_{x}^{\alpha}u, {}_{x}D_{b}^{\alpha}u$ and ${}_{c}D_{y}^{\alpha}u, {}_{y}D_{d}^{\alpha}u$, i.e.,
\begin{equation*}
\begin{aligned}
&({}_{a}D_{x}^{\alpha} +{}_{x}D_{b}^{\alpha})u_{i,j}^{n}
=({}_{L}D_{h_{x}}^{\alpha}+{}_{R}D_{h_{x}}^{\alpha})u_{i,j}^{n}+O(h_{x}^{2}),\\
&({}_{c}D_{y}^{\alpha} +{}_{y}D_{d}^{\alpha})u_{i,j}^{n}
=({}_{L}D_{h_{y}}^{\alpha}+{}_{R}D_{h_{y}}^{\alpha})u_{i,j}^{n}+O(h_{y}^{2}),
\end{aligned}
\end{equation*}
where the above Gr\"unwald difference operators are derived in \cite{tian2015class} with $(p,q)=(1,0)$,
\begin{equation*}
\begin{aligned}
{}_{L}D_{h_{x}}^{\alpha}u(x_{i})& =\frac{1}{h_{x}^{\alpha}}\sum_{k=0}^{i+1}\omega^{(\alpha)}_{k}u(x_{i-k+1}),
~{}_{R}D_{h_{x}}^{\alpha}u(x_{i}) =\frac{1}{h_{x}^{\alpha}}\sum_{k=0}^{N_{x}-i+1}\omega^{(\alpha)}_{k}u(x_{i+k-1}),\\
{}_{L}D_{h_{y}}^{\alpha}u(y_{j}) &=\frac{1}{h_{y}^{\alpha}}\sum_{k=0}^{j+1}\omega^{(\alpha)}_{k}u(y_{j-k+1}),
~{}_{R}D_{h_{y}}^{\alpha}u(y_{j}) =\frac{1}{h_{y}^{\alpha}}\sum_{k=0}^{N_{y}-j+1}\omega^{(\alpha)}_{k}u(y_{j+k-1}).
\end{aligned}
\end{equation*}
The coefficients $\omega^{(\alpha)}_{k}$ are defined as follows
$$\omega_{0}^{(\alpha)}=\frac{\alpha}{2}g_{0}^{(\alpha)},~\omega_{k}^{(\alpha)}=
\frac{\alpha}{2}g_{k}^{(\alpha)}+\frac{2-\alpha}{2}g_{k-1}^{(\alpha)},~k \ge 1,$$
where $g_{k}^{\alpha}=(-1)^{k}\binom \alpha k$.

Then, multiplying equation (\ref{s1:e3}) with $\tau$, we obtain
\begin{equation}\label{s1:e4}
\begin{aligned}
&\Bigl(1+K_{u}\tau({}_{L}D_{h_{x}}^{\alpha}+{}_{R}D_{h_{x}}^{\alpha}
+{}_{L}D_{h_{y}}^{\alpha}+{}_{R}D_{h_{y}}^{\alpha})\Bigr)u_{i,j}^{n+1} \\
&\quad =\Bigl(1-K_{u}\tau({}_{L}D_{h_{x}}^{\alpha}+{}_{R}D_{h_{x}}^{\alpha}
+{}_{L}D_{h_{y}}^{\alpha}+{}_{R}D_{h_{y}}^{\alpha})\Bigr)u_{i,j}^{n} \\
&\quad \quad-\frac{\tau}{4}u_{i,j}^{n+\frac{1}{2}}(3v_{i,j}^{n}-v_{i,j}^{n-1})^{2}
+F\tau(1-u_{i,j}^{n+\frac{1}{2}})+\tau \epsilon_{i,j}^{n},
\end{aligned}
\end{equation}
where $|\epsilon_{i,j}^{n} |\leq c(\tau^{2}+h_{x}^{2}+h_{y}^{2})$ is the truncation error. We define
\begin{equation*}
\delta_{z}^{\alpha}={}_{L}D_{h_{z}}^{\alpha}+{}_{R}D_{h_{z}}^{\alpha},
~z=x,y.
\end{equation*}
Therefore, the equation $(\ref{s1:e4})$ can be rewritten as 
\begin{equation}\label{s1:e5}
\begin{aligned}
\Bigl(1+\tau K_{u}( \delta_{x}^{\alpha}+ \delta_{y}^{\alpha}) \Bigr)u_{i,j}^{n+1} 
 &=\Bigl(1-\tau K_{u}( \delta_{x}^{\alpha}+\delta_{y}^{\alpha})\Bigr)u_{i,j}^{n} 
 -\frac{\tau}{4}u_{i,j}^{n+\frac{1}{2}}(3v_{i,j}^{n}-v_{i,j}^{n-1})^{2}\\
&\quad +F\tau(1-u_{i,j}^{n+\frac{1}{2}})+\tau \epsilon_{i,j}^{n}.
\end{aligned}
\end{equation}
Using the Taylor expansion, we have
\begin{equation}\label{s1:e6}\begin{aligned}
(\tau K_{u})^{2}\delta_{x}^{\alpha}\delta_{y}^{\alpha}(u_{i,j}^{n+1}-u_{i,j}^{n})
&=\tau^{3}\Bigl((K_{u}~{}_{a}D_{x}^{\alpha}+K_{u}~{}_{x}D_{b}^{\alpha})
(K_{u}~{}_{c}D_{y}^{\alpha}+K_{u}~{}_{y}D_{d}^{\alpha})u_{t}\Bigr)_{i,j}^{n+\frac{1}{2}}\\
&+O(\tau^{5}+\tau^{3}(h_{x}^2+h_{y}^2)).
\end{aligned}\end{equation}
Adding $(\ref{s1:e6})$ to $(\ref{s1:e5})$, we obtain the alternative direction iteration (ADI) scheme as follows
{\small\begin{equation}\label{s1:e7}
\begin{aligned}
(1+\tau K_{u} \delta_{x}^{\alpha})(1+\tau K_{u}\delta_{y}^{\alpha})u_{i,j}^{n+1}
&=(1-\tau K_{u} \delta_{x}^{\alpha})(1-\tau K_{u}\delta_{y}^{\alpha})u_{i,j}^{n}-
\frac{\tau}{4}u_{i,j}^{n+\frac{1}{2}}(3v_{i,j}^{n}-v_{i,j}^{n-1})^{2}\\
&\quad +F\tau(1-u_{i,j}^{n+\frac{1}{2}})+O(\tau^{3}+\tau^{3} h_{x}^{2}+\tau^{3} h_{y}^{2}).
\end{aligned}
\end{equation}}
Similarly, we can obtain the discretization scheme for $v$
{\small\begin{equation}\label{s1:ee7}
\begin{aligned}
(1+\tau K_{v} \delta_{x}^{\alpha})(1+\tau K_{v}\delta_{y}^{\alpha})v_{i,j}^{n+1}
&=(1-\tau K_{v} \delta_{x}^{\alpha})(1-\tau K_{v}\delta_{y}^{\alpha})v_{i,j}^{n}+
\frac{\tau}{4}u_{i,j}^{n+\frac{1}{2}}(3v_{i,j}^{n}-v_{i,j}^{n-1})^{2}\\
&\quad -\tau(F+\kappa)v_{i,j}^{n+\frac{1}{2}}+O(\tau^{3}+\tau^{3} h_{x}^{2}+\tau^{3} h_{y}^{2}),
\end{aligned} 
\end{equation}}
where $K_{v}=\frac{\mu_{v}}{4\cos(\frac{\pi \alpha}{2})}$.

We replace $u_{i,j}^{n}, v_{i,j}^{n}$ and $u_{i,j}^{n+\frac{1}{2}}$ by the numerical approximations $U_{i,j}^{n}, V_{i,j}^{n}$ and $U_{i,j}^{n+\frac{1}{2}}$,  and obtain an ADI finite difference scheme for the model $(\ref{s1:e1})$ for $n=0$,
\begin{equation}\label{s1:eee8}
\begin{aligned}
(1+\tau K_{u} \delta_{x}^{\alpha})(1+\tau K_{u} \delta_{y}^{\alpha})U_{i,j}^{1}
&=(1-\tau K_{u} \delta_{x}^{\alpha})(1-\tau K_{u}\delta_{y}^{\alpha})U_{i,j}^{0}+F\tau(1-U_{i,j}^{\frac{1}{2}})\\
&\quad -\tau U_{i,j}^{\frac{1}{2}}(V^{0}_{i,j})^2,
\end{aligned}
\end{equation}
\begin{equation}\label{s1:e9}
\begin{aligned}
(1+\tau K_{v} \delta_{x}^{\alpha})(1+\tau K_{v} \delta_{y}^{\alpha})V_{i,j}^{1}
&=(1-\tau K_{v} \delta_{x}^{\alpha})(1-\tau K_{v}\delta_{y}^{\alpha})V_{i,j}^{0}-\tau(F+\kappa)V_{i,j}^{\frac{1}{2}}\\
&\quad +\tau U_{i,j}^{\frac{1}{2}}(V^{0}_{i,j})^{2},
\end{aligned}
\end{equation}
and for $1\leq n\leq M-1$,
\begin{equation}\label{s1:e8}
\begin{aligned}
(1+\tau K_{u} \delta_{x}^{\alpha})(1+\tau K_{u} \delta_{y}^{\alpha})U_{i,j}^{n+1}
&=(1-\tau K_{u} \delta_{x}^{\alpha})(1-\tau K_{u}\delta_{y}^{\alpha})U_{i,j}^{n}+F\tau(1-U_{i,j}^{n+\frac{1}{2}})\\
&\quad -\frac{\tau}{4}U_{i,j}^{n+\frac{1}{2}}(3V_{i,j}^{n}-V_{i,j}^{n-1})^{2},
\end{aligned}
\end{equation}
\begin{equation}\label{s1:e10}
\begin{aligned}
(1+\tau K_{v} \delta_{x}^{\alpha})(1+\tau K_{v} \delta_{y}^{\alpha})V_{i,j}^{n+1}
&=(1-\tau K_{v} \delta_{x}^{\alpha})(1-\tau K_{v}\delta_{y}^{\alpha})V_{i,j}^{n}-\tau(F+\kappa)V_{i,j}^{n+\frac{1}{2}}\\
&\quad +\frac{\tau}{4}U_{i,j}^{n+\frac{1}{2}}(3V_{i,j}^{n}-V_{i,j}^{n-1})^{2}.
\end{aligned}
\end{equation}

We define the following matrices 
$$U^{n}_{*}=[U_{i,j}^{n}], \quad V^{n}_{*}=[V_{i,j}^{n}], \quad \text{for} ~1\leq i \leq N_{x}-1,~ 1\leq j \leq N_{y}-1.$$
Therefore, for $0\leq n\leq M-1$, we can rewrite the numerical scheme into the matrix form as follows
\begin{equation}\label{s1:ee9}
\begin{aligned}
&\Bigl(I+\frac{\tau K_{u}}{h_{x}^{\alpha}}B\Bigr)U^{n+1}_{*}\Bigl(I+\frac{\tau K_{u}}{h_{y}^{\alpha}}B\Bigr)
=\Bigl(I-\frac{\tau K_{u}}{h_{x}^{\alpha}}B\Bigr)U^{n}_{*}\Bigl(I-\frac{\tau K_{u}}{h_{y}^{\alpha}}B\Bigr)+H^{n+\frac{1}{2}},\\
&\Bigl(I+\frac{\tau K_{v}}{h_{x}^{\alpha}}B\Bigr)V^{n+1}_{*}\Bigl(I+\frac{\tau K_{v}}{h_{y}^{\alpha}}B\Bigr)
=\Bigl(I-\frac{\tau K_{v}}{h_{x}^{\alpha}}B\Bigr)V^{n}_{*}\Bigl(I-\frac{\tau K_{v}}{h_{y}^{\alpha}}B\Bigr)+G^{n+\frac{1}{2}},
\end{aligned}
\end{equation}
where 
\begin{equation*}\begin{aligned}
&H_{i,j}^{n+\frac{1}{2}}=\left.\Biggl\{\begin{array}{ll}-\tau U_{i,j}^{\frac{1}{2}}(V^{0}_{i,j})^{2}+F\tau(1-U_{i,j}^{\frac{1}{2}}),&n=0,\\[2mm]
-\frac{\tau}{4}U_{i,j}^{n+\frac{1}{2}}(3V_{i,j}^{n}-V_{i,j}^{n-1})^{2}+F\tau(1-U_{i,j}^{n+\frac{1}{2}}),&1\leq n\leq M-1,\end{array}\right.\\
&G_{i,j}^{n+\frac{1}{2}}=\left.\Biggl\{\begin{array}{ll}\tau U_{i,j}^{\frac{1}{2}}(V^{0}_{i,j})^{2}-\tau(F+\kappa)V_{i,j}^{\frac{1}{2}},& n=0,\\[2mm]
\frac{\tau}{4}U_{i,j}^{n+\frac{1}{2}}(3V_{i,j}^{n}-V_{i,j}^{n-1})^{2}-\tau(F+\kappa)V_{i,j}^{n+\frac{1}{2}},&1\leq n\leq M-1,\end{array}\right.\\
&B=A+A^{T},
\end{aligned}
\end{equation*}
with 
\begin{equation*}
A=\left( \begin{array}{ccccc}
\omega_{1}^{(\alpha)} &\omega_{0}^{(\alpha)}&&&\\
\omega_{2}^{(\alpha)}&\omega_{1}^{(\alpha)}&\omega_{0}^{(\alpha)}&&\\
\vdots&\omega_{2}^{(\alpha)}&\omega_{1}^{(\alpha)}&\ddots&\\
\omega_{n-2}^{(\alpha)}&\cdots&\ddots&\ddots&\omega_{0}^{(\alpha)}\\
\omega_{n-1}^{(\alpha)}&\omega_{n-2}^{(\alpha)}&\cdots&\omega_{2}^{(\alpha)}&\omega_{1}^{(\alpha)}
\end{array}\right).
\end{equation*}

So we can use the alternative direction iteration (ADI) method to solve $(\ref{s1:ee9})$ efficiently. In addition, since the coefficient matrices in this numerical scheme are Toeplitz matrices, a class of fast algorithms introduced in \cite{wang2012fast,wangdu2014fast} can be applied to accelerate the computation.
\subsection{Stability analysis}\label{s3}

In this part, we consider the stability of the time semi-discrete scheme of the fractional GS model $(\ref{s1:e1})$. We represent $u(x,y,t^{n}), v(x,y,t^{n})$ by the notations $U^{n}, V^{n}$, which are the solutions of the time semi-discrete scheme. Then the time semi-discrete scheme can be written as follows
{\small\begin{equation}\label{s3:e10} 
\frac{U^{n+1}-U^{n}}{\tau}=\left.\Biggl\{\begin{array}{ll}
-\mu_{u}(-\Delta)^{\frac{\alpha}{2}}U^{\frac{1}{2}}
-U^{\frac{1}{2}}(V^{0})^{2}+F(1-U^{\frac{1}{2}}),&n=0,\\[2mm]
-\mu_{u}(-\Delta)^{\frac{\alpha}{2}}U^{n+\frac{1}{2}}
-U^{n+\frac{1}{2}}(V^{*,n+\frac{1}{2}})^{2}+F(1-U^{n+\frac{1}{2}}), &1\leq n \leq M-1, \end{array}\right.
\end{equation}}
{\small\begin{equation}\label{s3:e11} 
\frac{V^{n+1}-V^{n}}{\tau}=\left.\Biggl\{\begin{array}{ll}
-\mu_{v}(-\Delta)^{\frac{\alpha}{2}}V^{\frac{1}{2}}
+U^{\frac{1}{2}}(V^{0})^{2}-(F+\kappa)V^{\frac{1}{2}},&n=0,\\[2mm]
-\mu_{v}(-\Delta)^{\frac{\alpha}{2}}V^{n+\frac{1}{2}}
+U^{n+\frac{1}{2}}(V^{*,n+\frac{1}{2}})^{2}-(F+\kappa)V^{n+\frac{1}{2}},&1\leq n \leq M-1, \end{array}\right.
\end{equation}}
where $U^{n+\frac{1}{2}}=\frac{U^{n}+U^{n+1}}{2},~V^{n+\frac{1}{2}}=\frac{V^{n}+V^{n+1}}{2}$ and $V^{*,n+\frac{1}{2}}=\frac{3V^{n}-V^{n-1}}{2}$ with $1\leq n\leq M-1$. 

\begin{theorem}\label{s4:l1}
The time semi-discrete scheme of the fractional GS model is unconditionally stable for $1< \alpha \leq2$, and the following estimates hold for $0\leq n\leq M-1$,
\begin{equation}\label{tt:4}
\|U^{n+1}\|^{2} \leq \|u_{0}\|^2+FT|\Omega|,
\end{equation}
\begin{equation}\label{tt:5}\begin{aligned}
\|W^{n+1}\|^{2}\leq & \|u_{0}+v_{0}\|^{2}+(\kappa T+\frac{|\mu_{v}-\mu_{u}|^{2}}{2\mu_{u}\mu_{v}|\cos(\frac{\pi \alpha}{2})|^2})\|u_{0}\|^{2}\\
&+(1+\kappa T+\frac{|\mu_{v}-\mu_{u}|^{2}}{2\mu_{u}\mu_{v}|\cos(\frac{\pi \alpha}{2})|^2})FT|\Omega|,
\end{aligned}\end{equation}
where $W^{n}=U^{n}+V^{n}$, which is the solution of the following equation 
\begin{equation}\label{tt:6}\begin{aligned}
\frac{W^{n+1}-W^{n}}{\tau}=&-\mu_{v}(-\Delta)^{\frac{\alpha}{2}}W^{n+\frac{1}{2}}-(F+\kappa)W^{n+\frac{1}{2}}\\
&+\Big[(\mu_{v}-\mu_{u})(-\Delta)^{\frac{\alpha}{2}}U^{n+\frac{1}{2}}+\kappa U^{n+\frac{1}{2}}+F\Bigr].
\end{aligned}\end{equation}
\end{theorem}

The proof is presented in Appendix \ref{adx2}.

%%%%%%%%%%%%%%%%%%%%%%%%%%%%%%%%%%
\section{Numerical experiments}\label{s6}
\setcounter{equation}{0}

In this section, we carry out numerical experiments to verify the accuracy of the numerical method proposed in Section \ref{s1} and use the scheme in numerical simulations for the fractional GS model to study the pattern formation.

\subsection{Convergence tests}
In this subsection, we conduct two numerical examples runs for the benchmark problems to test the accuracy of the numerical scheme.
\begin{example}\label{s6:e1}
We consider the following fractional diffusion problem
\begin{equation}\label{s6:ee1}\begin{aligned}
&\frac{\partial{u(x,y,t)}}{\partial t}=-(-\Delta)^{\frac{\alpha}{2}}u(x,y,t)+f(x,y,t),&&(x,y)\in \Omega, t\in [0,T],\\
&u(x,y,0)=x^{4}(1-x)^{4}y^{4}(1-y)^{4},&&  (x,y)\in \Omega,\\
&u(x,y,t)|_{\partial \Omega}=0, && t\in [0,T].
\end{aligned}\end{equation}
The domain is $\Omega=(0,1)^{2}$. The source term is
\begin{equation*}\begin{aligned}
f(x,y,t)=\Bigl[&\Bigl(-x^{4}(1-x)^{4}y^{4}(1-y)^{4}\Bigr)+Ky^{4}(1-y)^{4}\Bigl(\frac{\Gamma(5)}{\Gamma(5-\alpha)}
(x^{4-\alpha}+(1-x)^{4-\alpha})\\
&-\frac{4\Gamma(6)}{\Gamma(6-\alpha)}(x^{5-\alpha}+(1-x)^{5-\alpha})+\frac{6\Gamma(7)}{\Gamma(7-\alpha)}(x^{6-\alpha}+(1-x)^{6-\alpha})\\
&-\frac{4\Gamma(8)}{\Gamma(8-\alpha)}(x^{7-\alpha}+(1-x)^{7-\alpha})+\frac{\Gamma(9)}{\Gamma(9-\alpha)}(x^{8-\alpha}+(1-x)^{8-\alpha})\Bigr)\\
&+K\Bigl(\frac{\Gamma(5)}{\Gamma(5-\alpha)}
(y^{4-\alpha}+(1-y)^{4-\alpha})-\frac{4\Gamma(6)}{\Gamma(6-\alpha)}(y^{5-\alpha}+(1-y)^{5-\alpha})\\
&+\frac{6\Gamma(7)}{\Gamma(7-\alpha)}(y^{6-\alpha}+(1-y)^{6-\alpha})-\frac{4\Gamma(8)}{\Gamma(8-\alpha)}(y^{7-\alpha}+(1-y)^{7-\alpha})\\
&+\frac{\Gamma(9)}{\Gamma(9-\alpha)}(y^{8-\alpha}+(1-y)^{8-\alpha})\Bigr)x^{4}(1-x)^{4}\Bigr]e^{-t},\\
\end{aligned}\end{equation*}
where the coefficient $K=(2\cos(\pi \alpha/2))^{-1}$.
The exact solution of this fractional equation is $u(x,y,t)=e^{-t}x^{4}(1-x)^{4}y^{4}(1-y)^{4}$.
\end{example}

In our simulations, we perform two types of numerical tests. Firstly, we use a fine time step $\tau=\frac{1}{3000}$ and refine $h=h_{x}=h_{y}$ from $\frac{1}{16}$ to $\frac{1}{256}$ to observe the spatial convergence rates. We present the numerical results in Table \ref{table1} and observe the second-order accuracy in space for any $\alpha$. In addition, we use a fine spatial mesh size $h=\frac{1}{1024}$ and refine $\tau$ from $\frac{1}{5}$ to $\frac{1}{25}$ to observe the temporal convergence rates. We present the numerical results in Table \ref{table2} and observe the second-order accuracy in time for any $\alpha$.
\begin{table}[h]
\renewcommand{\arraystretch}{1.2}
\centering  
\begin{tabular}{p{1cm} p{1cm}<{\centering} p{3cm}<{\centering} p{2cm}<{\centering}}  
\hline
$\alpha$ &$h$ &$\frac{\|u^{n}-U^{n}\|}{\|u^{n}\|}$ &rate\\ \hline  
1.2 &$\frac{1}{16}$ &0.023 &--\\
&$\frac{1}{32}$&0.0059&1.9704\\
&$\frac{1}{64}$&0.0015&1.9799\\
&$\frac{1}{128}$&3.7586E-04&1.9882\\
&$\frac{1}{256}$&9.4409E-05&1.9932\\ \hline        
1.5 &$\frac{1}{16}$ &0.0216&--\\ 
&$\frac{1}{32}$&0.0055&1.9815\\
&$\frac{1}{64}$&0.0014&1.9867\\
&$\frac{1}{128}$&3.4652E-04&1.9919\\
&$\frac{1}{256}$&8.6947E-05&1.9947\\ \hline
1.8&$\frac{1}{16}$&0.0171&--\\
&$\frac{1}{32}$&0.0043&1.9968\\
&$\frac{1}{64}$&0.0011&1.9962\\
&$\frac{1}{128}$&2.6846e-04&1.9970\\      
&$\frac{1}{256}$&6.7289e-05&1.9963\\ \hline
\end{tabular}
\caption{Spatial $L^{2}$ errors and their corresponding convergence rates for Example $\ref{s6:e1}$ at $t=1$ for $\alpha=1.2, 1.5, 1.8$ with time step $\tau =\frac{1}{3000}$.\label{table1}}
\end{table}

\begin{table}[htbp]
\renewcommand{\arraystretch}{1.2}
\centering  
\begin{tabular}{p{1cm} p{1cm}<{\centering} p{3cm}<{\centering} p{2cm}<{\centering}}  
\hline
$\alpha$&$\tau$  &$\frac{\|u^{n}-U^{n}\|}{\|u^{n}\|}$ &rate\\ \hline  
1.2 &$\frac{1}{5}$ &0.0196 &--\\
&$\frac{1}{10}$&0.0048&2.0399\\
&$\frac{1}{15}$&0.0021&2.0088\\
&$\frac{1}{20}$&0.0012&1.9986\\
&$\frac{1}{25}$&7.6011E-04&1.9912\\ \hline        
1.5 &$\frac{1}{5}$ &0.0340&--\\ 
&$\frac{1}{10}$&0.0080&2.0850\\
&$\frac{1}{15}$&0.0035&2.0228\\
&$\frac{1}{20}$&0.0020&2.0079\\
&$\frac{1}{25}$&0.0013&2.0003\\ \hline
1.8&$\frac{1}{5}$&0.0613&--\\
&$\frac{1}{10}$&0.0138&2.1536\\
&$\frac{1}{15}$&0.0060&2.0519\\
&$\frac{1}{20}$&0.0034&2.0231\\      
&$\frac{1}{25}$&0.0021&2.0115\\ \hline
\end{tabular}
\caption{Temporal $L^{2}$ errors and their corresponding convergence rates for Example $\ref{s6:e1}$ at $t=1$ for $\alpha=1.2, 1.5, 1.8$ with spatial partition $h=\frac{1}{1024}$.\label{table2}}
\end{table}

\begin{example}\label{s6:e2}
We investigate the accuracy of the numerical scheme for the same problem  $(\ref{s6:ee1})$ with the initial condition $u(x,y,0)=sin(\pi x)sin(\pi y)$ and the boundary condition $u(x,y,t)|_{\partial \Omega}=0, t>0$. The source term is $f(x,y,t)=1$. 
\end{example}

\begin{figure}[htp]
\begin{center}
\includegraphics[width=3.8in,height=2.75in]{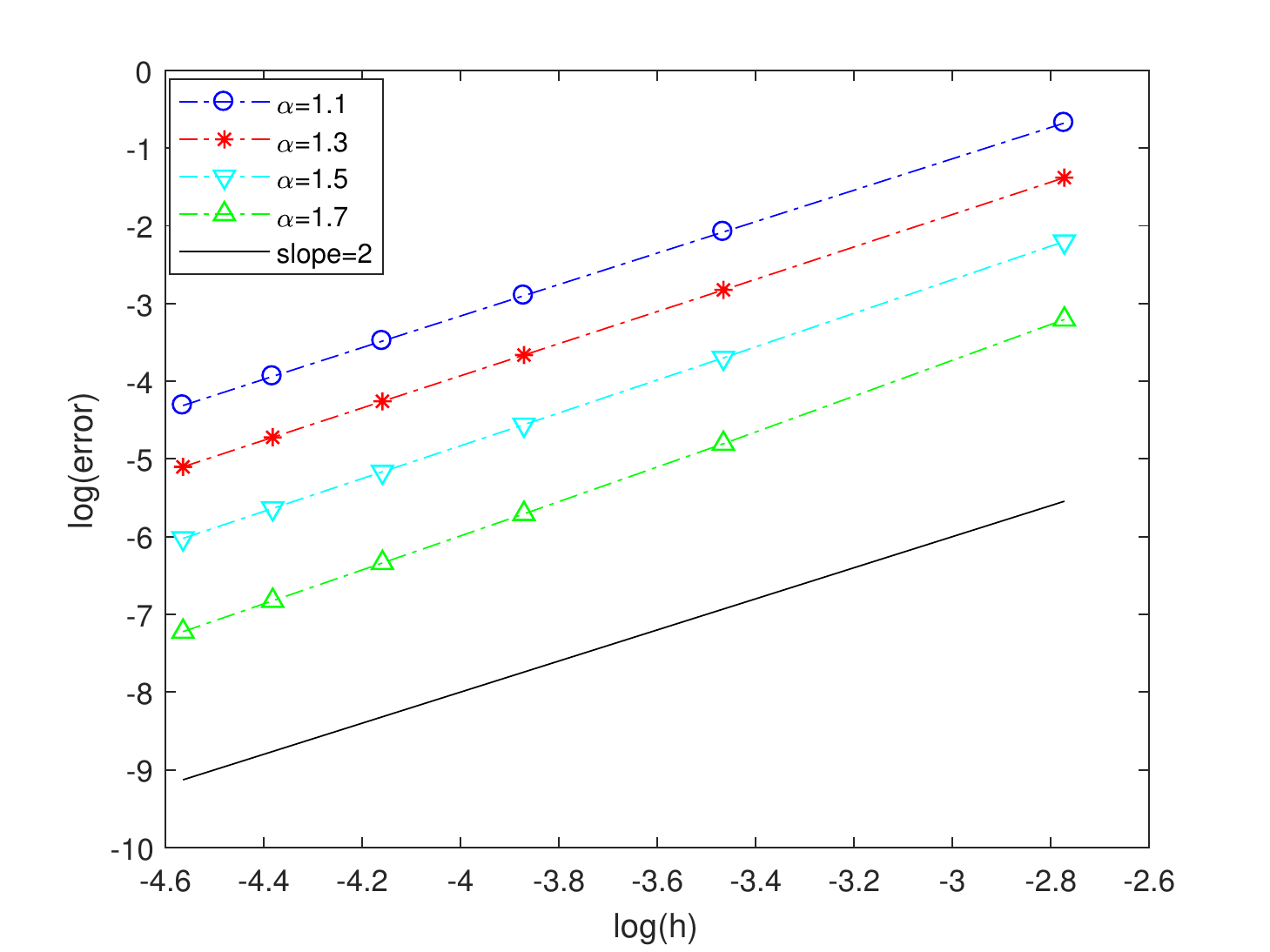} 
\caption{\label{Fig. 2}$L^2$ errors between numerical solutions and reference solutions $(h=\frac{1}{1024})$ as a function of $h$ in log-log scale for Example \ref{s6:e2} at $t=1$ for $\alpha=1.1, 1.3, 1.5, 1.7$ with $\tau=h$. }
%\label{default}
\end{center}
\end{figure}

In this case, the exact solution is unknown. We use the numerical results in finer spatial partition $(h=1/1024)$ as the reference results to study the convergence rates in space and time. We choose $h=h_{x}=h_{y}$ and $\tau=h$. In Fig. $\ref{Fig. 2}$, we plot the $L^2$ errors between the numerical results and reference solutions as the function of $h$ in log-log scale at $t=1$ for several fractional orders $\alpha=1.1, 1.3, 1.5, 1.7$. We observe that the numerical scheme has second-order accuracy both in time and space.

\subsection{The numerical simulations of the fractional GS model}\label{gss}
We carry out the numerical simulations to study the dynamics of the fractional GS model with a perturbation to the spatially homogeneous steady state. The spatially initial condition is 
\begin{equation*}
(u,v)=\left.\Biggl\{\begin{array}{ll} (1,0), \quad(x,y)\in \Omega \backslash O_{c}, \\[2mm]
(\frac{1}{2},\frac{1}{4}), \quad(x,y) \in O_{c},\end{array}\right.
\end{equation*}
where $\Omega=(-1,2)^2$ and $O_{c}=\{(x,y)|(x-0.5)^2+(y-0.5)^2\leq 0.04^2\}$. This initial condition is a perturbation of the steady state $u_{*}=1,v_{*}=0$ imposed with the zero boundary condition. The spatial mesh size is chosen as $h=h_{x}=h_{y}=\frac{1}{1024}$ and time step is $\tau=0.1$. We choose $\mu_{u}=2\times10^{-5}, \mu_v=\mu_u/2, F=0.03$ and vary $\kappa$ in a range (see subsection \ref{2.1}). The model is known to generate different mechanisms of pattern formation when the ratio of diffusion coefficients $\mu_u/\mu_v>1$.

\begin{figure}[h!]
\centering
{\subfloat[$t=1000$]{
\begin{minipage}{.17\textwidth}\centering
\includegraphics[width=1.0\textwidth]{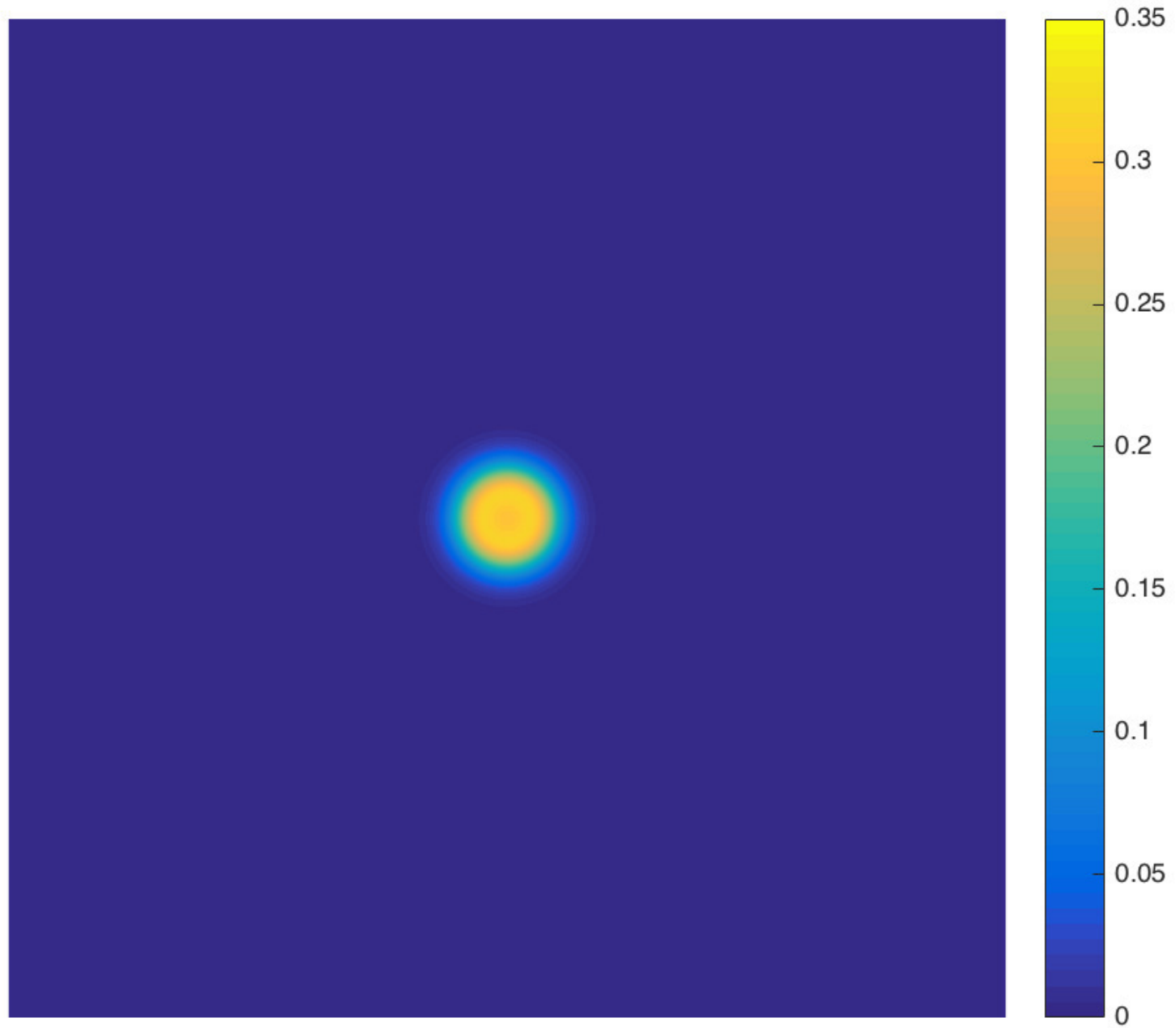}
\end{minipage}}
 \subfloat[$t=6200$]{
\begin{minipage}{.17\textwidth}\centering
\includegraphics[width=1.0\textwidth]{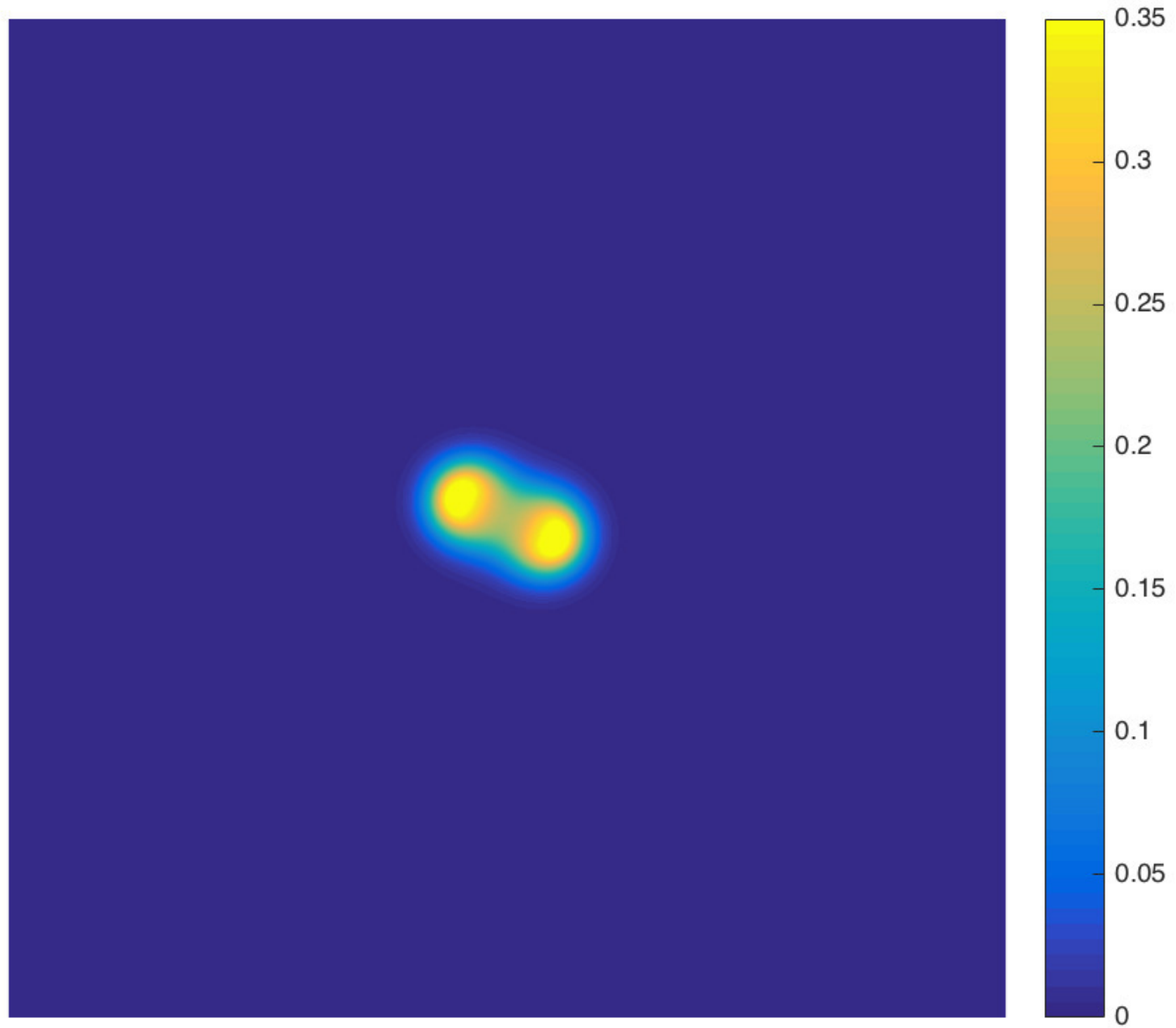}
\end{minipage}}
 \subfloat[$t=6800$]{
\begin{minipage}{.17\textwidth}\centering
\includegraphics[width=1.0\textwidth]{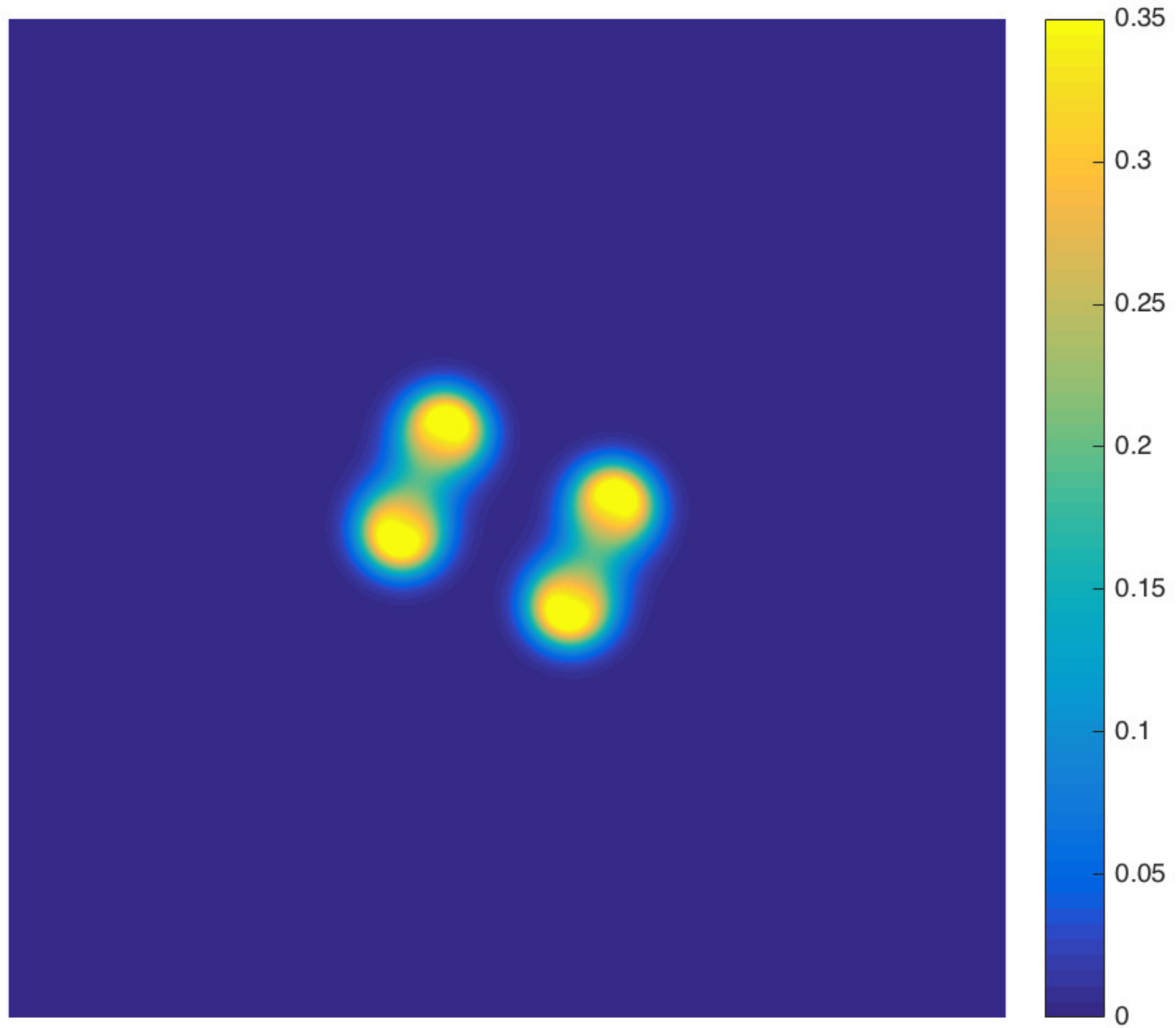}
\end{minipage}}
 \subfloat[$t=9400$]{
\begin{minipage}{.17\textwidth}\centering
\includegraphics[width=1.0\textwidth]{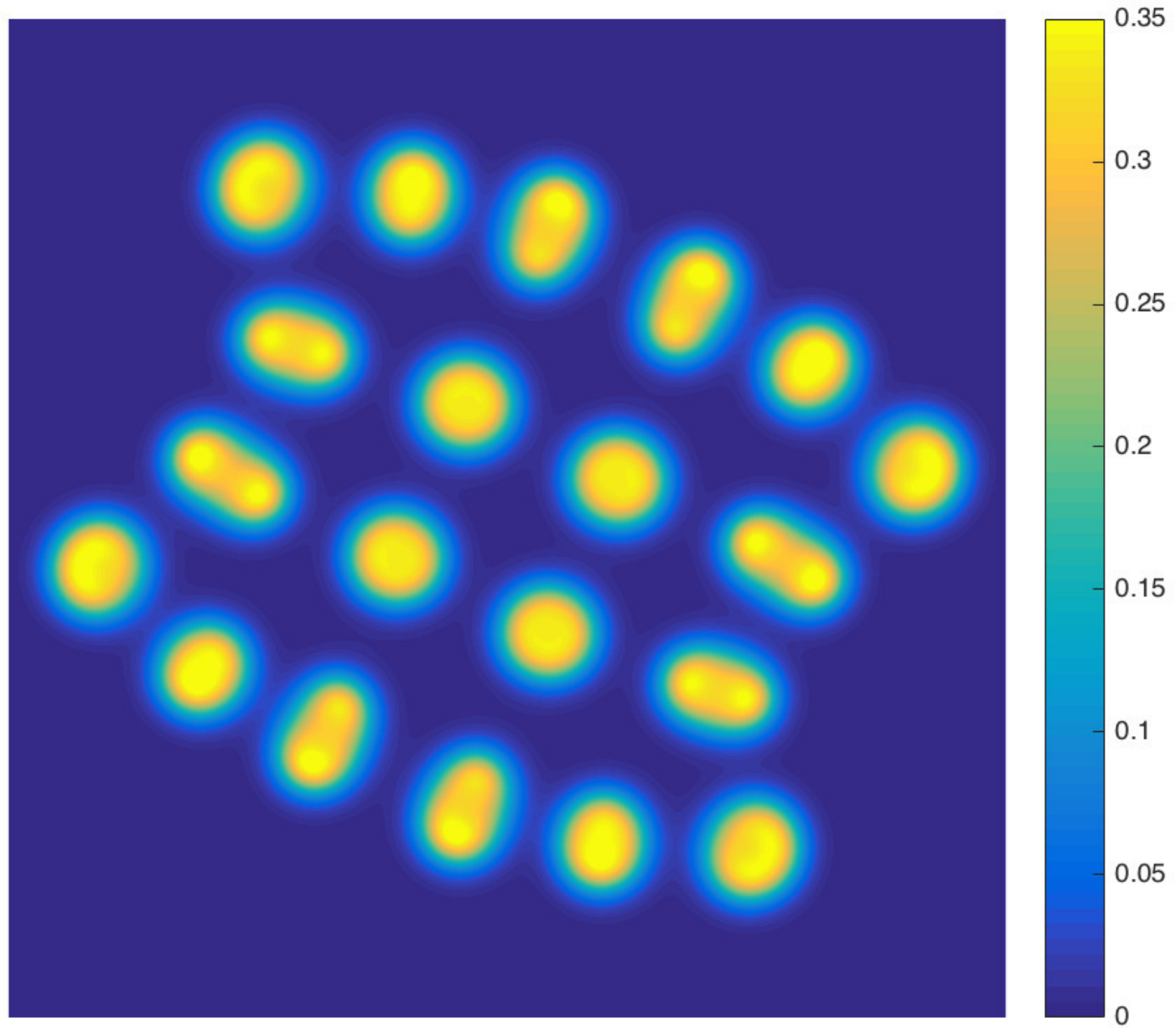}\\
\end{minipage}}
 \subfloat[$t=30000$]{
\begin{minipage}{.17\textwidth}\centering
\includegraphics[width=1.0\textwidth]{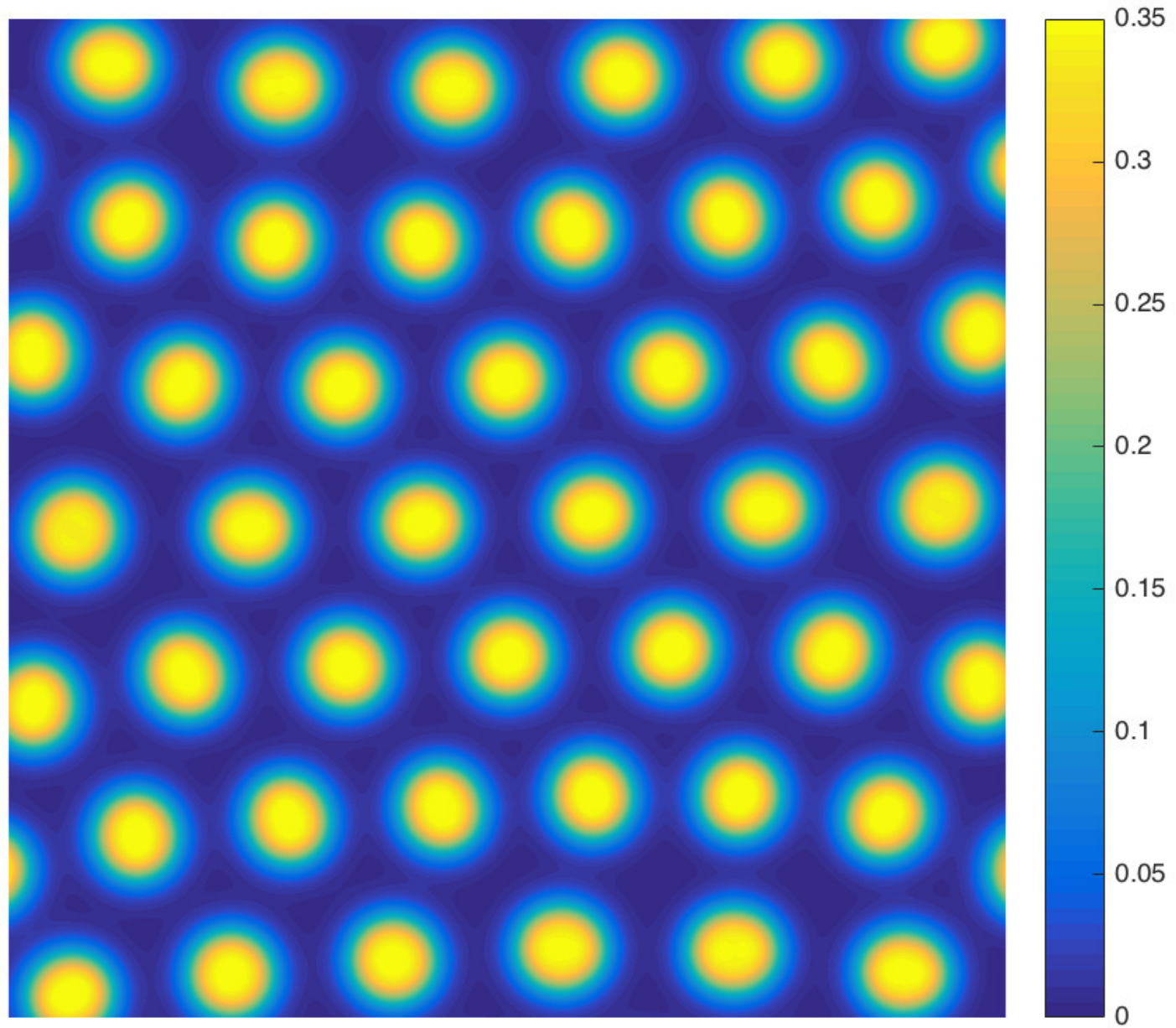}\\
\end{minipage}}}
{\subfloat[$t=200$]{
\begin{minipage}{.17\textwidth}\centering
\includegraphics[width=1.0\textwidth]{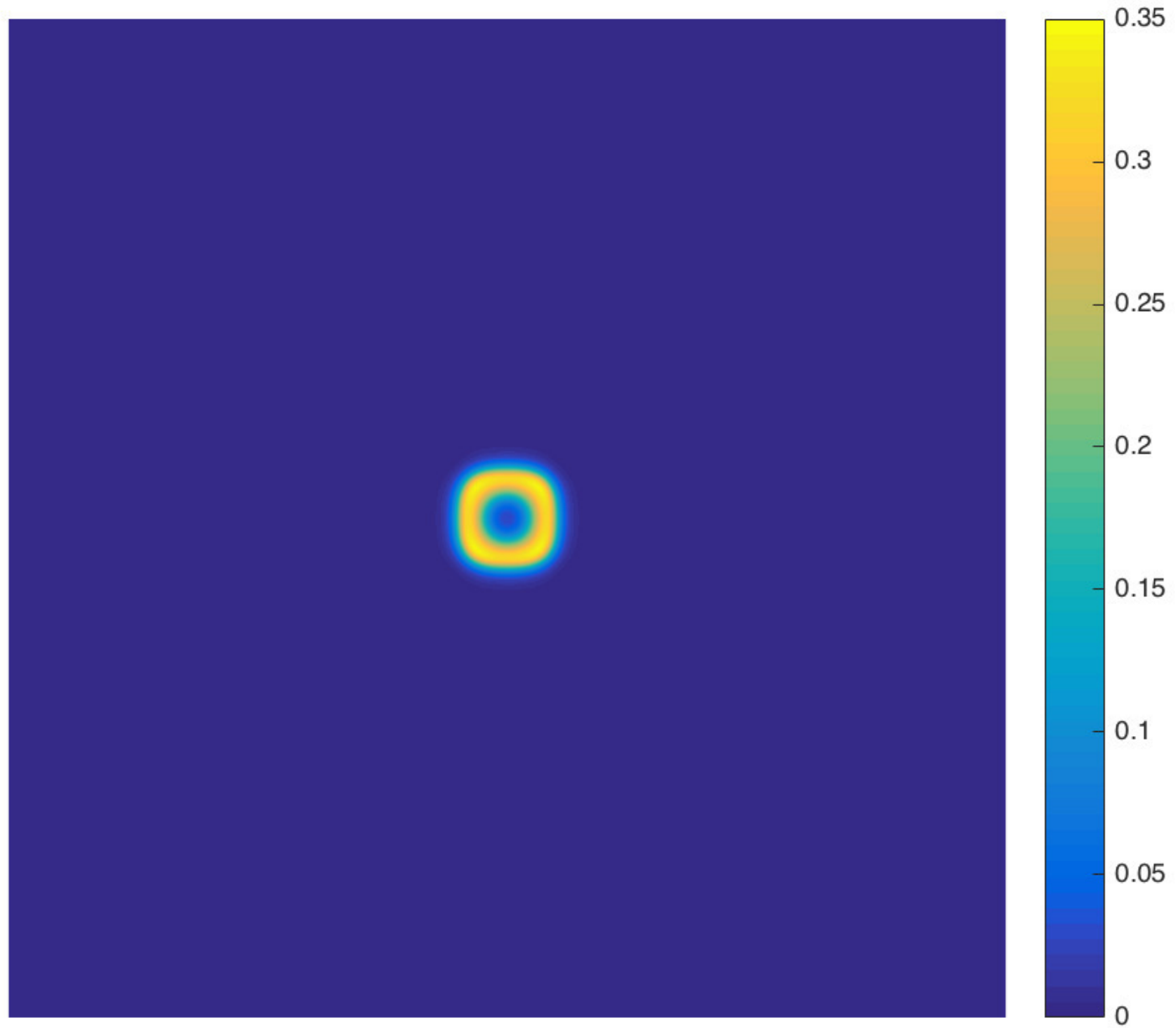}
\end{minipage}}
 \subfloat[$t=2000$]{
\begin{minipage}{.17\textwidth}\centering
\includegraphics[width=1.0\textwidth]{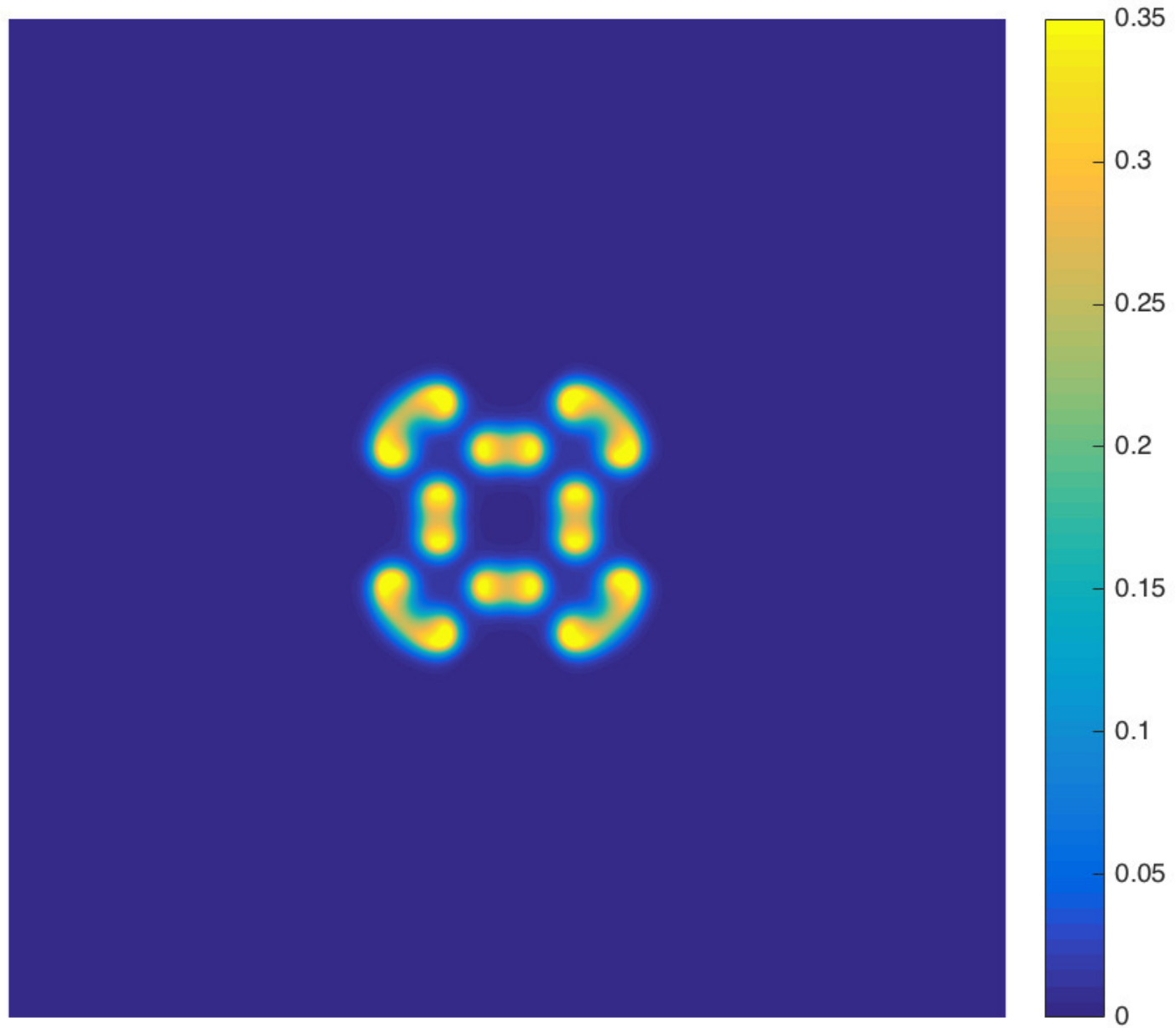}
\end{minipage}}
 \subfloat[$t=3400$]{
\begin{minipage}{.17\textwidth}\centering
\includegraphics[width=1.0\textwidth]{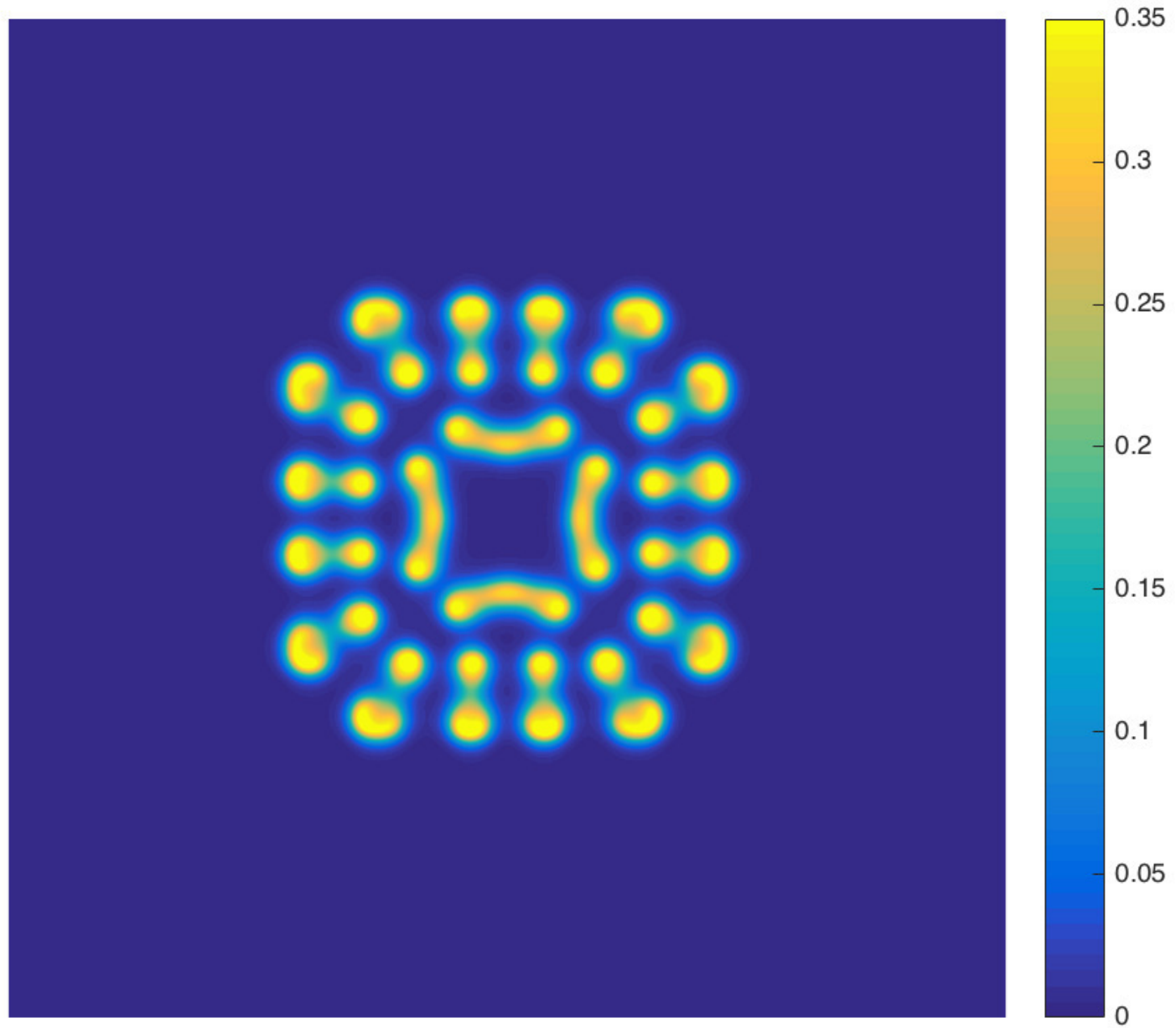}
\end{minipage}}
 \subfloat[$t=6000$]{
\begin{minipage}{.17\textwidth}\centering
\includegraphics[width=1.0\textwidth]{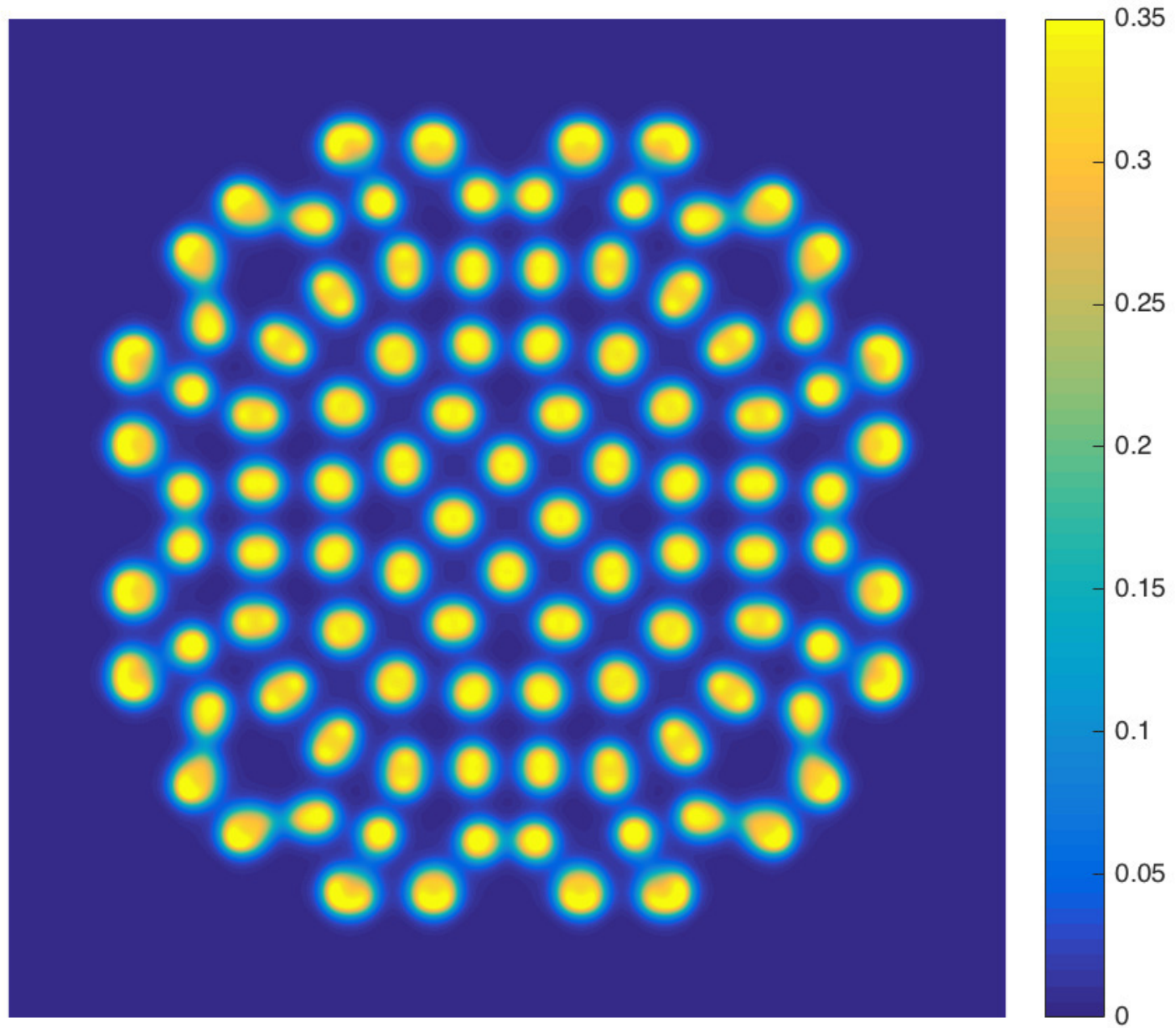}\\
\end{minipage}}
 \subfloat[$t=30000$]{
\begin{minipage}{.17\textwidth}\centering
\includegraphics[width=1.0\textwidth]{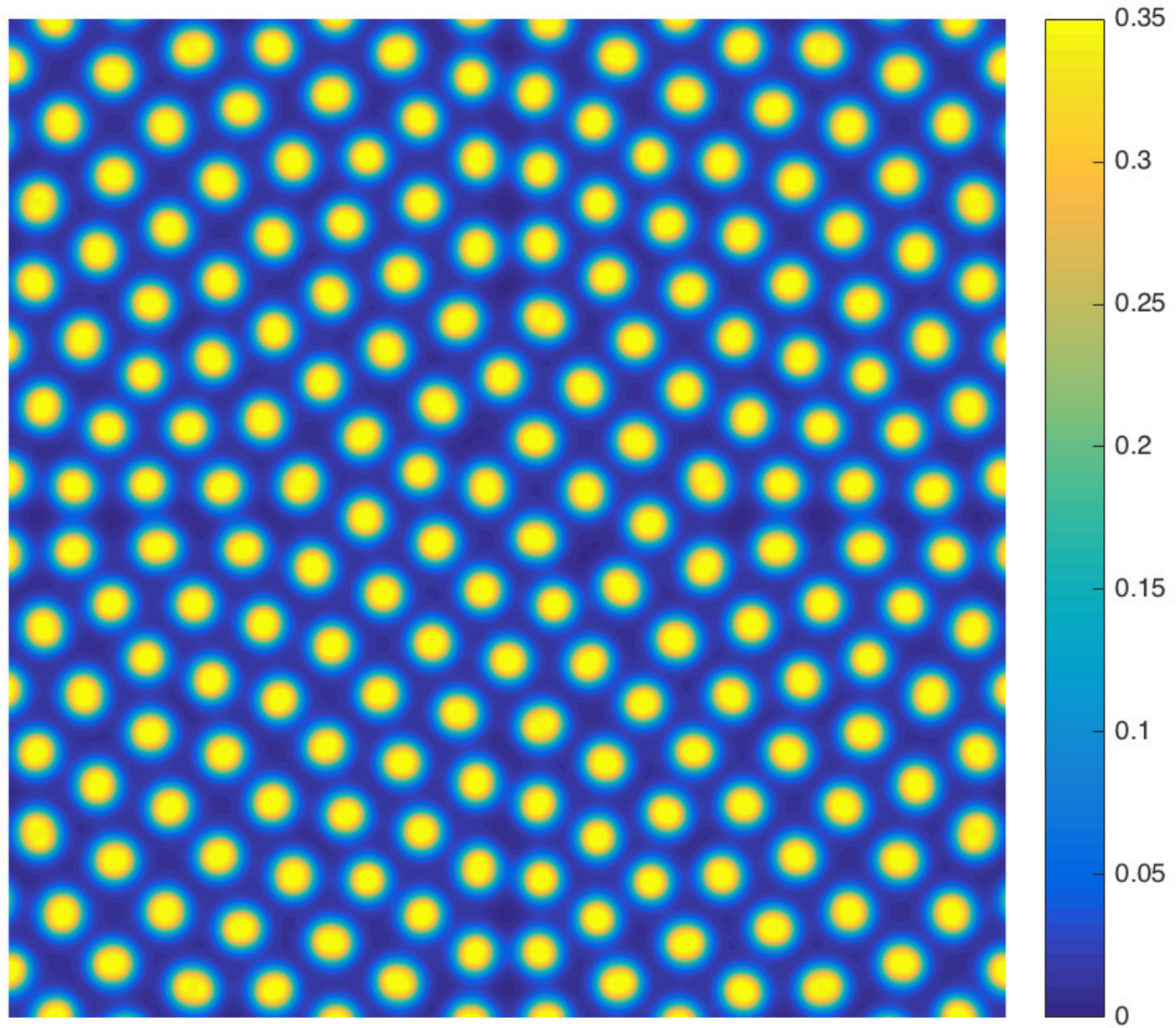}\\
\end{minipage}}}
{\subfloat[$t=2000$]{
\begin{minipage}{.17\textwidth}\centering
\includegraphics[width=1.0\textwidth]{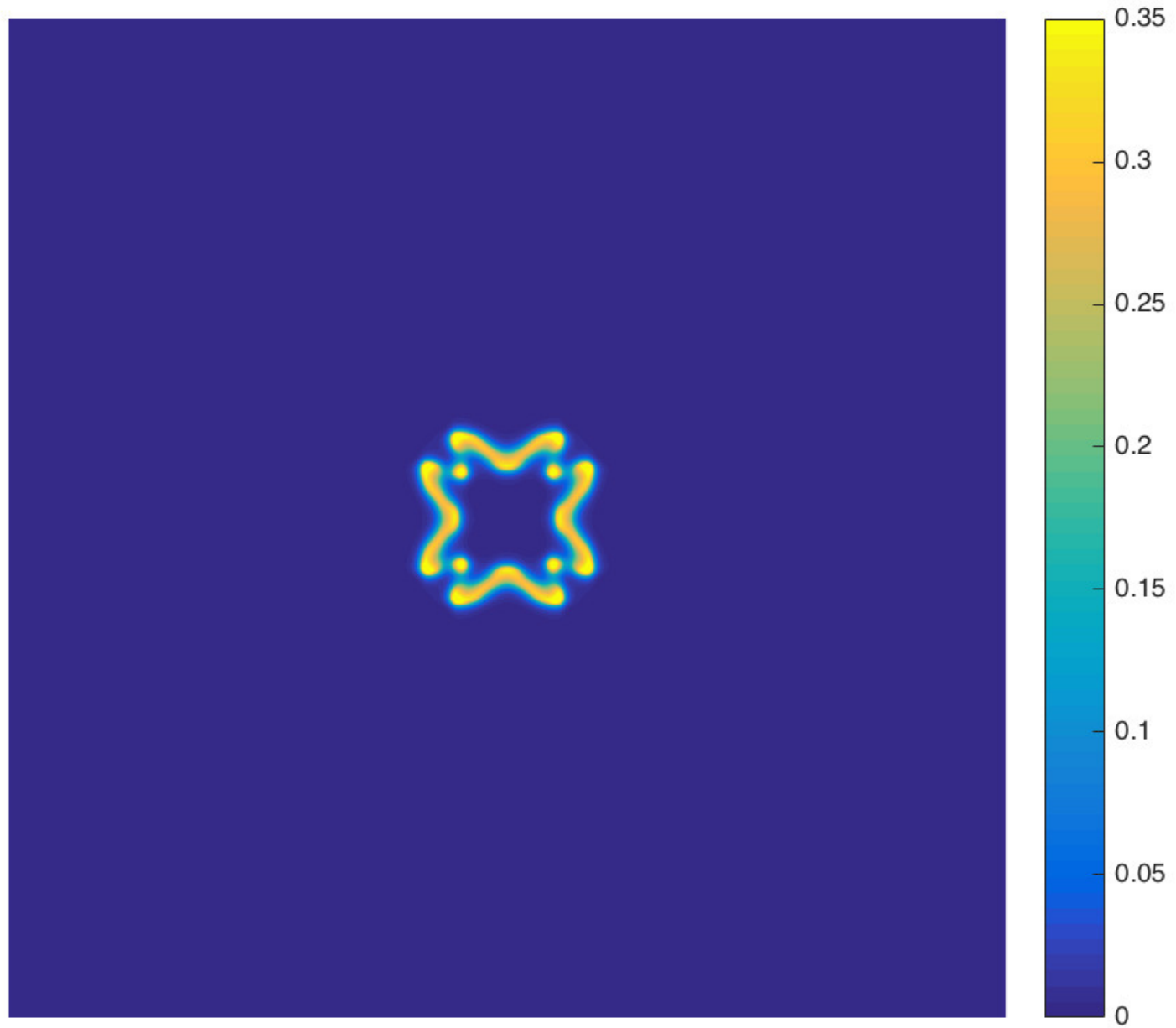}
\end{minipage}}
 \subfloat[$t=3000$]{
\begin{minipage}{.17\textwidth}\centering
\includegraphics[width=1.0\textwidth]{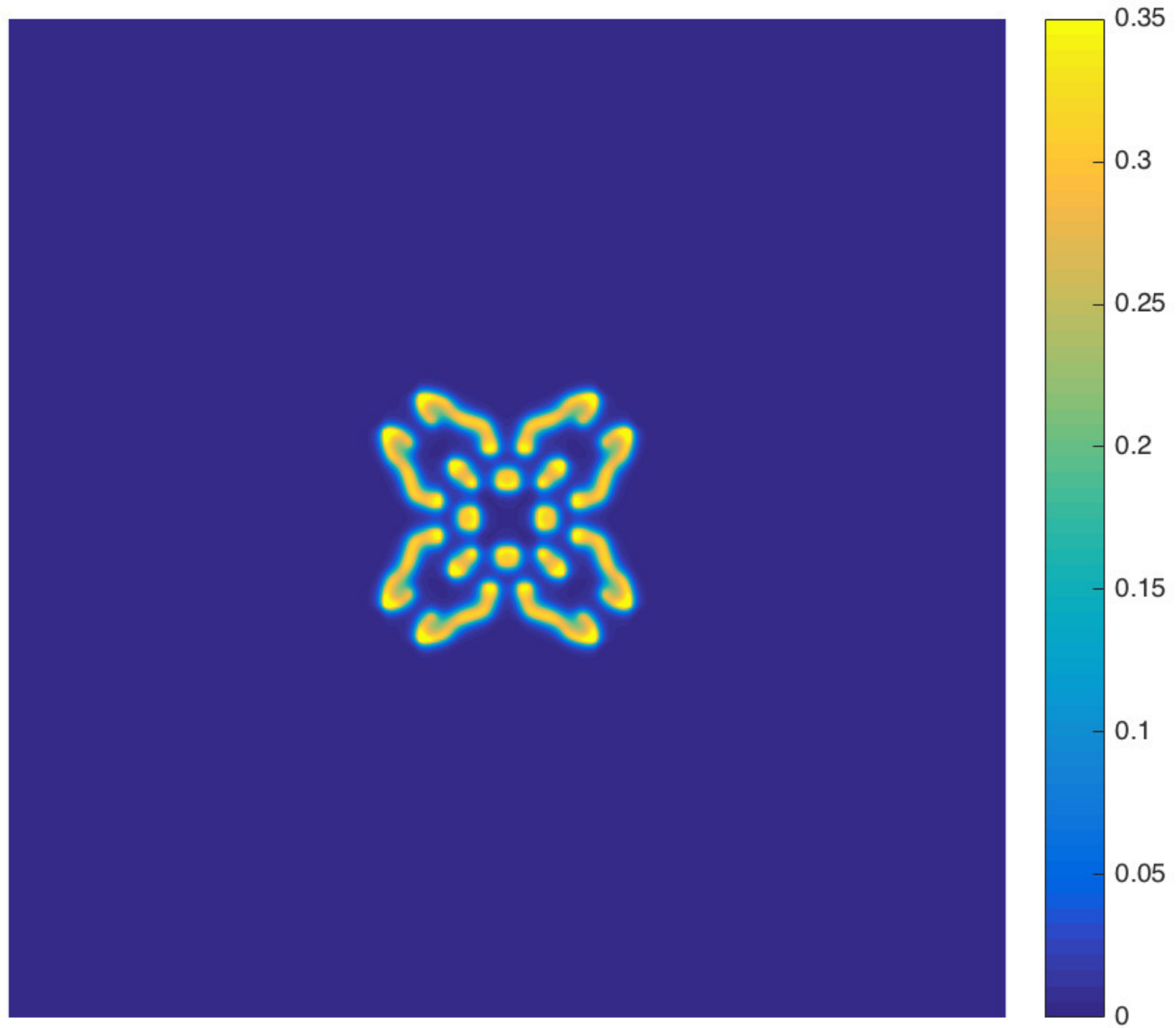}
\end{minipage}}
 \subfloat[$t=5000$]{
\begin{minipage}{.17\textwidth}\centering
\includegraphics[width=1.0\textwidth]{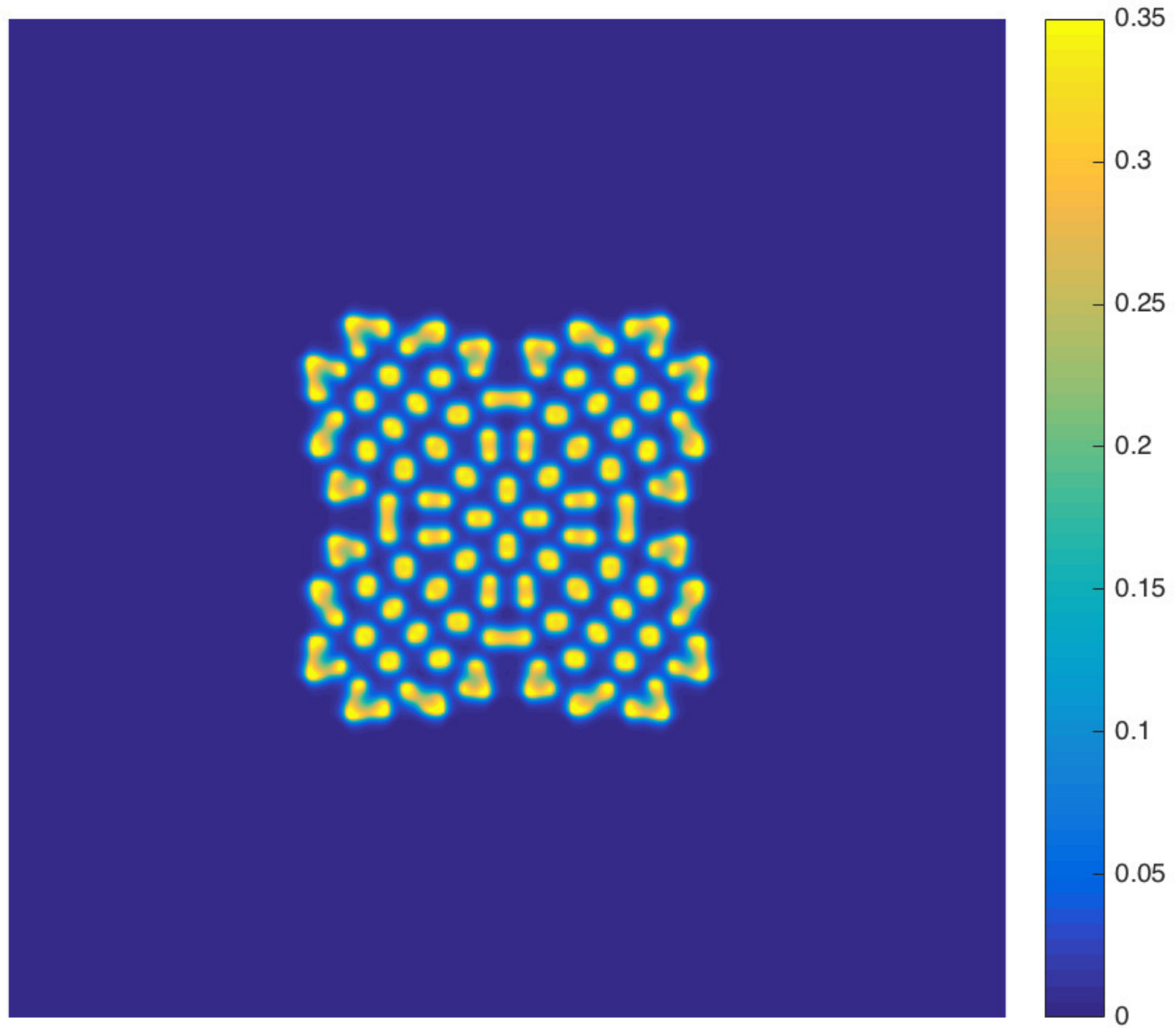}
\end{minipage}}
 \subfloat[$t=9000$]{
\begin{minipage}{.17\textwidth}\centering
\includegraphics[width=1.0\textwidth]{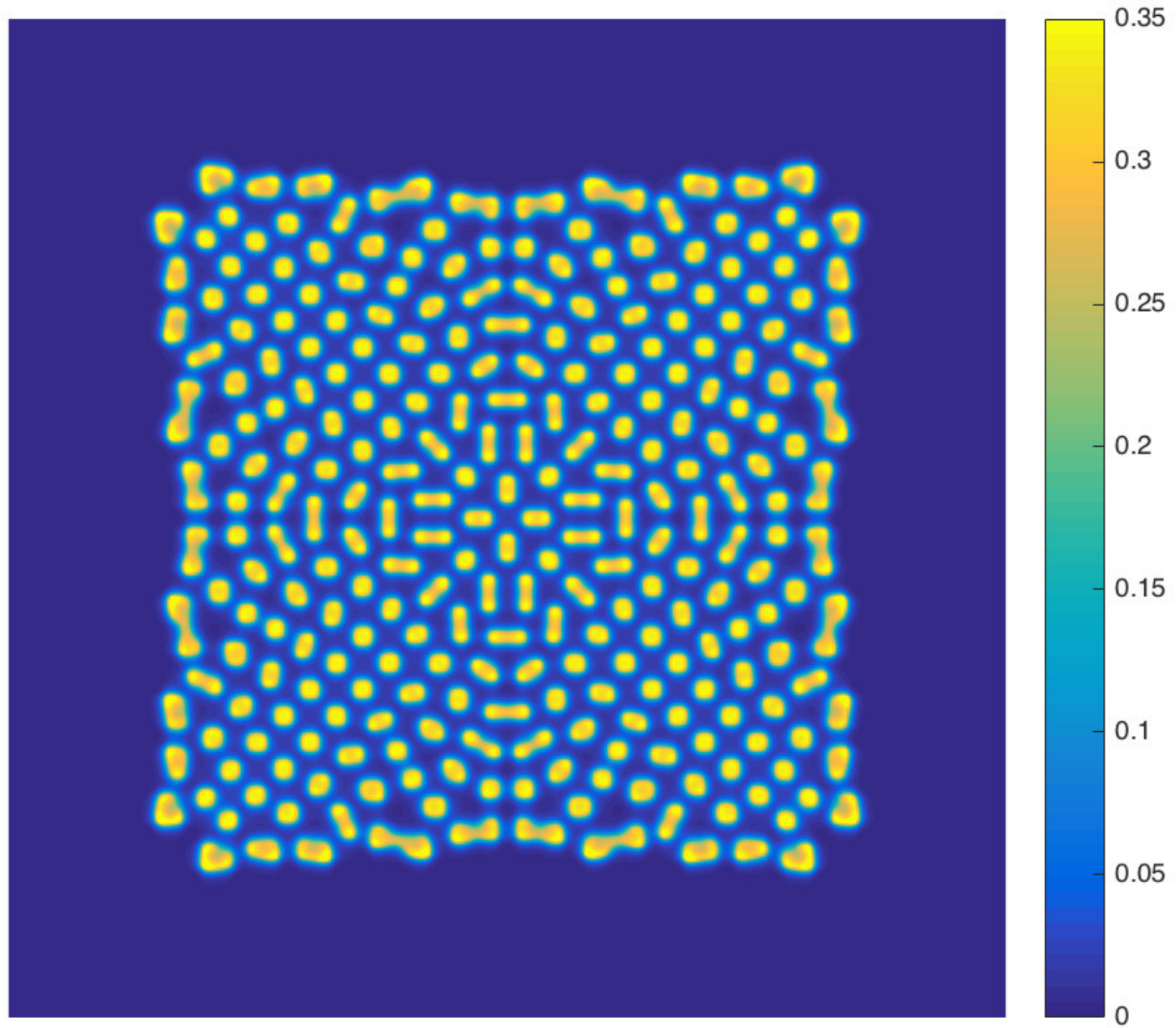}\\
\end{minipage}}
 \subfloat[$t=30000$]{
\begin{minipage}{.17\textwidth}\centering
\includegraphics[width=1.0\textwidth]{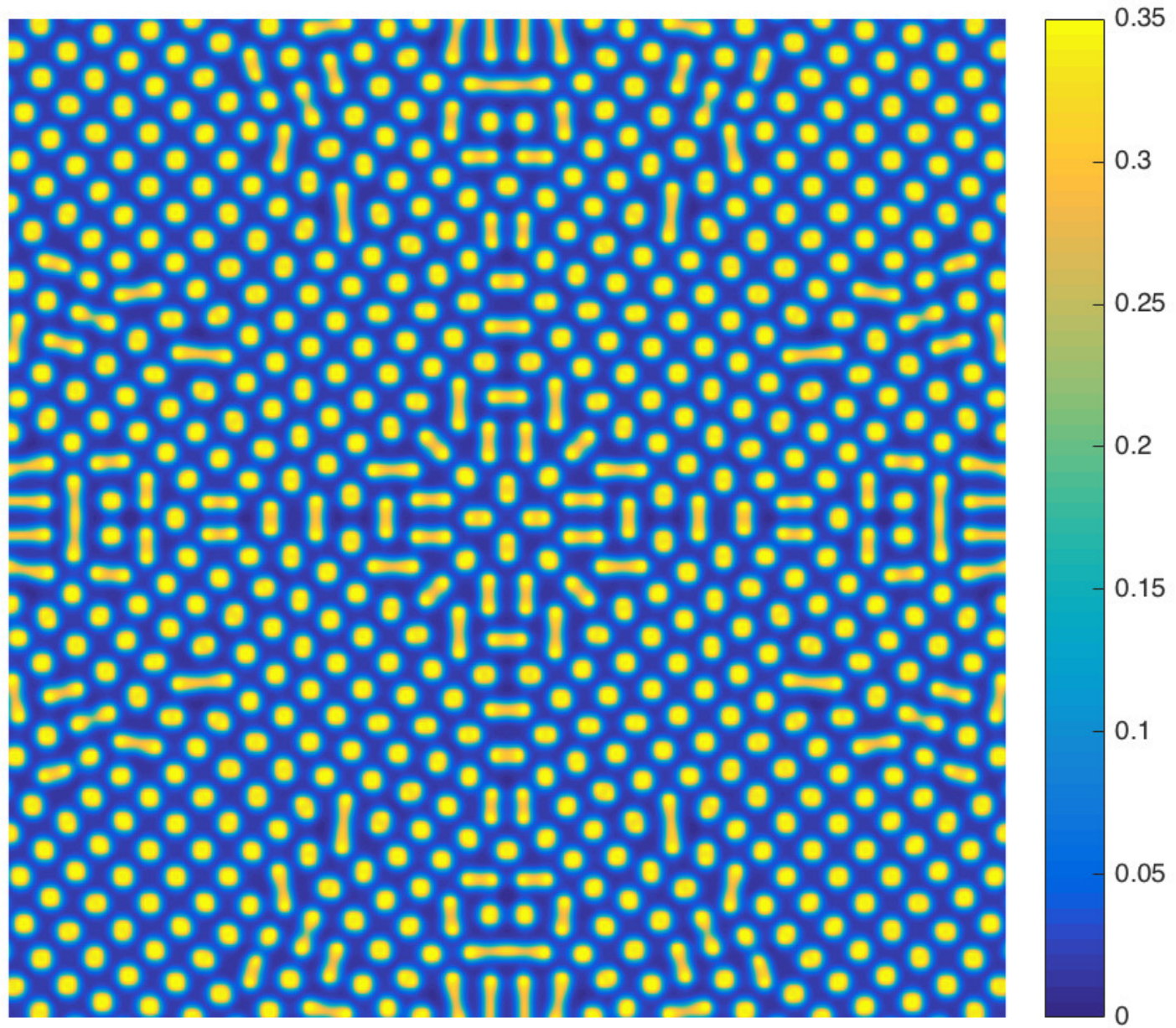}\\
\end{minipage}}}
\caption{Evolution of the solution $v$ in Eq. (\ref{s1:e1}) with $\kappa=0.063$: (a)-(e) $v$ contours are shown at different times from left to right for fractional order $\alpha=2.0$; (f)-(j) $v$ contours are shown at different times from left to right for fractional order $\alpha=1.7$; (k)-(o) $v$ contours are shown at different times from left to right for fractional order $\alpha=1.5$.}
\label{Fig. 3}
\end{figure}

\begin{figure}[!h]
\centering
{\subfloat[$t=200$]{
\begin{minipage}{.17\textwidth}\centering
\includegraphics[width=1.0\textwidth]{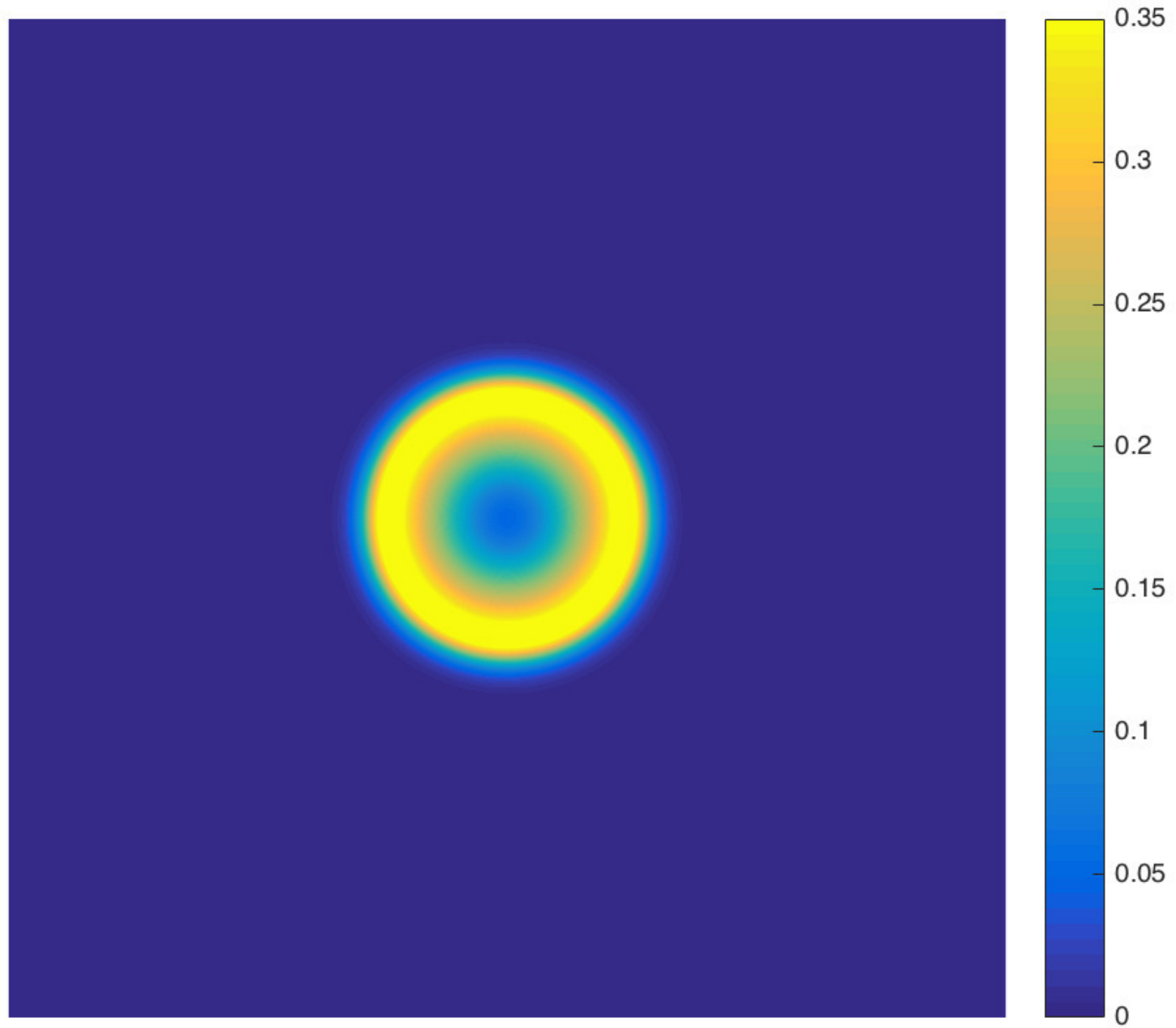}
\end{minipage}}
 \subfloat[$t=800$]{
\begin{minipage}{.17\textwidth}\centering
\includegraphics[width=1.0\textwidth]{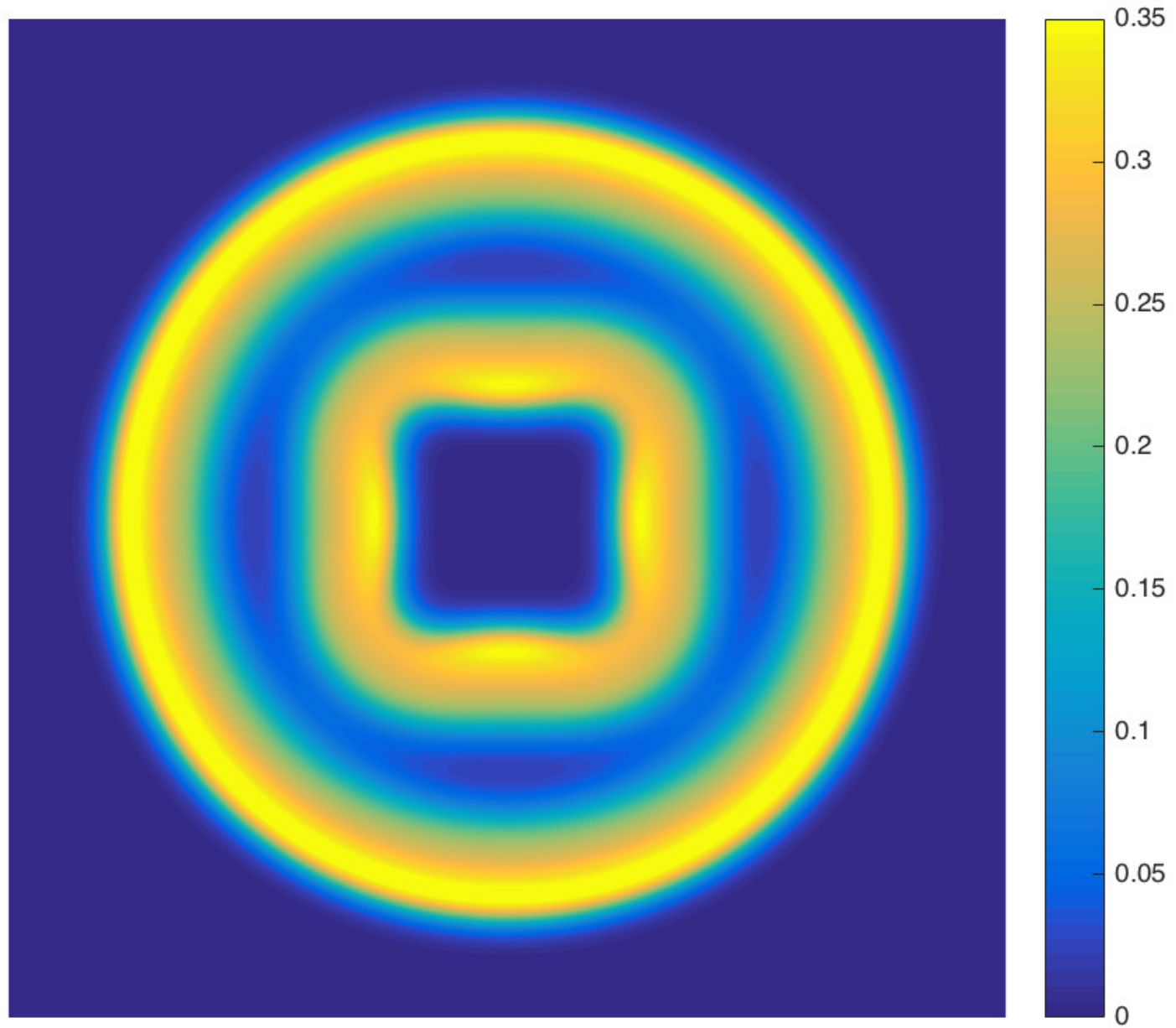}
\end{minipage}}
 \subfloat[$t=2000$]{
\begin{minipage}{.17\textwidth}\centering
\includegraphics[width=1.0\textwidth]{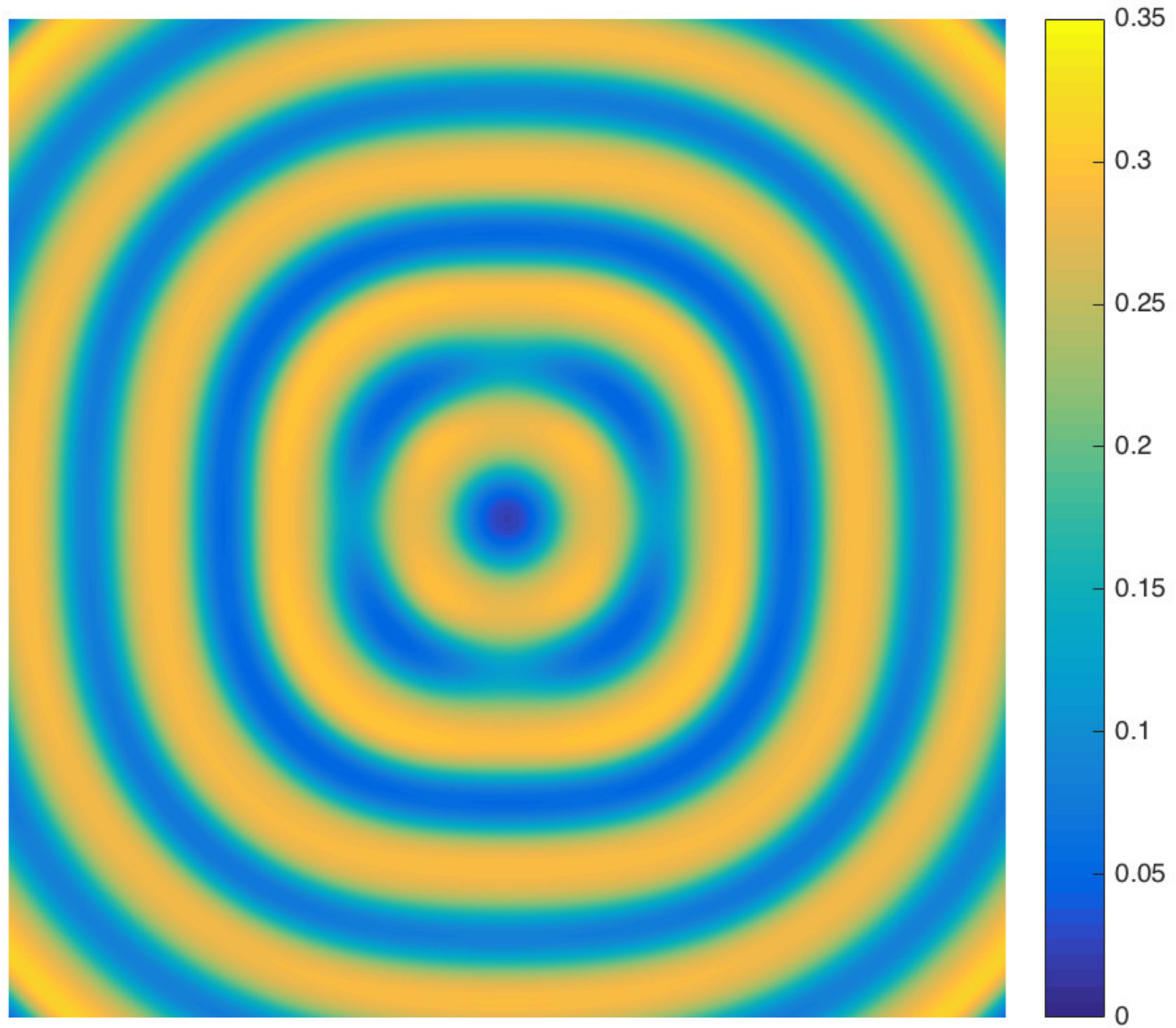}
\end{minipage}}
 \subfloat[$t=15000$]{
\begin{minipage}{.17\textwidth}\centering
\includegraphics[width=1.0\textwidth]{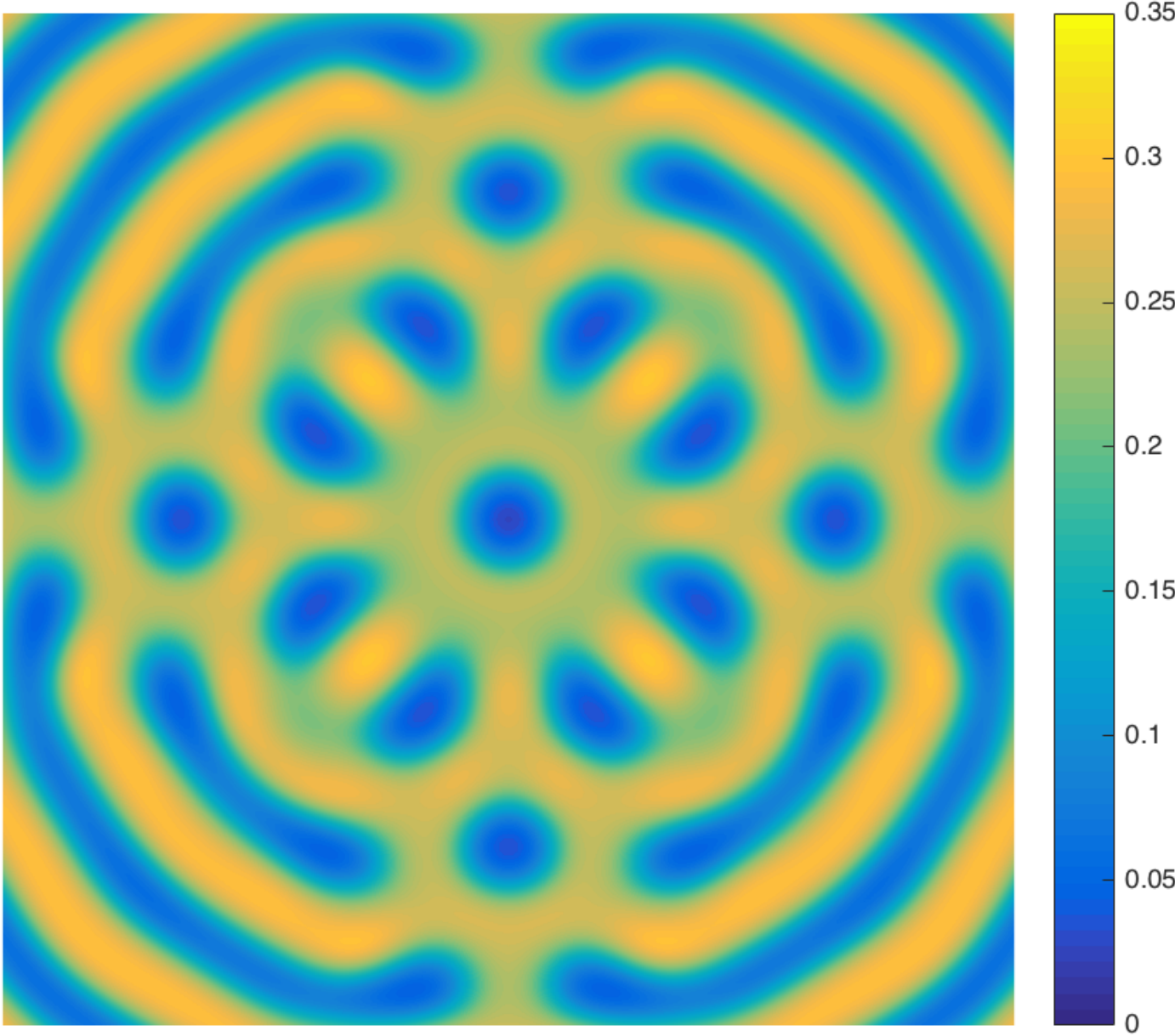}\\
\end{minipage}}
 \subfloat[$t=30000$]{
\begin{minipage}{.17\textwidth}\centering
\includegraphics[width=1.0\textwidth]{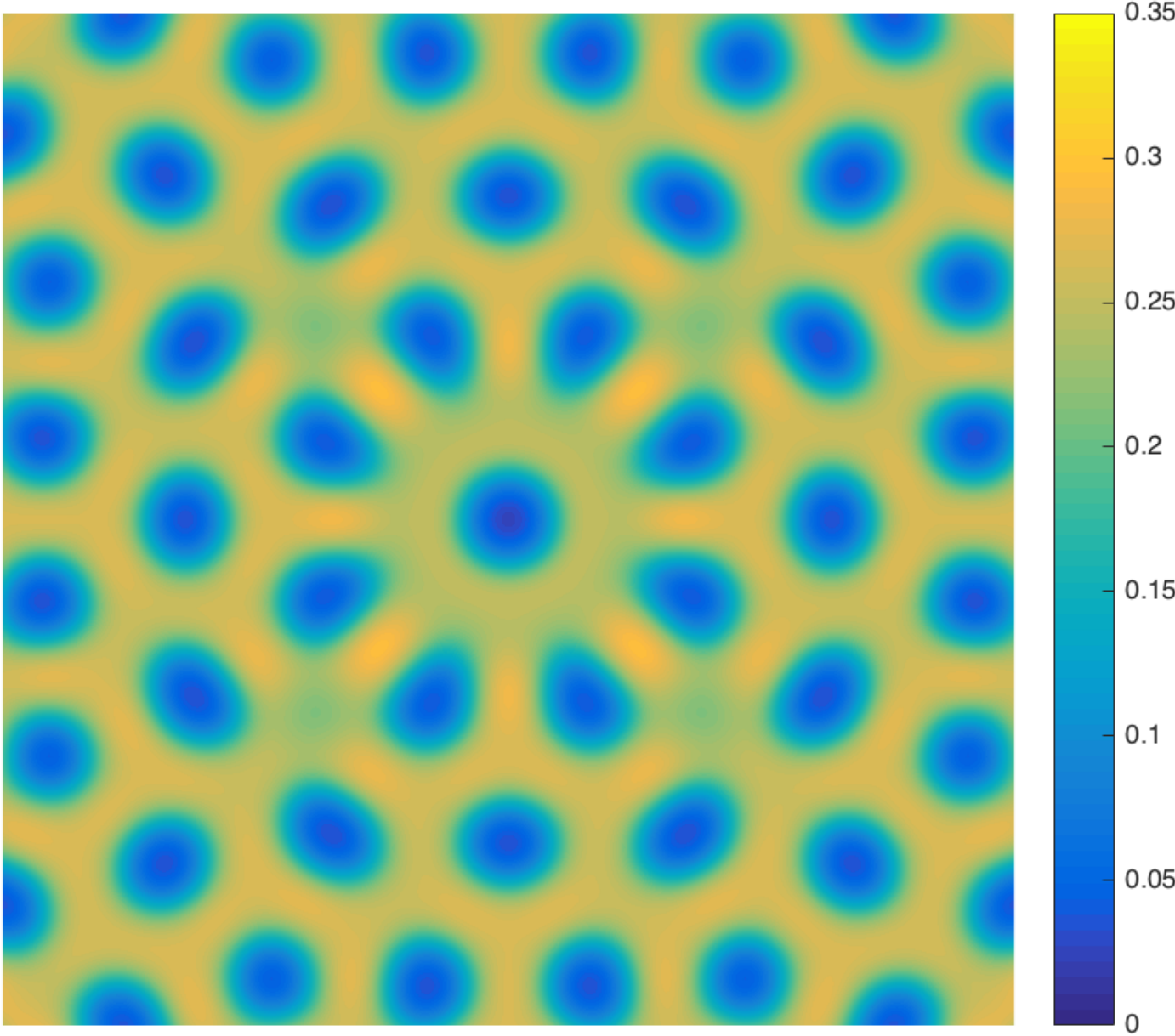}\\
\end{minipage}}}
{\subfloat[$t=400$]{
\begin{minipage}{.17\textwidth}\centering
\includegraphics[width=1.0\textwidth]{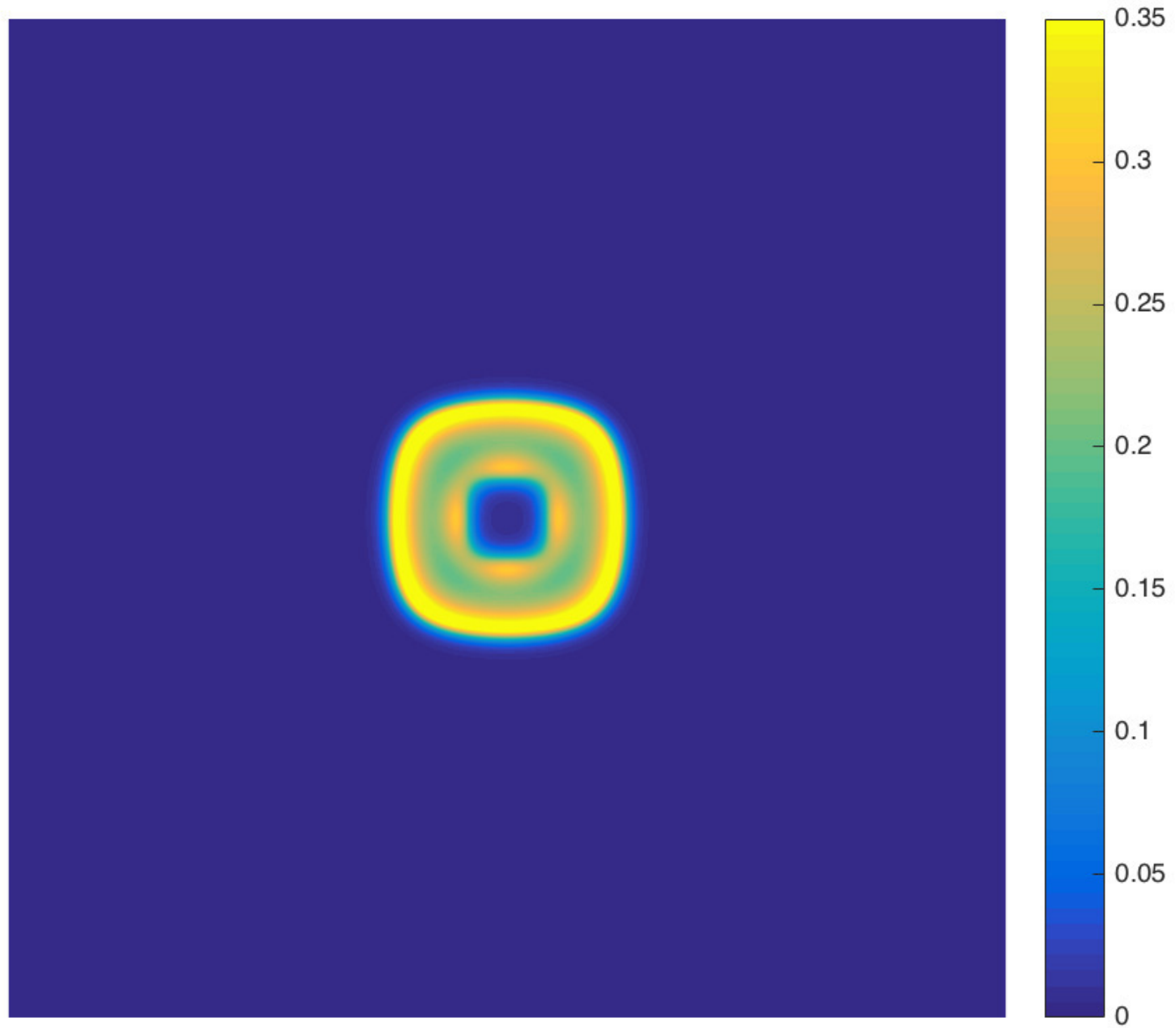}
\end{minipage}}
 \subfloat[$t=1600$]{
\begin{minipage}{.17\textwidth}\centering
\includegraphics[width=1.0\textwidth]{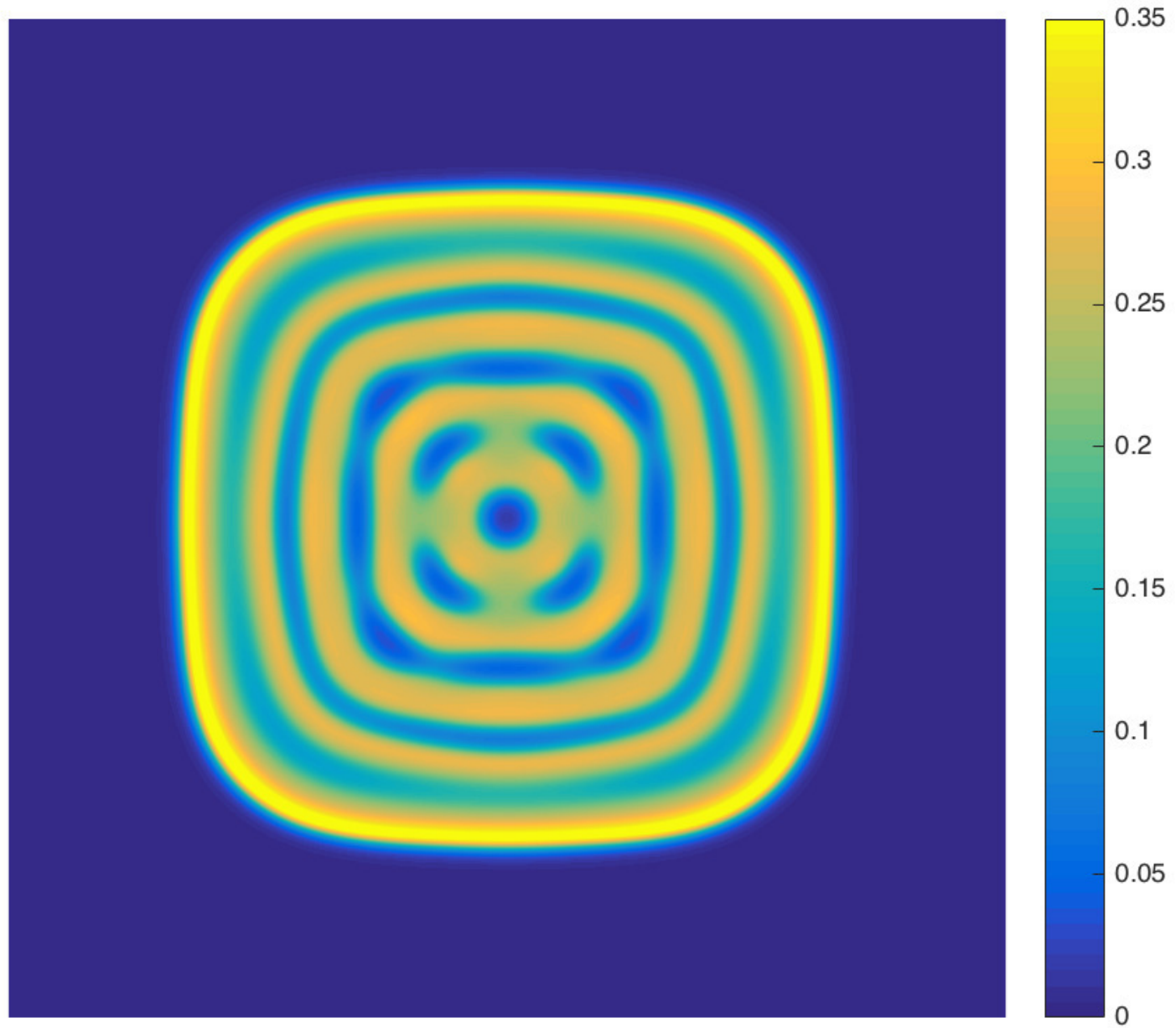}
\end{minipage}}
 \subfloat[$t=5000$]{
\begin{minipage}{.17\textwidth}\centering
\includegraphics[width=1.0\textwidth]{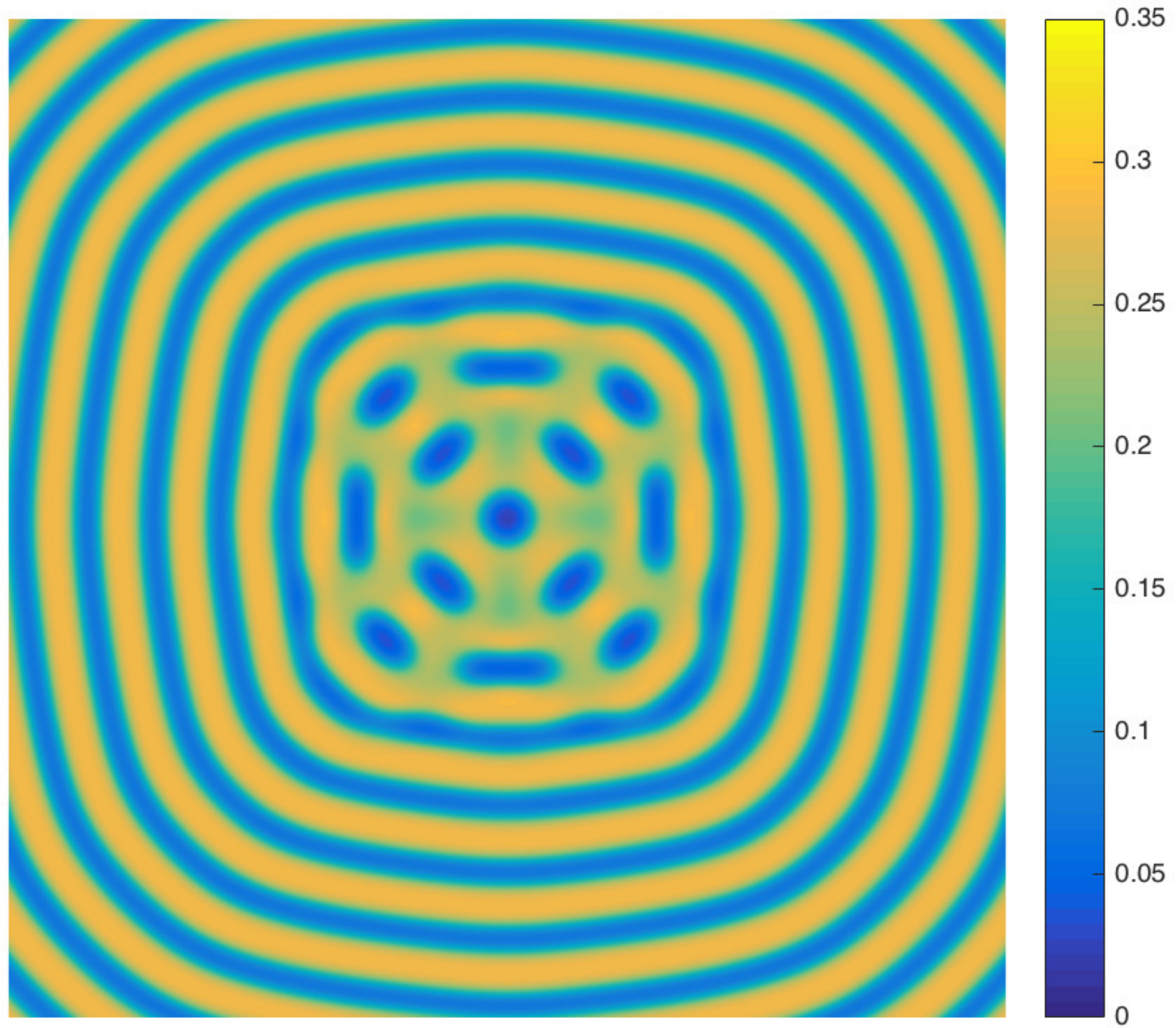}
\end{minipage}}
 \subfloat[$t=10000$]{
\begin{minipage}{.17\textwidth}\centering
\includegraphics[width=1.0\textwidth]{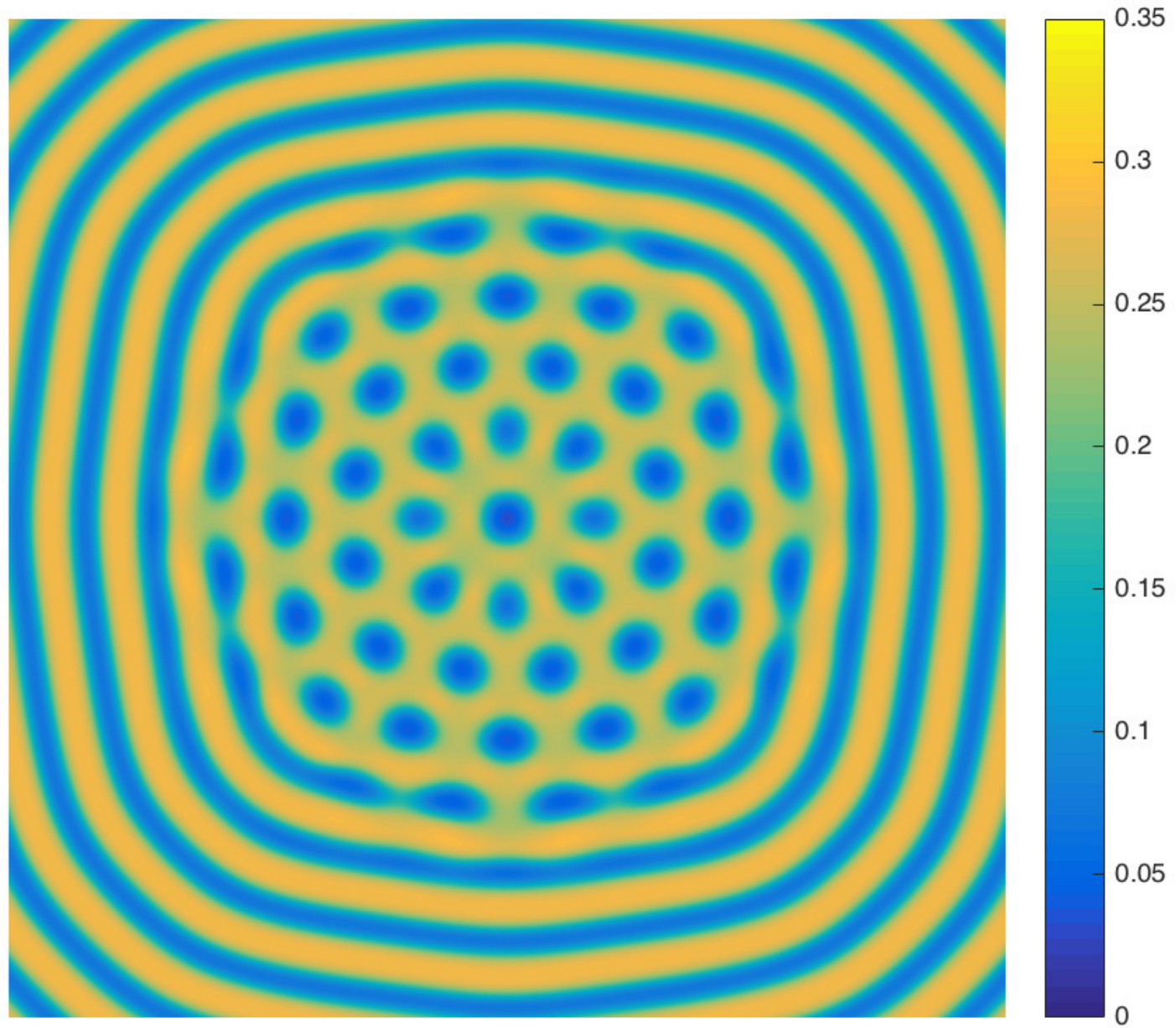}\\
\end{minipage}}
 \subfloat[$t=30000$]{
\begin{minipage}{.17\textwidth}\centering
\includegraphics[width=1.0\textwidth]{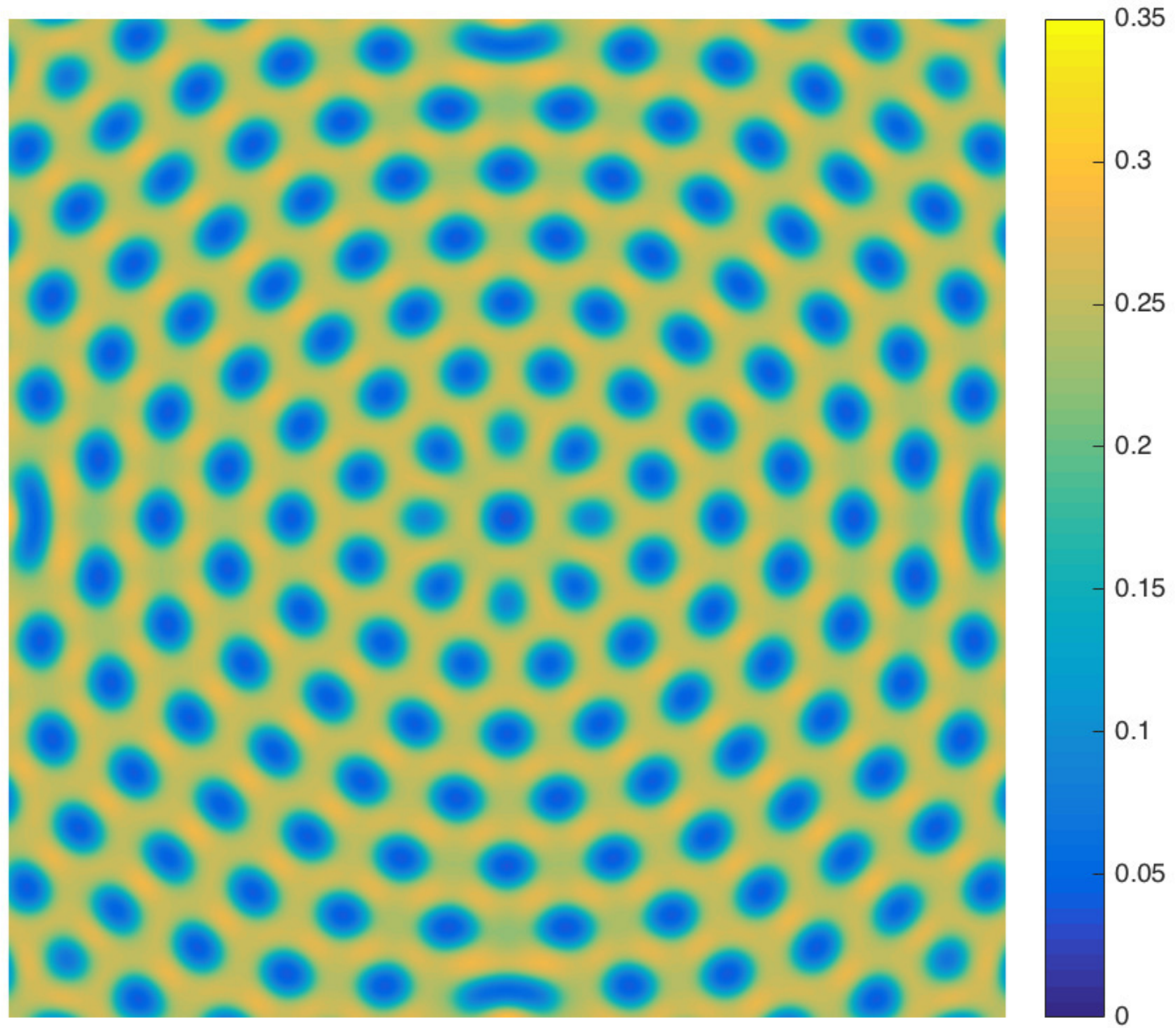}\\
\end{minipage}}}
{\subfloat[$t=800$]{
\begin{minipage}{.17\textwidth}\centering
\includegraphics[width=1.0\textwidth]{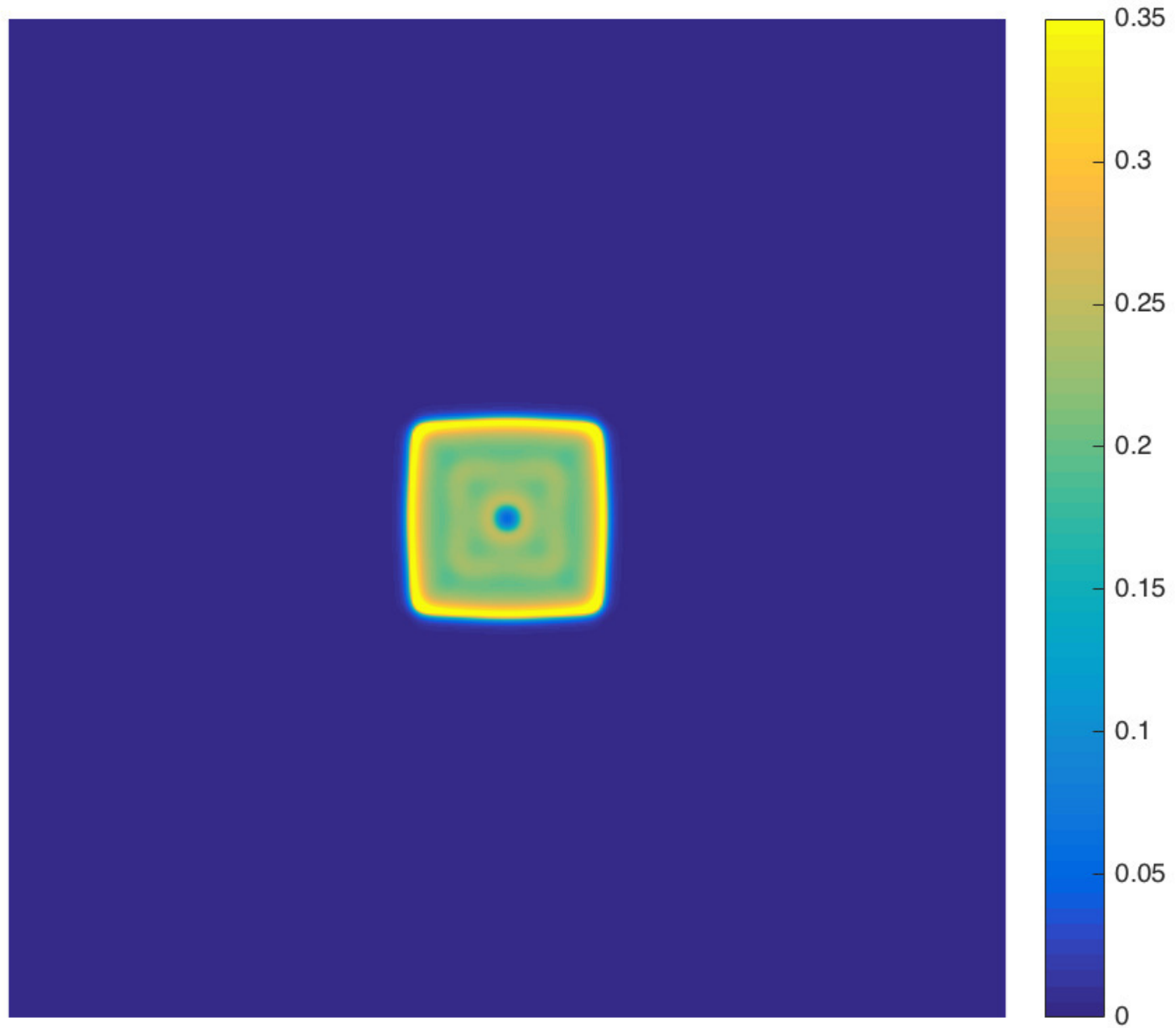}
\end{minipage}}
 \subfloat[$t=2600$]{
\begin{minipage}{.17\textwidth}\centering
\includegraphics[width=1.0\textwidth]{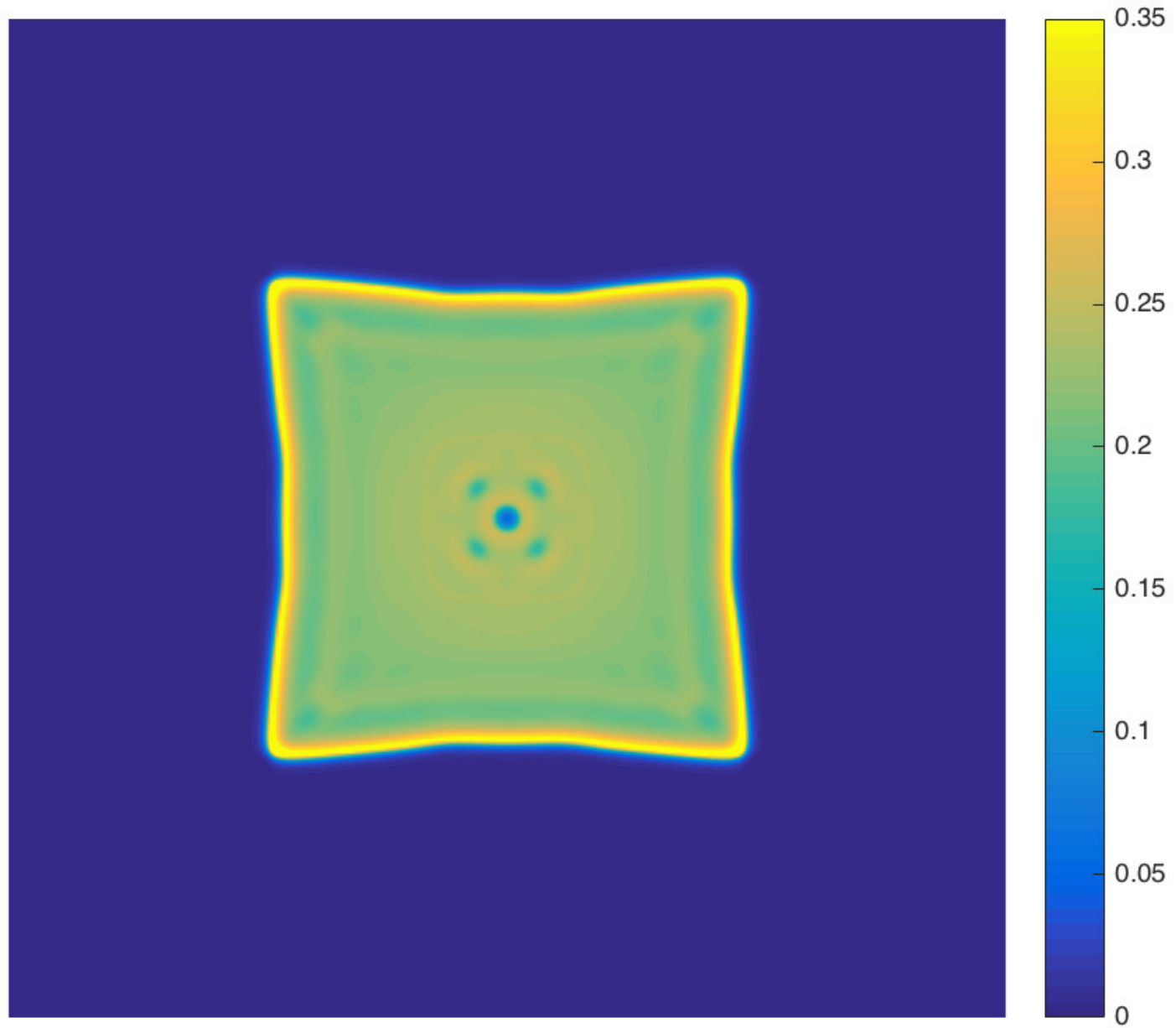}
\end{minipage}}
 \subfloat[$t=4000$]{
\begin{minipage}{.17\textwidth}\centering
\includegraphics[width=1.0\textwidth]{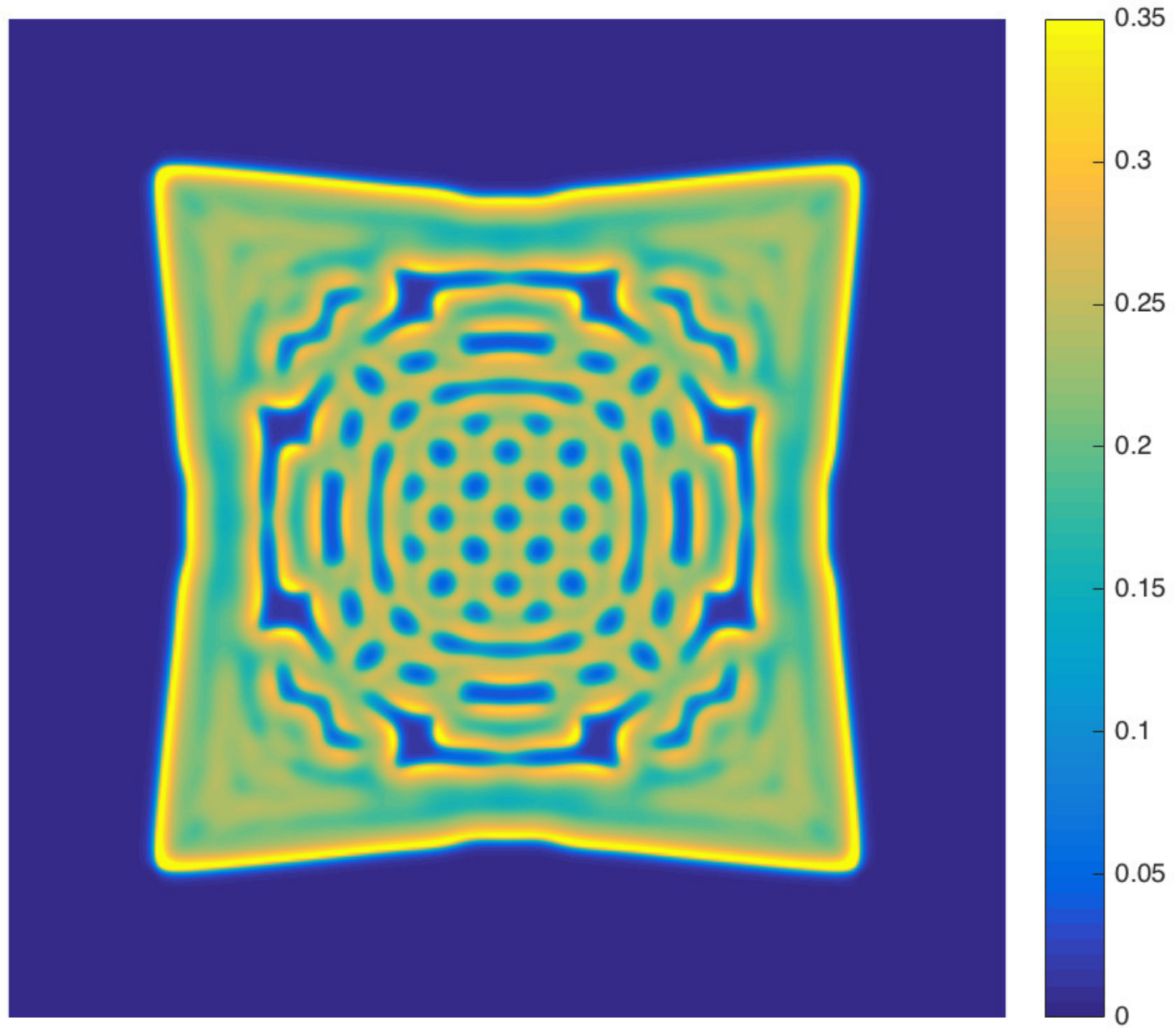}
\end{minipage}}
 \subfloat[$t=5800$]{
\begin{minipage}{.17\textwidth}\centering
\includegraphics[width=1.0\textwidth]{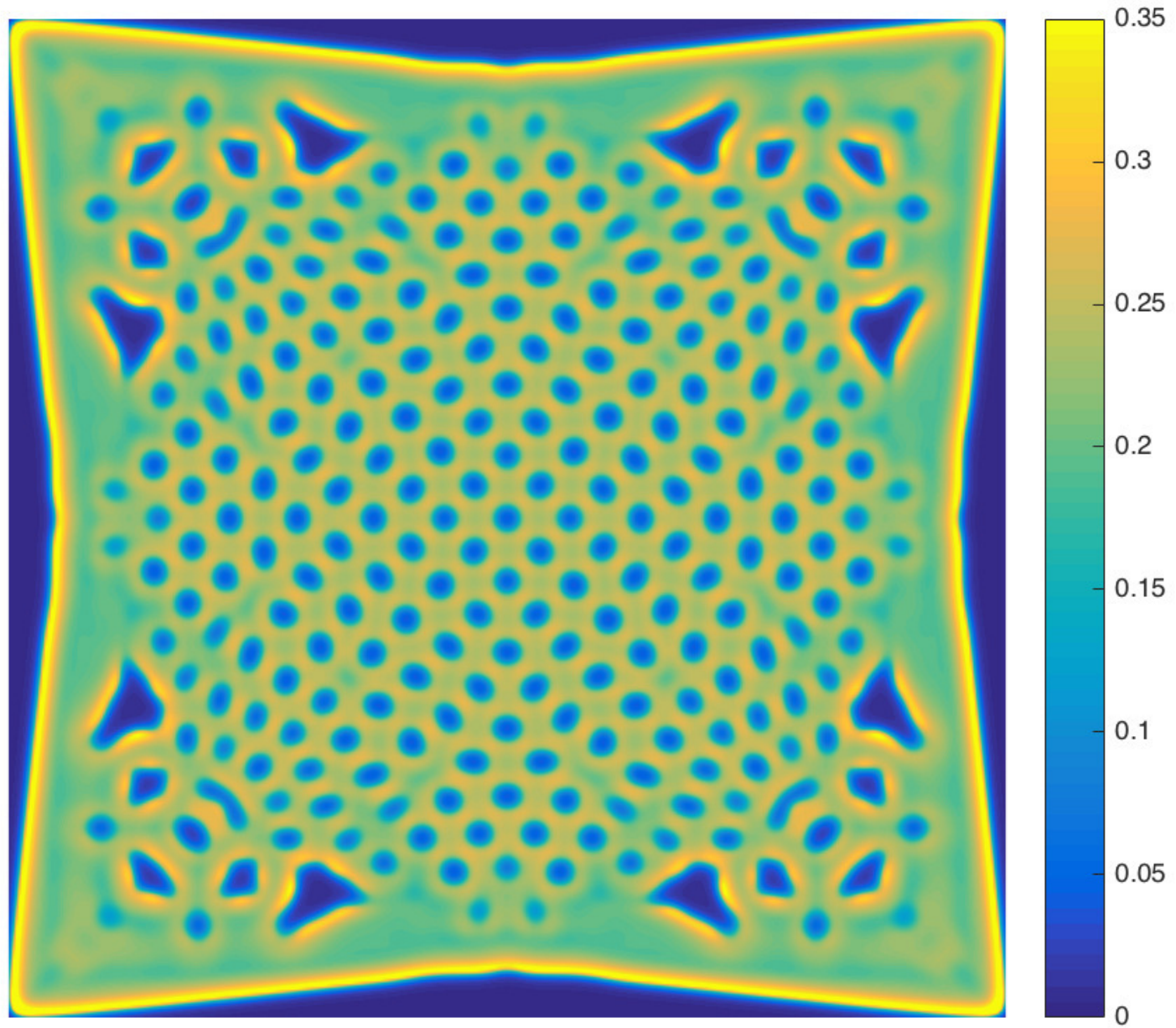}\\
\end{minipage}}
 \subfloat[$t=30000$]{
\begin{minipage}{.17\textwidth}\centering
\includegraphics[width=1.0\textwidth]{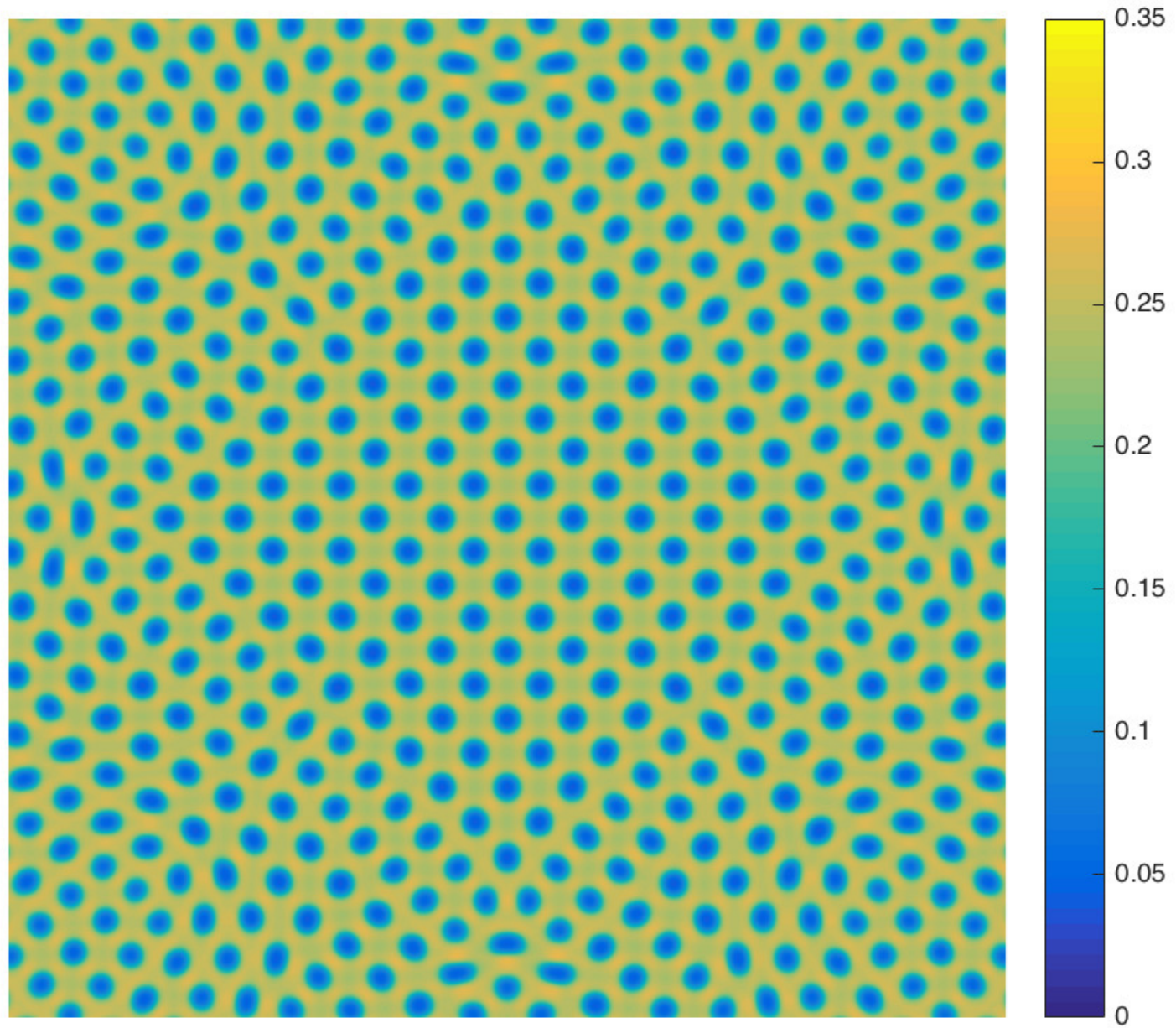}\\
\end{minipage}}}
\caption{Evolution of the solution $v$ in Eq. (\ref{s1:e1}) with $\kappa=0.055$: (a)-(e) $v$ contours are shown at different times from left to right for fractional order $\alpha=2.0$; (f)-(j) $v$ contours are shown at different times from left to right for fractional order $\alpha=1.7$; (k)-(o) $v$ contours are shown at different times from left to right for fractional order $\alpha=1.5$.}
\label{Fig. 4}
\end{figure}

Figs. $\ref{Fig. 3}$, $\ref{Fig. 4}$ show the evolutions of the numerical solution $v$ and summarize the effects of the super-diffusion for the fractional GS model. All the pictures in the two figures are snapshots from the numerical solutions $v$ in the domain $(0,1)^2$. Since the fractional order affects the speed of the diffusion, the speeds of pattern formation are different for different $\alpha$.

In Fig. $\ref{Fig. 3}$ with $\kappa=0.063$, under the influence of the standard diffusion $(\alpha=2.0)$, the GS model exhibits patterns of mitosis. However, in the super-diffusion case  $\alpha=1.7$, the replication pattern has completely changed. Furthermore, when the fractional order $\alpha=1.5$, the patterns present different behavior. We observe that the structure of patterns and the size of the spots are different. In Fig. $\ref{Fig. 4}$ with $\kappa=0.055$, the GS model with the normal diffusion $(\alpha=2.0)$ produces a circular wave propagating outward and form a structured pattern shown in picture (e). In addition, the reduction of the fractional order $\alpha=1.7$ affects the size of patterns with smaller spots. For smaller fractional order $(\alpha=1.5)$, we observe a new process of pattern formation. The process propagates outward until the whole area reaches the final steady state.

\subsection{Comparison of simulations with spectral collocation method}

In order to verify the accuracy of the numerical results obtained by the previously mentioned numerical method, we also use the spectral collocation method \cite{zeng2017generalized} in space discretization to solve the fractional GS model. In this simulation, the domain $\Omega$, parameters $\kappa, F$ and initial-boundary condition are the same as the aforementioned simulation in subsection \ref{gss}. We choose the collocation points $N=2000$ for the spectral collocation method, $h=\frac{1}{2000}$ for the difference scheme and time step $\tau=0.1$ for both methods. In addition, we also compute the RDFs $g(r)$ for steady spot patterns obtained by the different numerical methods to quantify the averaged distance between spot pairs.

The RDFs $g(r)$ describes the particle density which varies as a function of the distance from a reference particle. It is useful to measure the probability of finding a particle at a distance $r$ away from the given reference particle. In this paper, we calculate the distances between the center of the reference spots and that of the other spots. We then bin them into a histogram and normalize them with respect to the total spot numbers. In $g(r)$ function, the first peak value and its corresponding distance $r$ indicates the type of patterns and the fractional order effectively. Consequently, we not only can compare the steady patterns obtained by the aforementioned two numerical methods, but also can compare the corresponding RDFs to further guarantee the coherence of the two numerical results.

\begin{figure}[htpb]
\centering
 \subfloat[$\alpha=2.0, \kappa_{1}=0.063$.]{
\begin{minipage}{.31\textwidth}\centering
\includegraphics[width=1.0\textwidth]{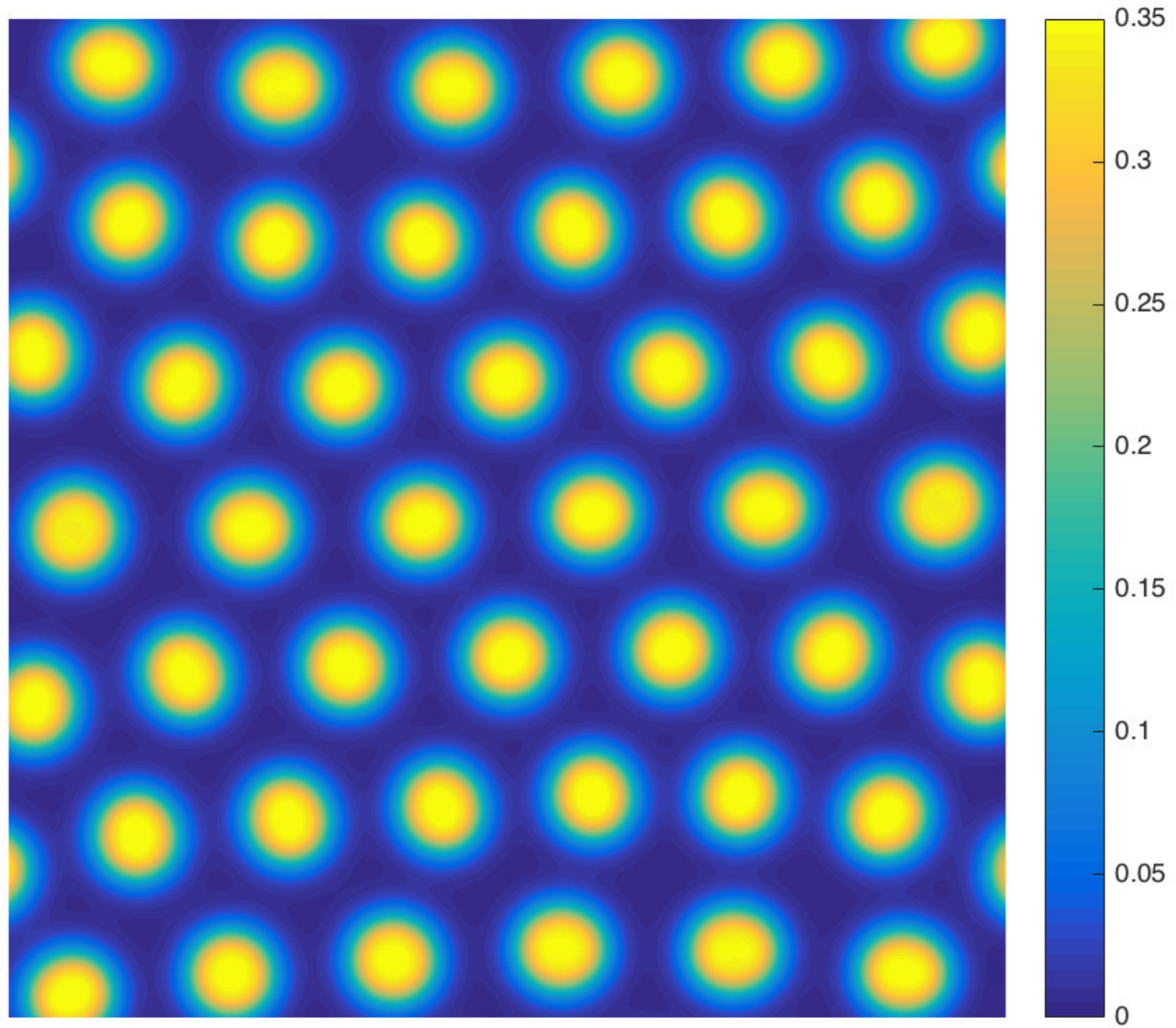}\\
\includegraphics[width=1.0\textwidth]{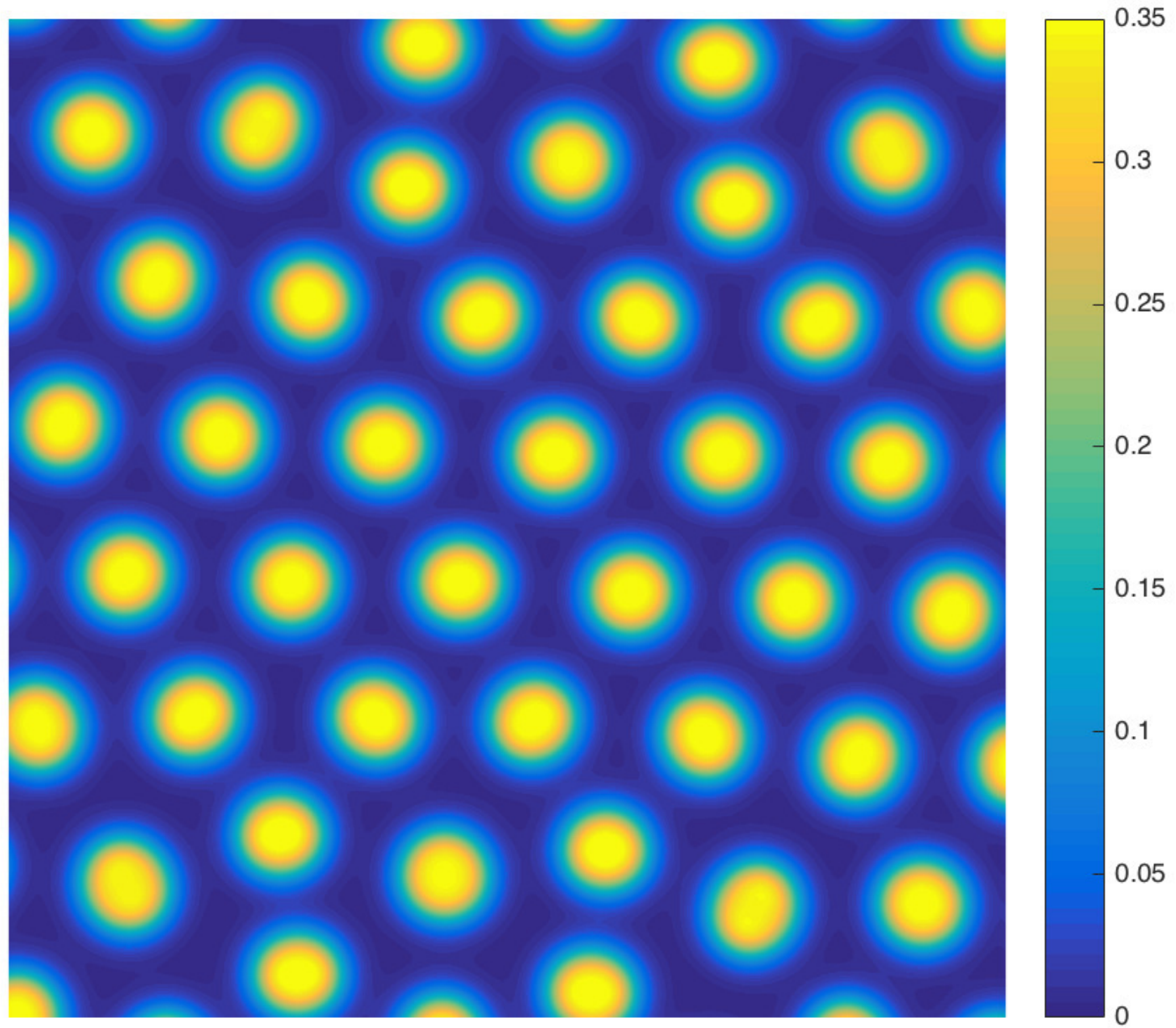}\\
\includegraphics[width=1.0\textwidth]{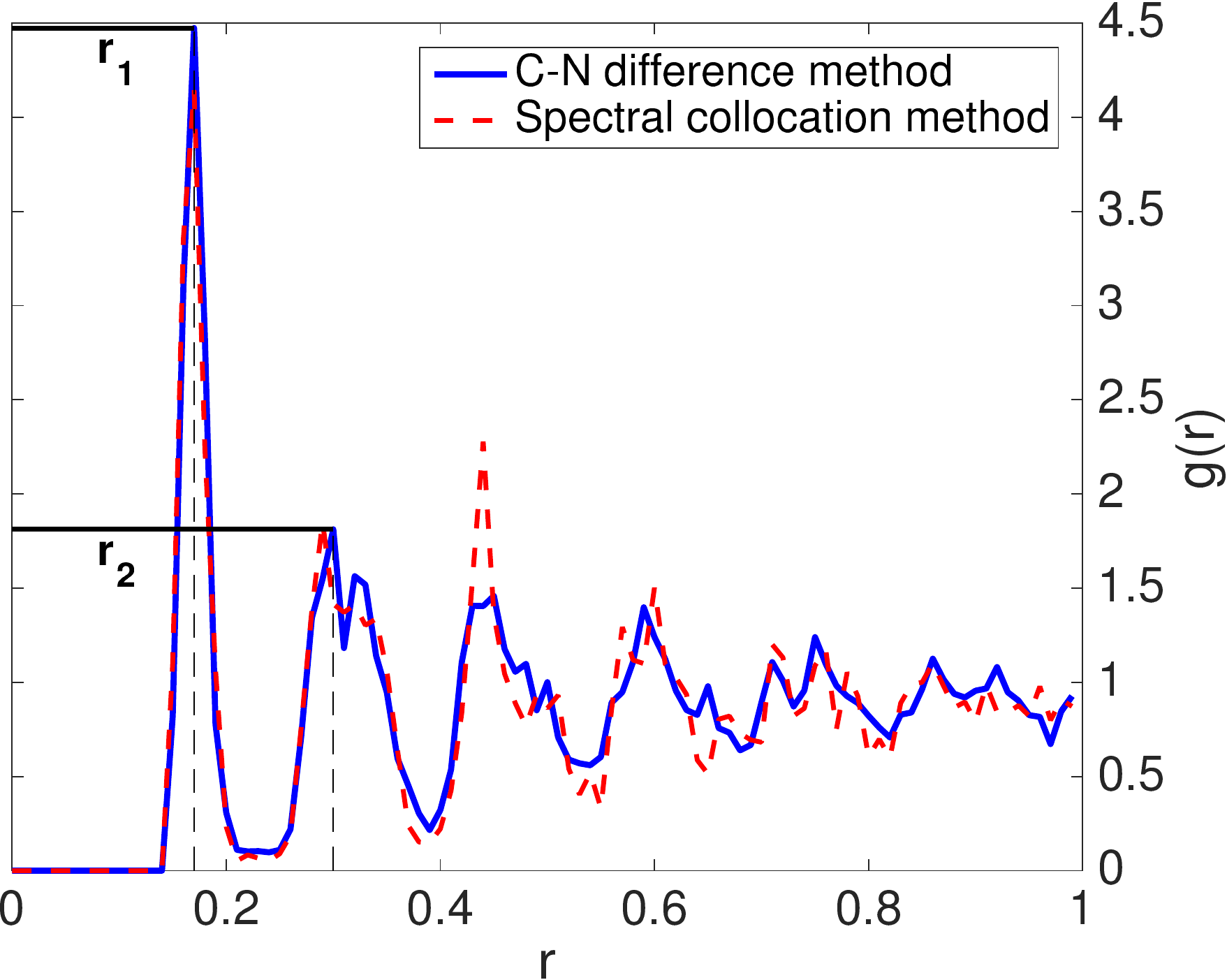}
\end{minipage}}
\subfloat[$\alpha=1.7, \kappa_{1}=0.063$.]{
\begin{minipage}{.31\textwidth}\centering
\includegraphics[width=1.0\textwidth]{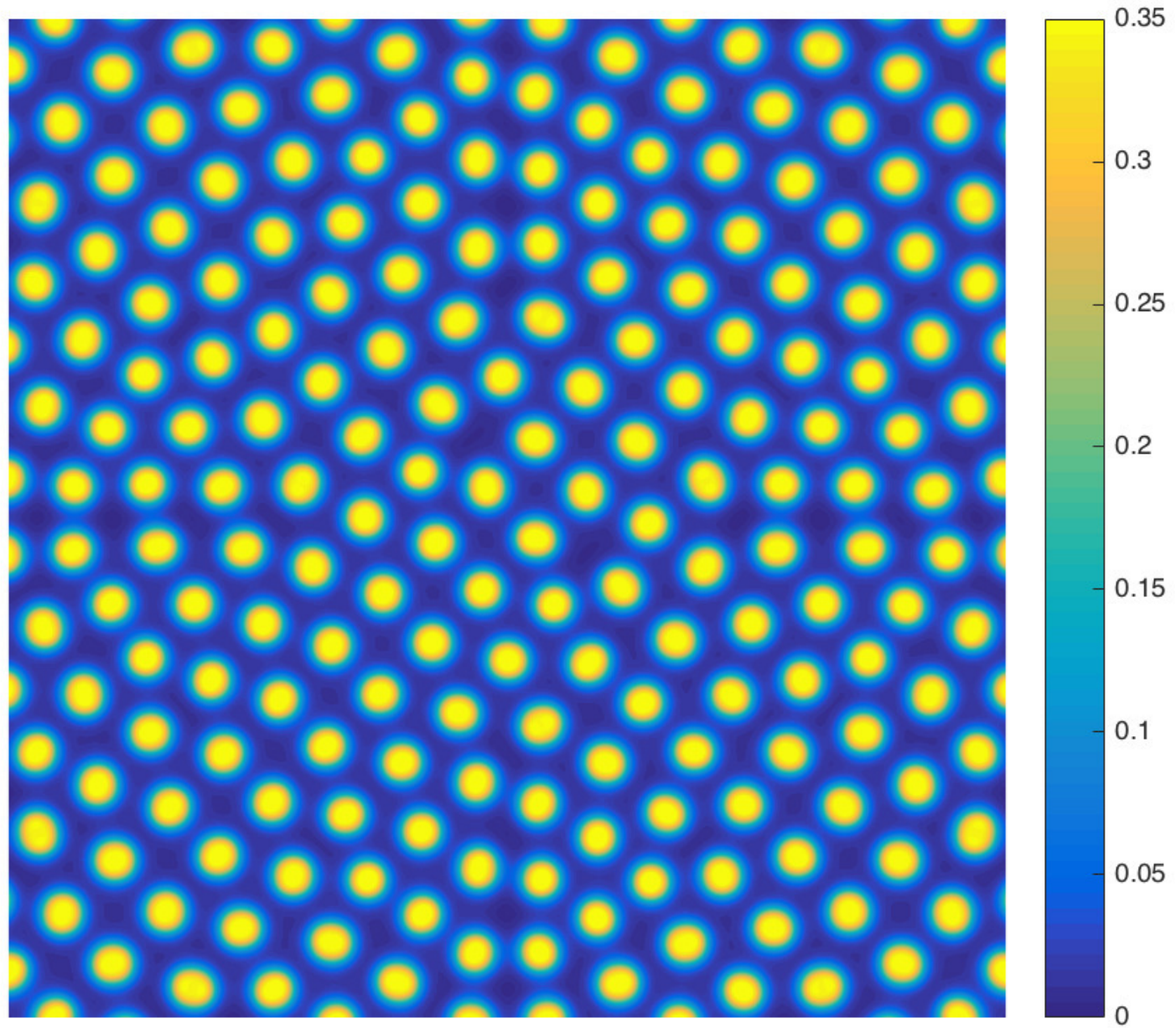}\\
\includegraphics[width=1.0\textwidth]{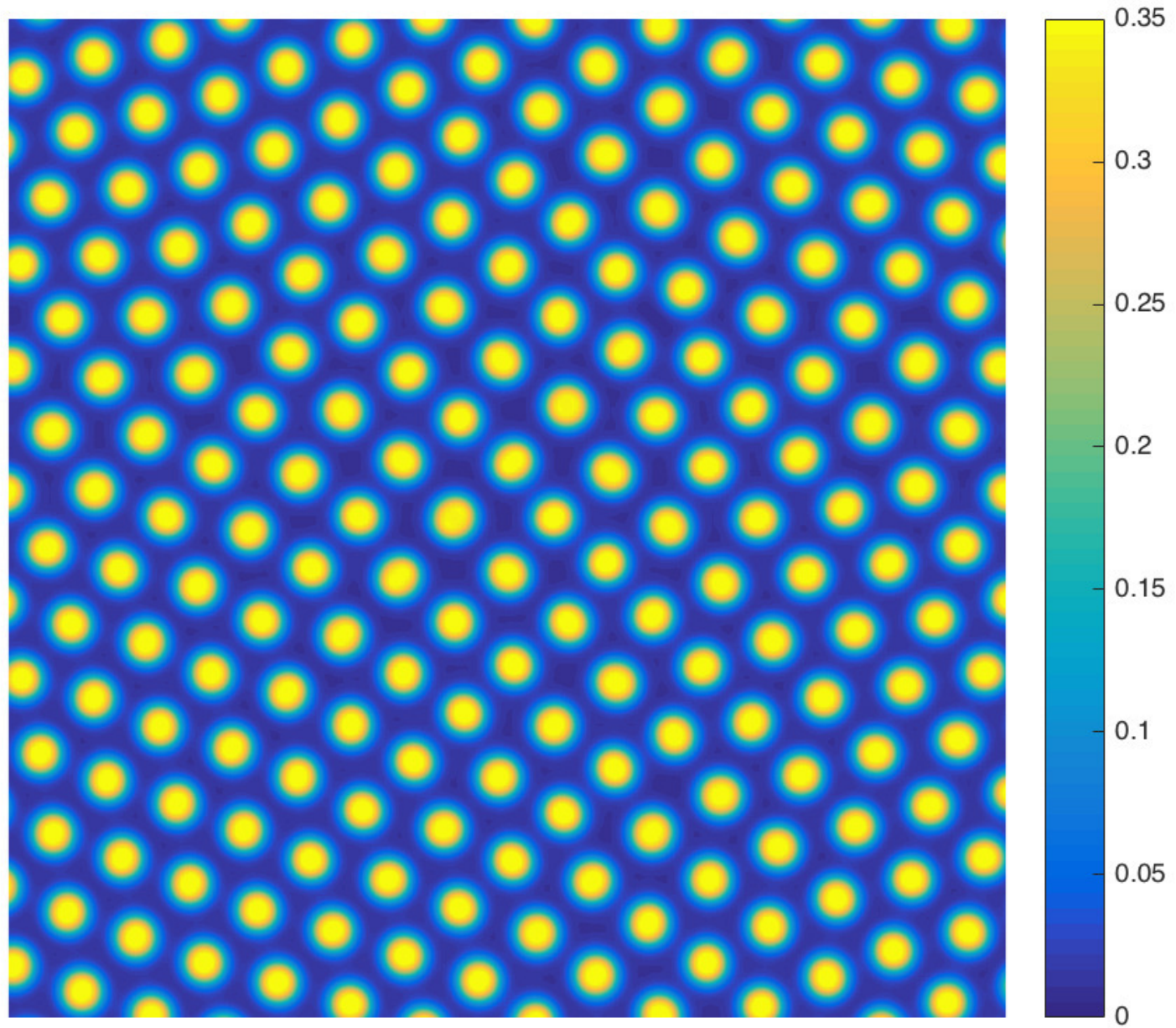}\\
\includegraphics[width=1.0\textwidth]{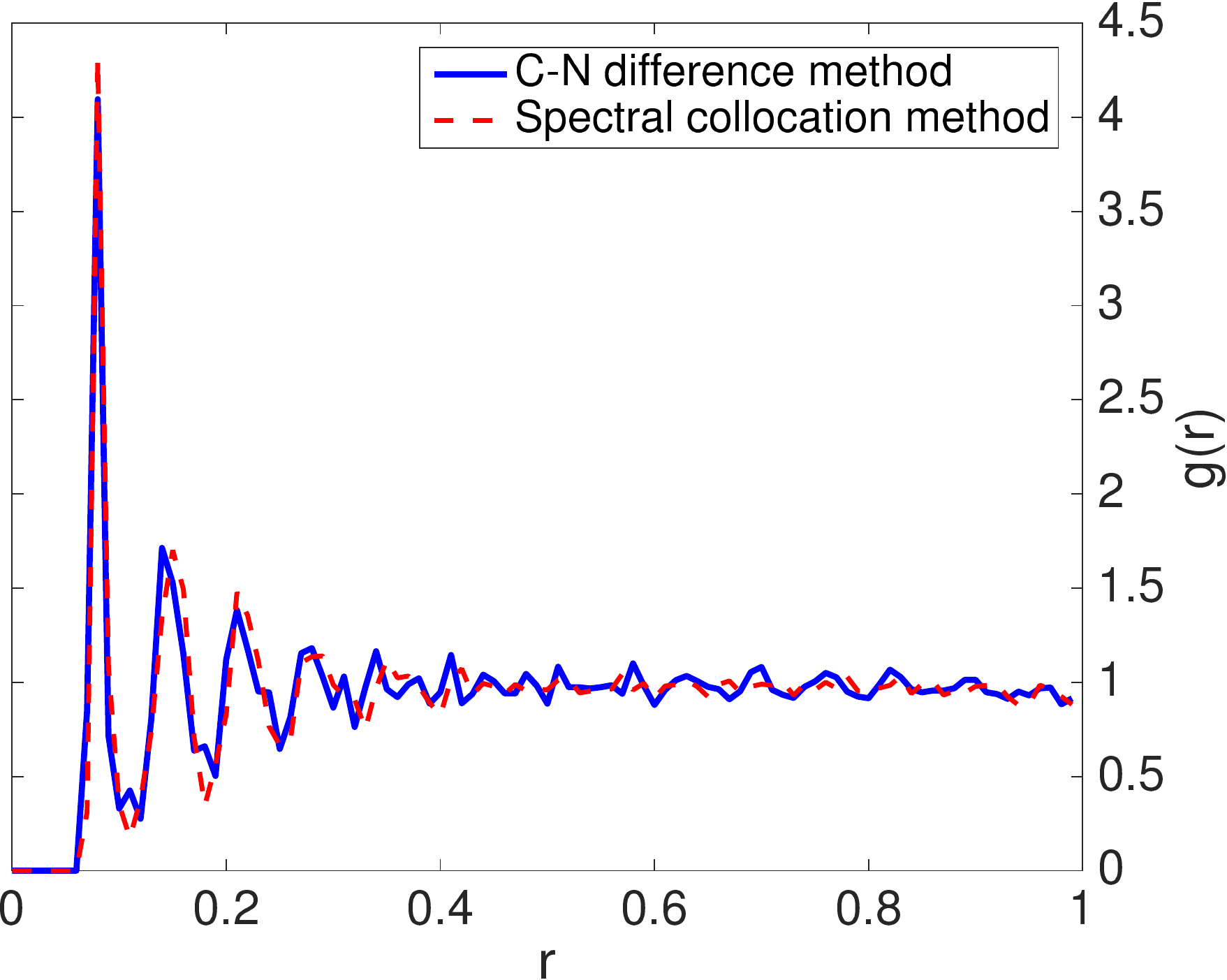}
\end{minipage}}
\subfloat[$\alpha=1.5, \kappa_{1}=0.063$.]{
\begin{minipage}{.31\textwidth}\centering
\includegraphics[width=1.0\textwidth]{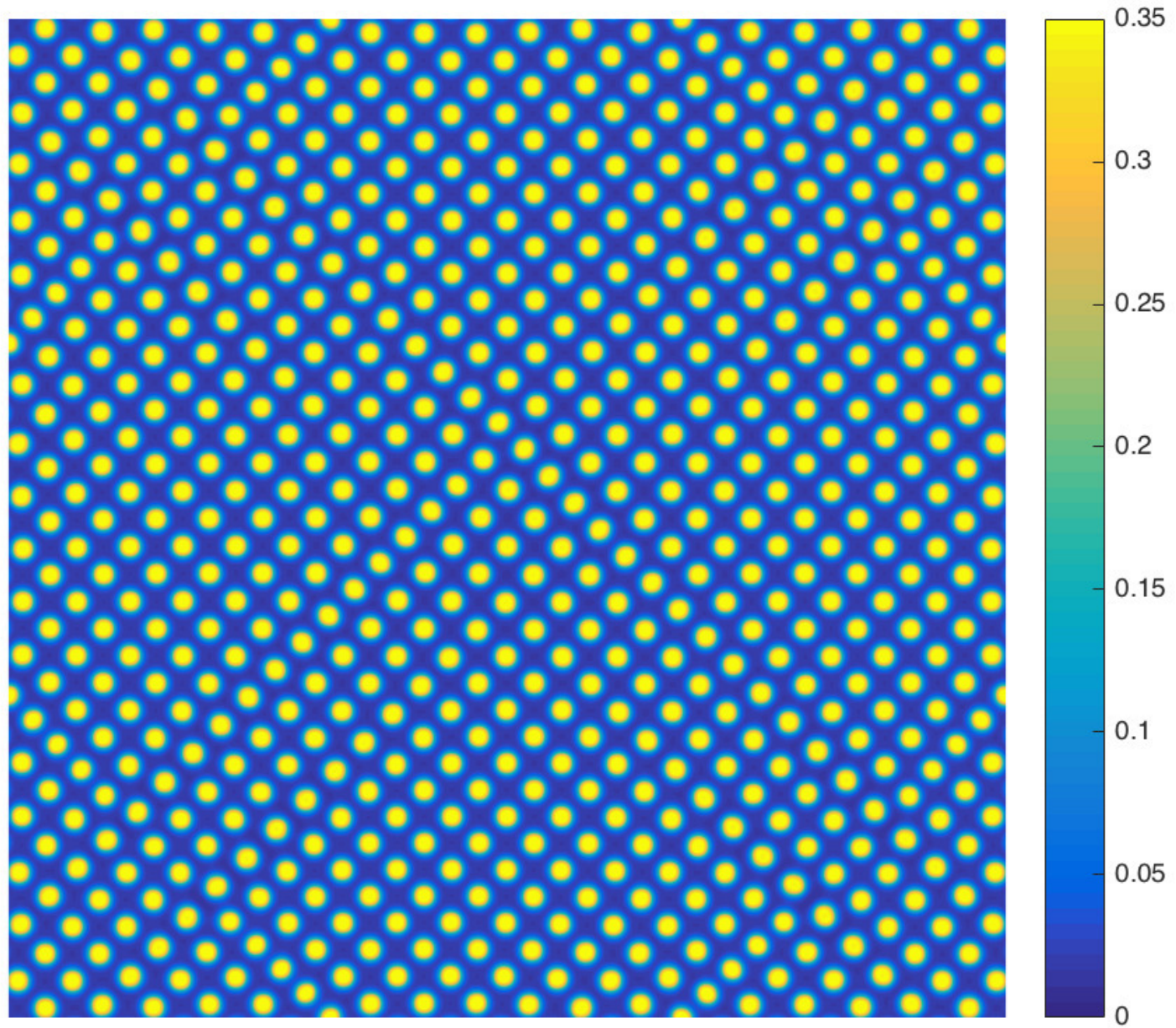}\\
\includegraphics[width=1.0\textwidth]{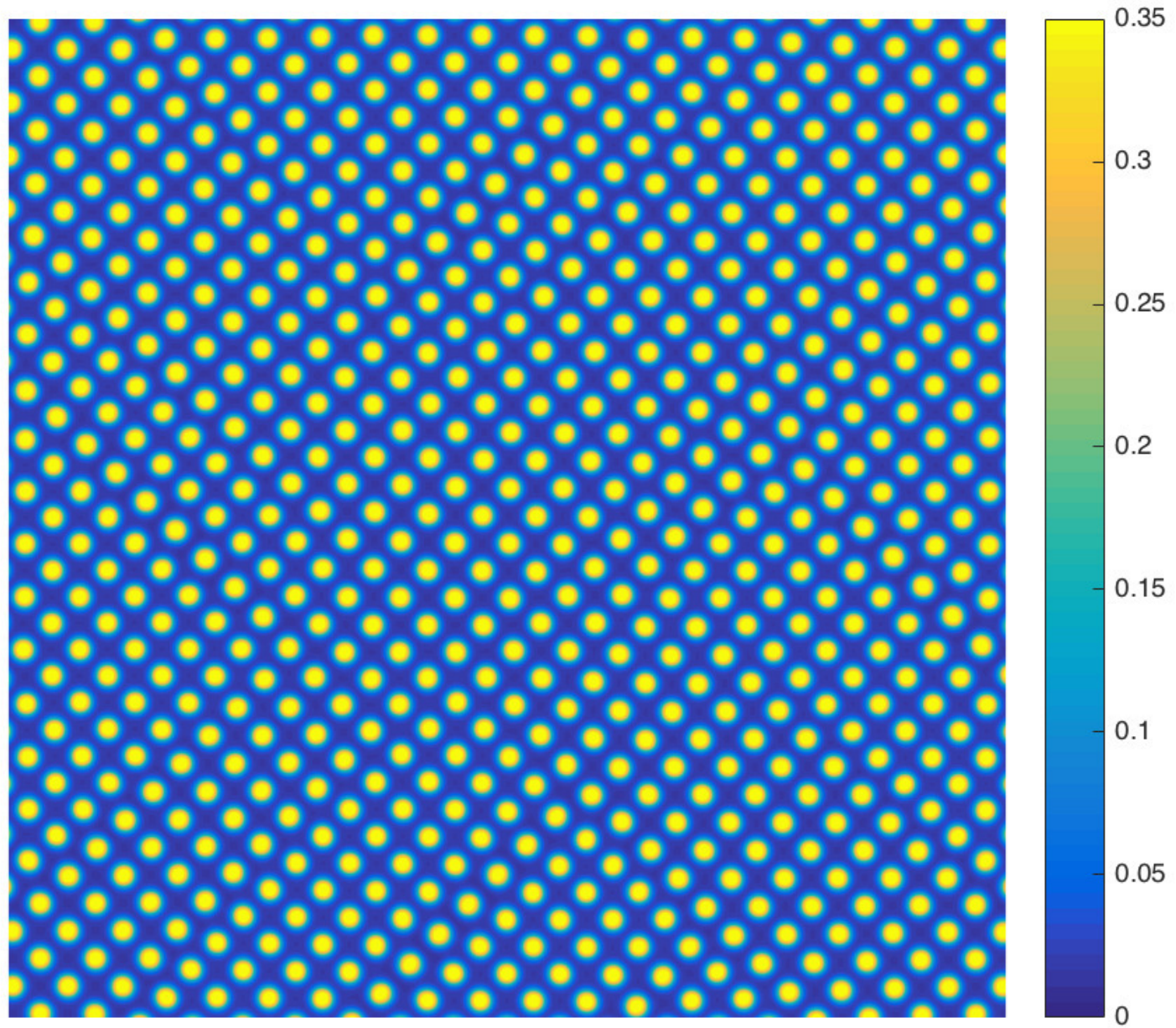}\\
\includegraphics[width=1.0\textwidth]{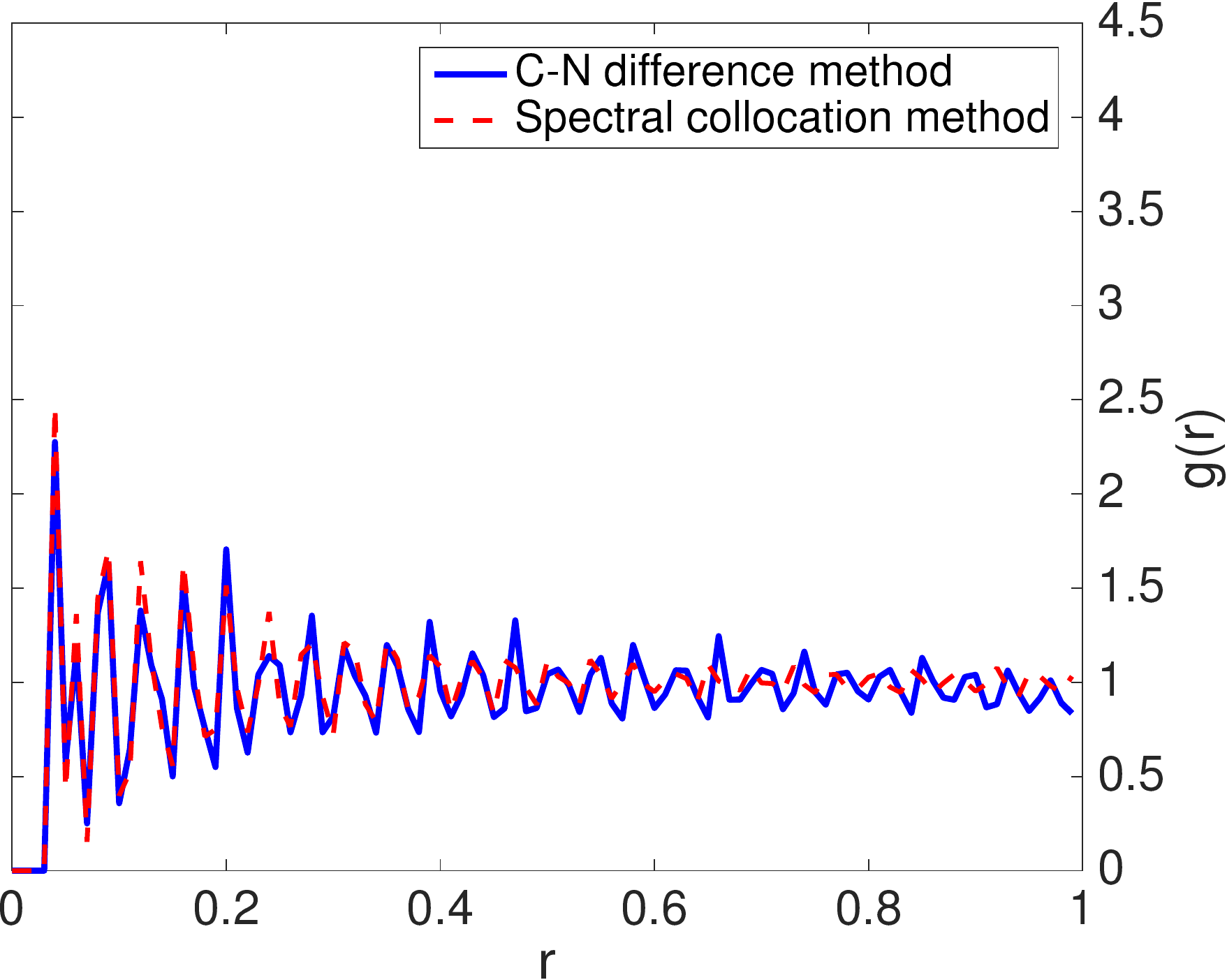}
\end{minipage}}
\caption{Steady patterns of the solution $v$ with $\kappa_{1}=0.063$ and their corresponding RDFs: $Top$: $v$ contours obtained by C-N difference scheme with weighted shifted Gr\"{u}nwald difference operators at steady states with different fractional orders. The fractional orders in (a)-(c) correspond to $\alpha=2.0, 1.7, 1.5$. $Middle$: $v$ contours obtained by spectral collocation method at steady states with different fractional orders.  The fractional orders in (a)-(c) correspond to $\alpha=2.0, 1.7, 1.5$. $Bottom$: RDFs $g(r)$: the blue line and the red dash line correspond to the variation of the spot density in steady spot  patterns obtained by C-N difference scheme and spectral collocation method with the same $\alpha$ and parameter $\kappa_{1}$, respectively. The different fractional orders in (a)-(c) correspond to $\alpha=2.0, 1.7, 1.5$.}
\label{Fig. a}
\end{figure}

\begin{figure}[htpb]
\centering
\subfloat[$\alpha=2.0, \kappa_{2}=0.055$.]{
\begin{minipage}{.31\textwidth}\centering
\includegraphics[width=1.0\textwidth]{d2055}\\
\includegraphics[width=1.0\textwidth]{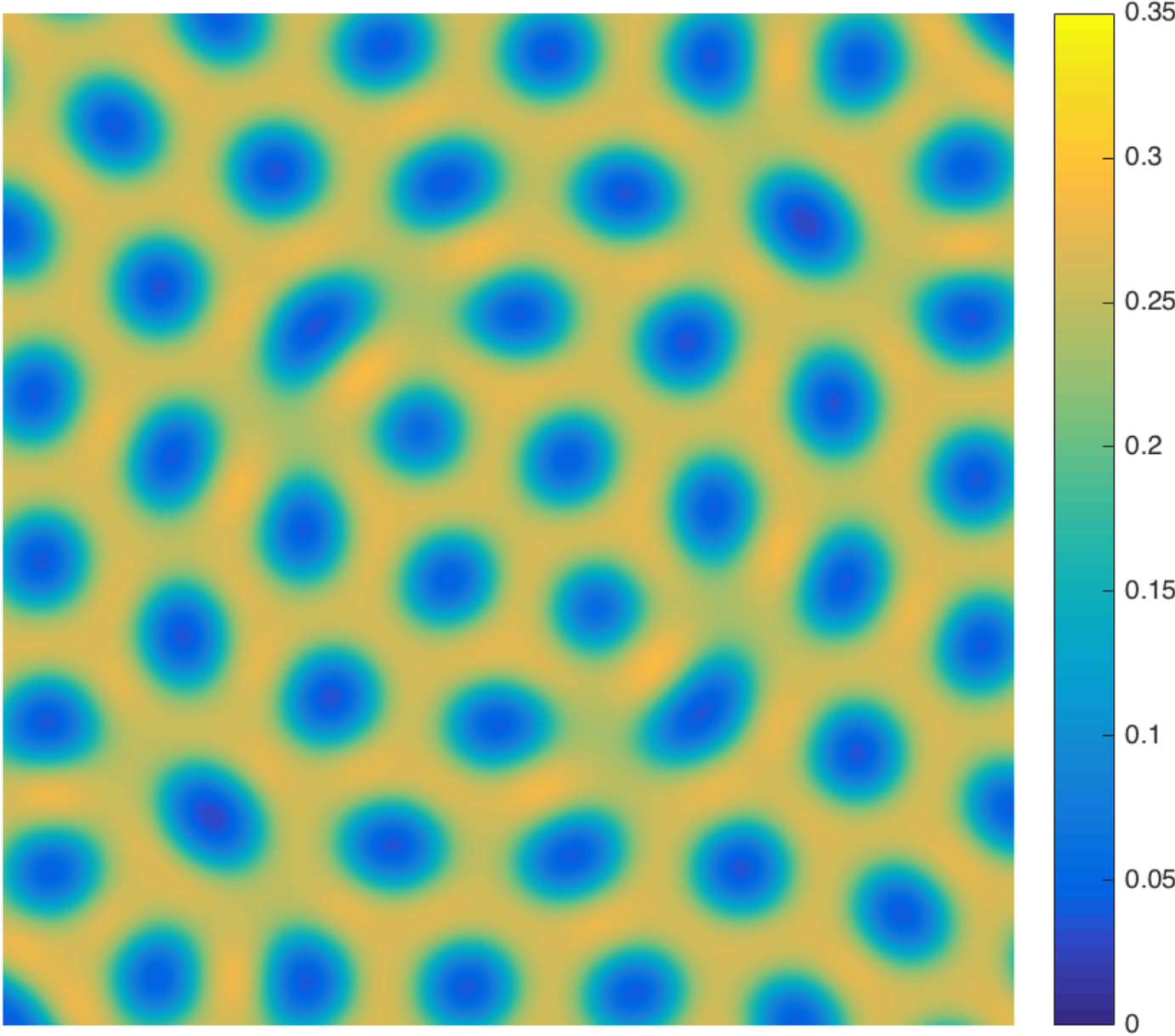}\\
\includegraphics[width=1.0\textwidth]{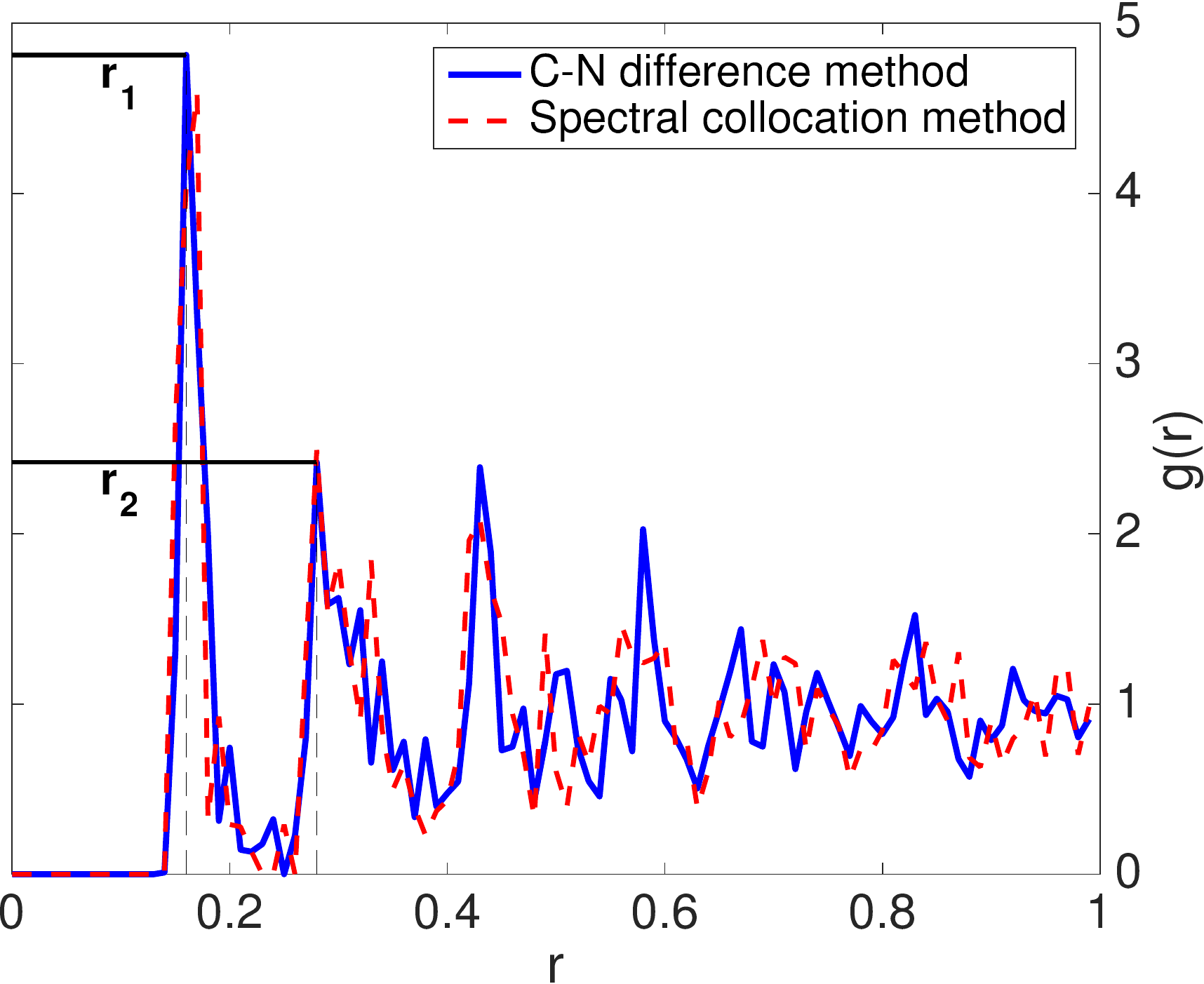}
\end{minipage}}
\subfloat[$\alpha=1.7, \kappa_{2}=0.055$.]{
\begin{minipage}{.31\textwidth}\centering
\includegraphics[width=1.0\textwidth]{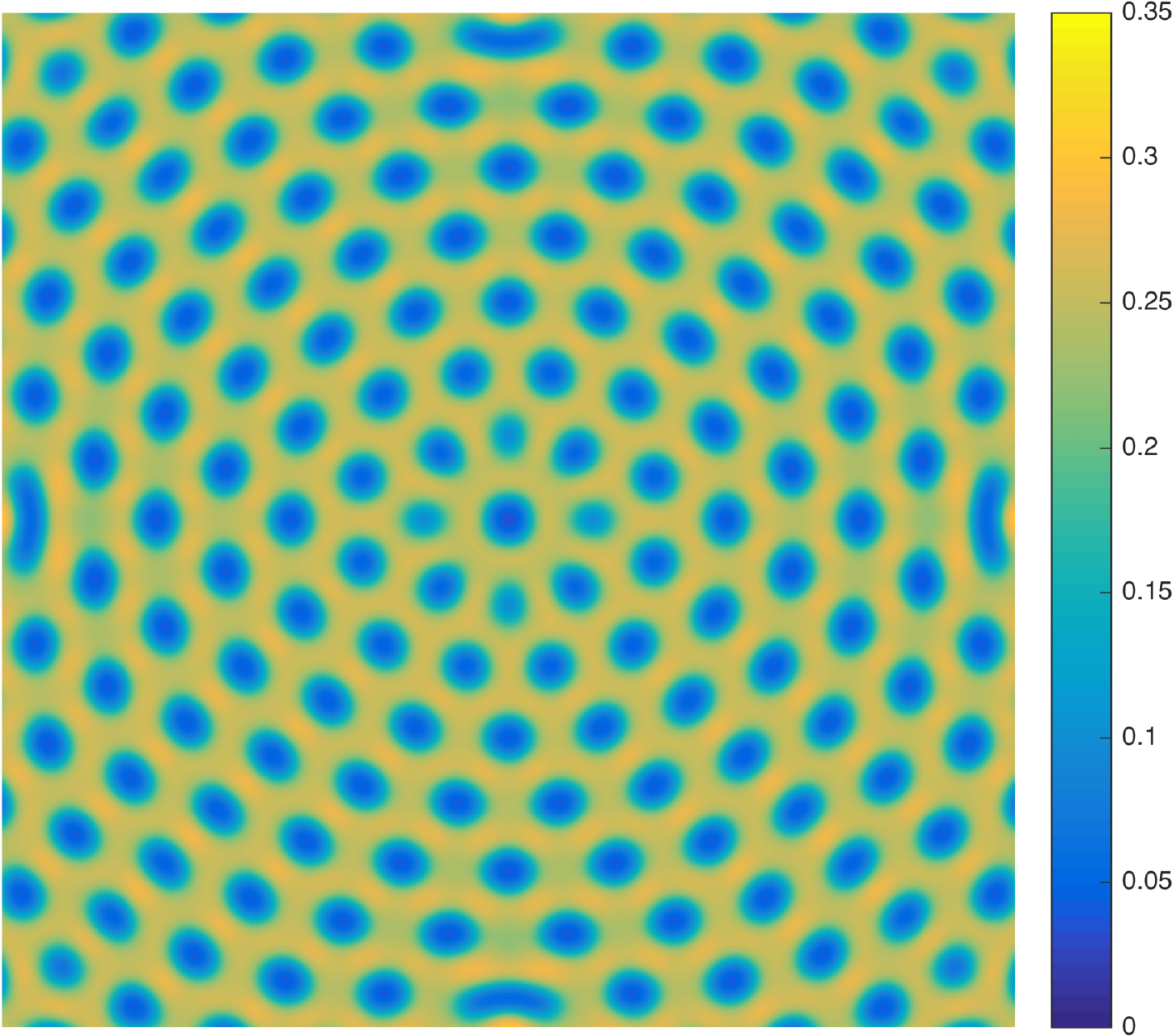}\\
\includegraphics[width=1.0\textwidth]{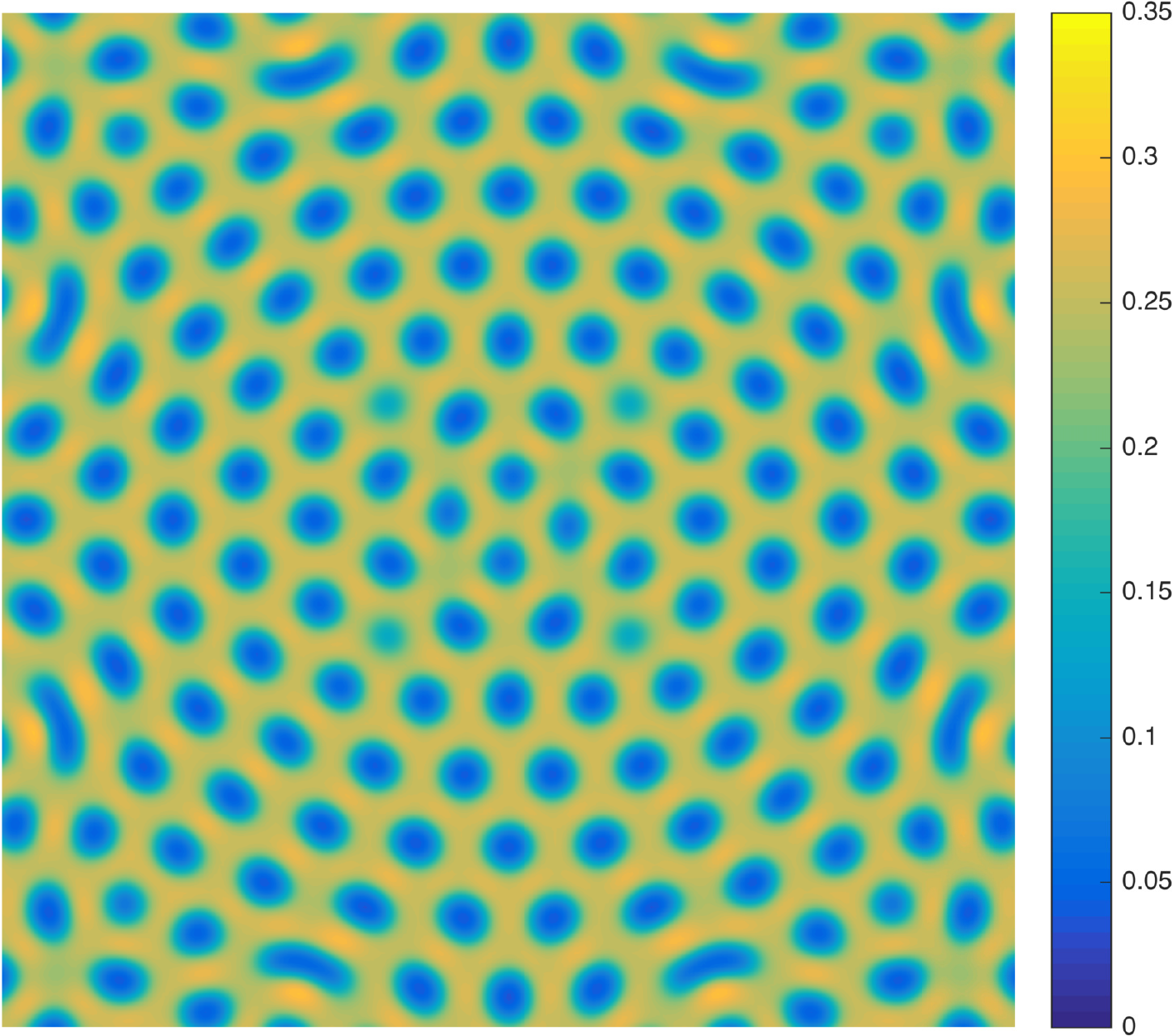}\\
\includegraphics[width=1.0\textwidth]{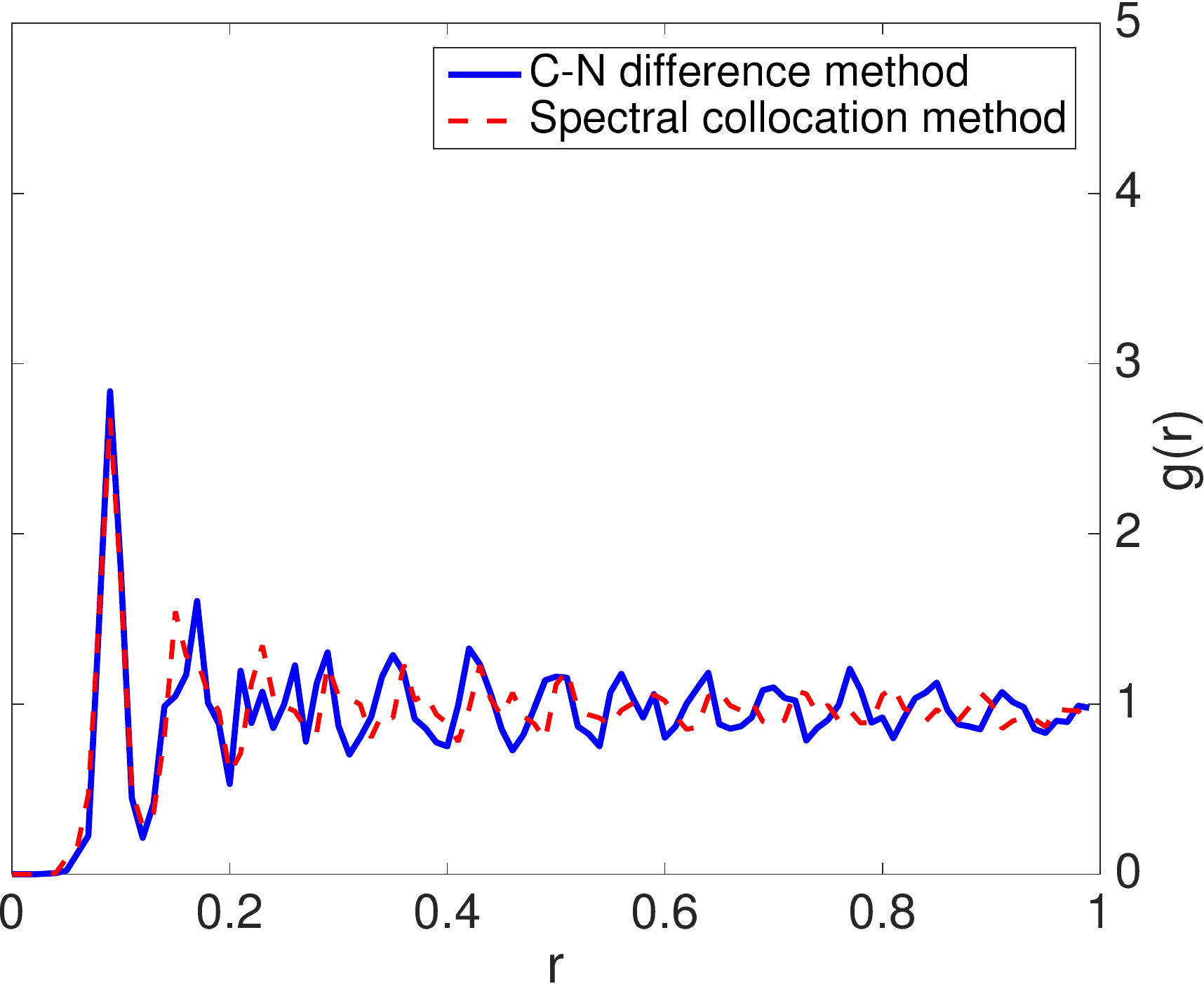}
\end{minipage}}
\subfloat[$\alpha=1.5, \kappa_{2}=0.055$.]{
\begin{minipage}{.31\textwidth}\centering
\includegraphics[width=1.0\textwidth]{d1555}\\
\includegraphics[width=1.0\textwidth]{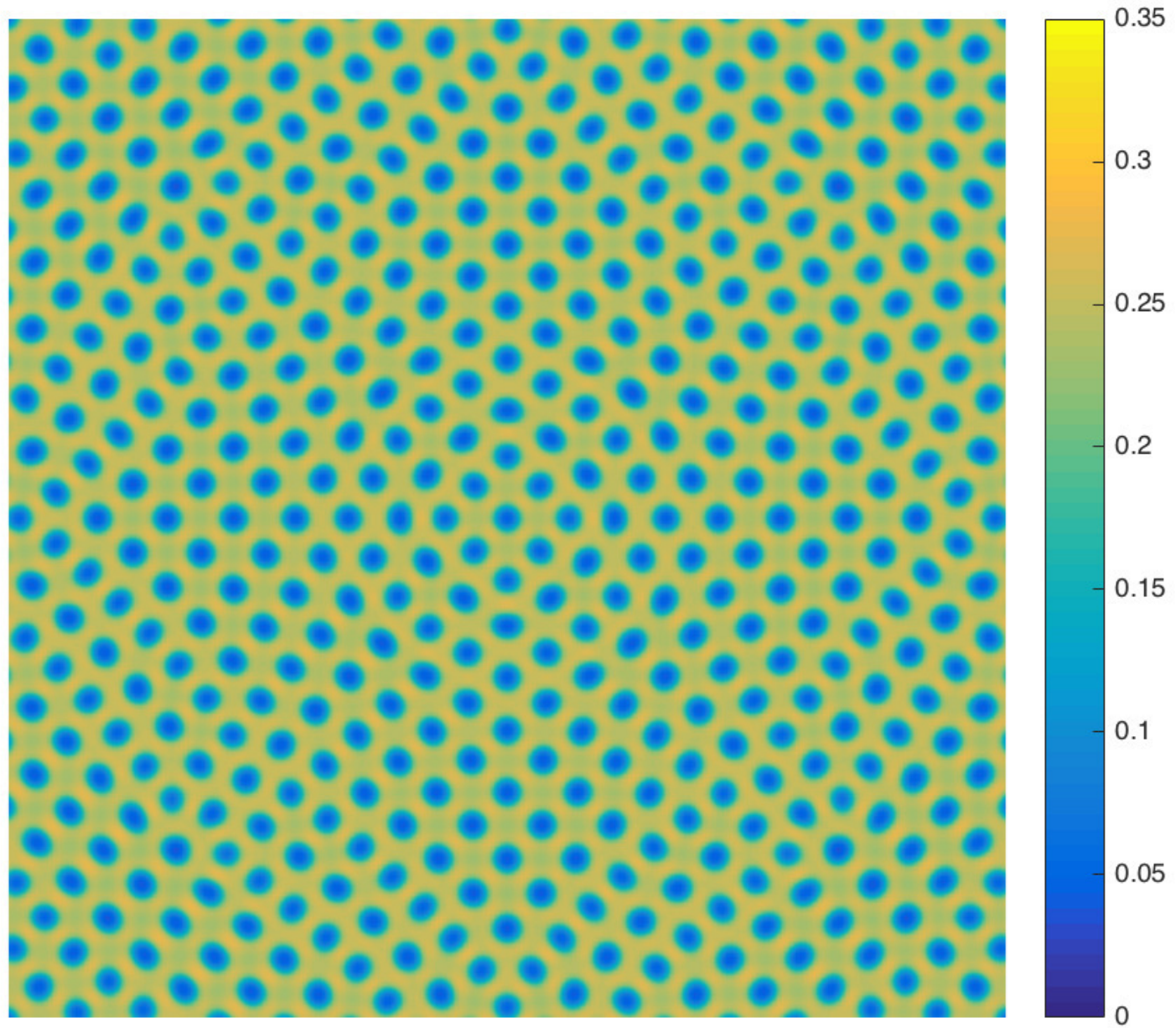}\\
\includegraphics[width=1.0\textwidth]{gr1555}
\end{minipage}}
\caption{ Steady patterns of the solution $v$ with $\kappa_{2}=0.055$ and their corresponding RDFs: $Top$: $v$ contours obtained by C-N difference scheme with weighted shifted Gr\"{u}nwald difference operators at steady states with different fractional orders. The fractional orders in (a)-(c) correspond to $\alpha=2.0, 1.7, 1.5$. $Middle$: $v$ contours obtained by spectral collocation method at steady states with different fractional orders.  The fractional orders in (a)-(c) correspond to $\alpha=2.0, 1.7, 1.5$. $Bottom$: RDFs $g(r)$: the blue line and the red dash line correspond to the variation of the spot density in steady spot patterns obtained by C-N difference scheme and spectral collocation method with the same $\alpha$ and parameter $\kappa_{2}$, respectively. The different fractional orders in (a)-(c) correspond to $\alpha=2.0, 1.7, 1.5$.}
\label{Fig. b}
\end{figure}

Figs. $\ref{Fig. a}$, $\ref{Fig. b}$ show the steady patterns of the numerical solutions $v$ obtained by different methods with parameter $\kappa_{1}=0.063, \kappa_{2}=0.055$ and their corresponding RDFs. All patterns in the two figures are snapshots from the solutions $v$ in the domain $(0,1)^2$. In Fig. $\ref{Fig. a}$, the top contours and the middle contours are the numerical results of the GS model by using the difference scheme and spectral collocation method with $\alpha=2.0, 1.7, 1.5$ and $\kappa_{1}=0.063$, respectively. The bottom figures are corresponding RDFs in which the blue line describes the spot density on the top contour, while the red dash line represents the spot density obtained by the collocation method. For $\alpha=2.0$, we observe that the numerical results obtained from the two different methods are almost the same. When the distance $r_{1}$ is 0.16, we compute that almost 4.5 spots are within this distance away from a spot in two contours. The spot number is around 1.8 with $r_{2}=0.3$. When $\alpha=1.7$, the steady patterns seem the same and the two RDFs are fitting well. We observe that $r_{1}$ decreases when $\alpha$ decreases. For $\alpha=1.5$, the final steady patterns are the same. In Fig. $\ref{Fig. b}$, due to the different $\kappa_{2}=0.055$, we obtain the different steady patterns from Fig \ref{Fig. a} with $\kappa_{1}=0.063$ for the same fractional orders. We observe that the steady patterns are almost the same obtained by the two different numerical methods for $\alpha=2.0, 1.7, 1.5$. The corresponding RDFs also further verify this point. When changing $\alpha$ from 2.0 to 1.5, the spot density increases which represents the decreasing of the distance between spot pairs. Consistently, the distance $r_{1}$ in the bottom figure decreases from 0.16 to 0.05. Consequently, the comparison between the steady patterns
and between their corresponding RDFs obtained by different numerical methods illustrates the accuracy of the difference scheme.

\subsection{The fractional power law of the RDFs}
In this subsection, we investigate the relationship between the fractional order $\alpha$ and $r_{1},r_{2}$ corresponding to the spot distance of the first and the second peak values in the RDFs. $r_{1},r_{2}$ are shown in Fig. $\ref{Fig. a}$ and Fig. $\ref{Fig. b}$. 

In our investigation, we simulate the fractional GS model with several fractional orders $\alpha=1.5, 1.6,\cdots, 2.0$ and different parameters $\kappa_{1}=0.063, \kappa_{2}=0.055$ by using the difference scheme. The domain $\Omega$, parameters $\mu_{u}, \mu_{v}$ and $F$ are the same as the simulations in subsection 5.2. The spatial mesh size is $h=\frac{1}{2000}$ and time step is $\tau=0.1$. Then we compute the corresponding RDFs with different fractional orders and estimate the spot distances $r_{1},r_{2}$. 
\begin{figure}[htp]
\centering
\subfloat[$\kappa_{1}=0.063.$]{
\begin{minipage}{.40\textwidth}\centering
\includegraphics[width=1.0\textwidth]{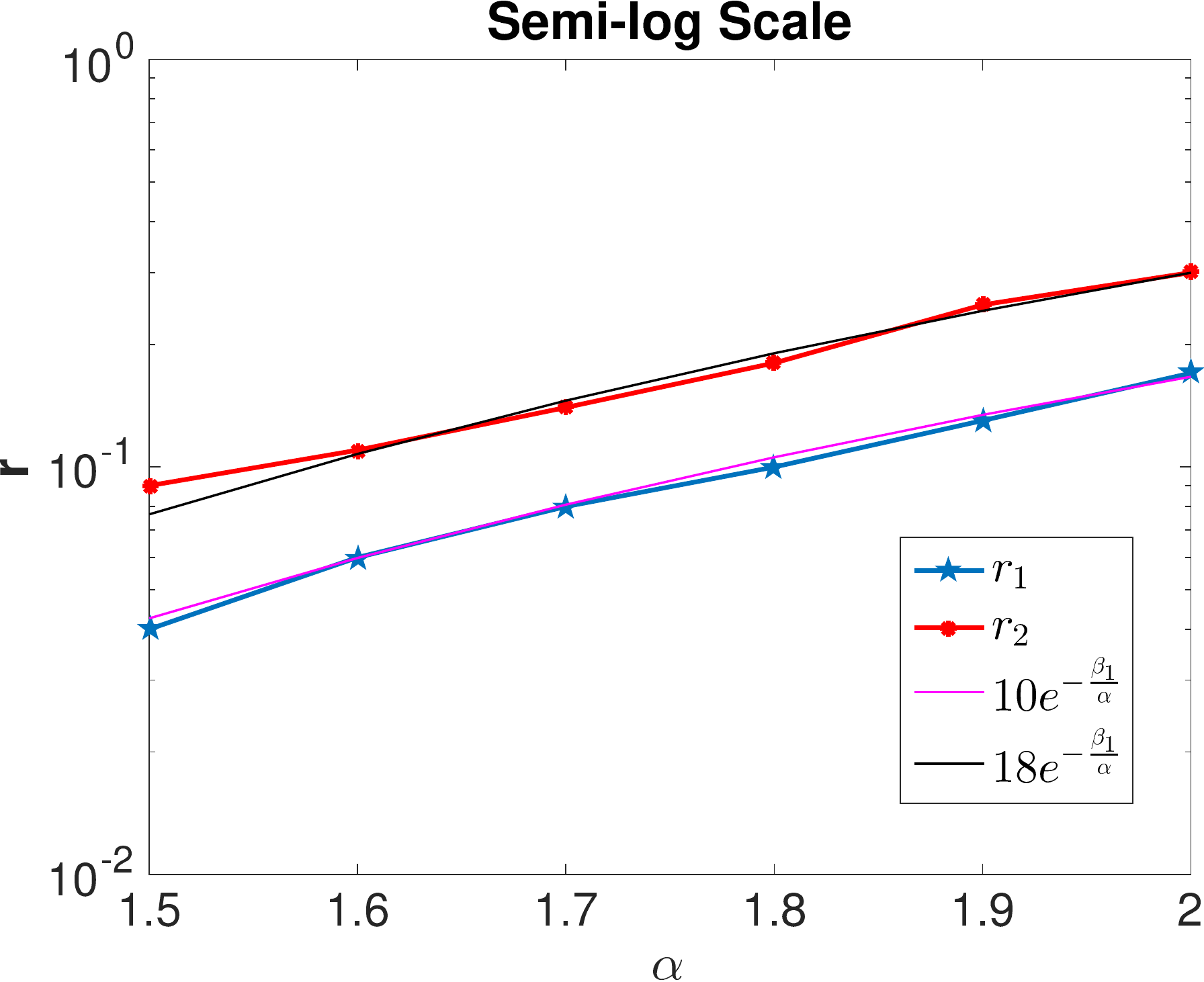}
\end{minipage}}
 \subfloat[$\kappa_{2}=0.055$.]{
\begin{minipage}{.40\textwidth}\centering
\includegraphics[width=1.0\textwidth]{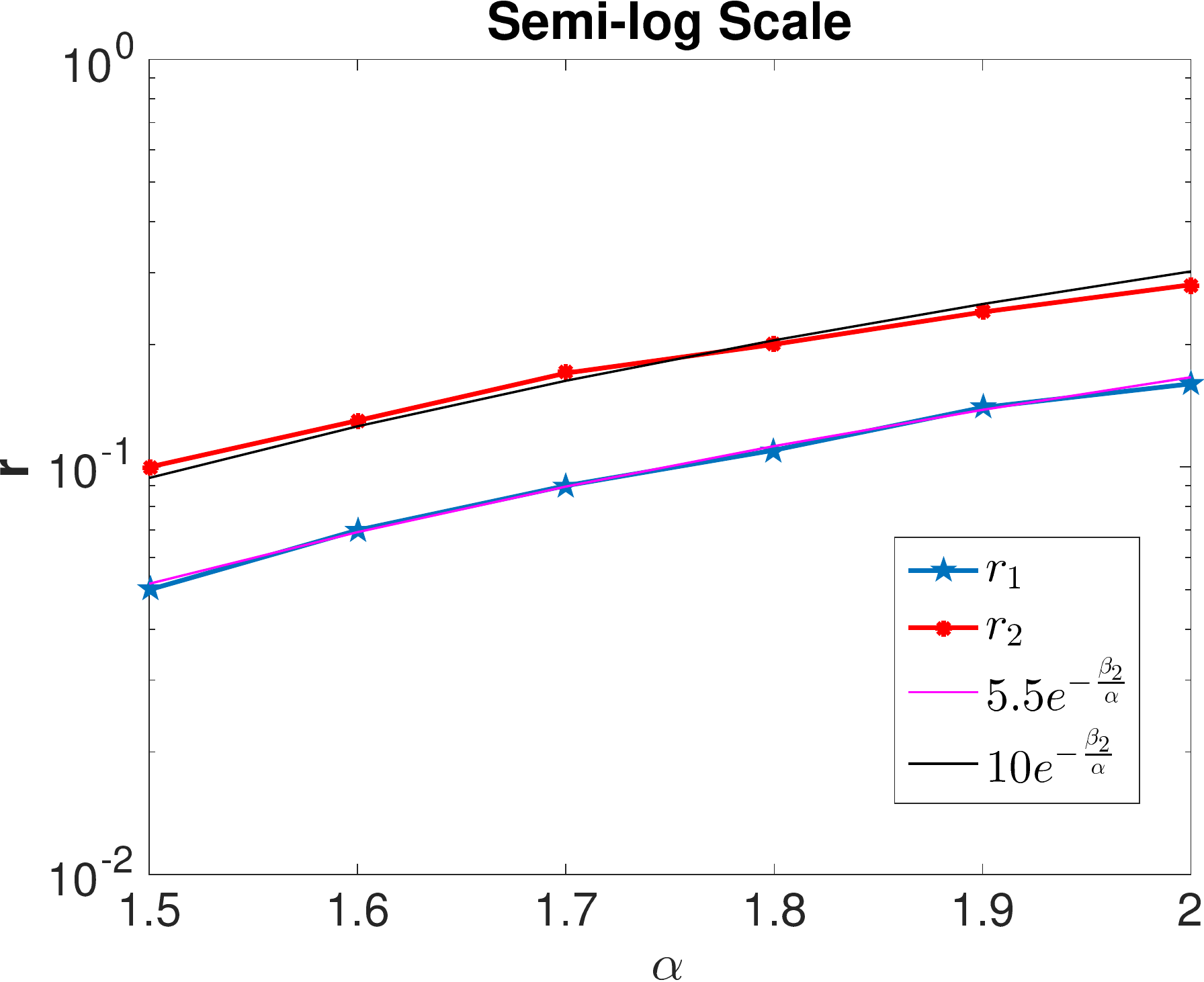}
\end{minipage}}
\caption{The semi-log scales between $\alpha$ and $r_{1},r_{2}$ with different parameters $\kappa_{1}=0.063,\kappa_{2}=0.055$: In the two figures, the light blue star line and the red dot line represent the semi-log relationships between $\alpha$ and $r_{1},r_{2}$, respectively. The pink line and the black line represent the log of the exponential functions with different constants $\beta_{1},\beta_{2}$.}
\label{Fig. end}
\end{figure}

Fig. $\ref{Fig. end}$ plots the spot distances $r_{1},r_{2}$ as functions of the fractional order $\alpha$, which can be fitted with an exponential function as follows
\begin{equation}\label{s5:s4}
r \sim O(e^{-\frac{\beta}{\alpha}}), \quad \quad \alpha \in (1,2],
\end{equation}
where the constant $\beta$ can be obtained numerically. We estimate that $\beta=\beta_{1}=8.19$ for $\kappa_{1}=0.063$, while $\beta=\beta_{2}=7.0$ for $\kappa_{2}=0.055$. This figure implies that there exists a scaling law between the fractional order $\alpha$ and the spot distances $r_{1},r_{2}$ which appear in the corresponding RDFs.

%%%%%%%%%%%%%%%%%%%%%%%%%%%%%%%%%%
\section{Summary}
\setcounter{equation}{0}

We develop a numerical algorithm for the GS model with fractional diffusion and investigate the formation of patterns. We analyze the homogeneous system without diffusion and obtain three steady states including one trivial point and additional two non-trivial points. In order to investigate the stability of steady states, we introduce perturbations to the system and derive the characteristic equation. By analyzing their eigenvalues, we observe that the trivial steady state is stable for all $F$ and $\kappa$, while the stability of the other two states depends on the values of the parameters. A Hopf bifurcation occurs when the parameters vary in a certain range. In addition, we prove the well-posedness of the fractional GS model. We use a C-N difference scheme for time discretization and weighted shifted Gr\"unwald difference operators approximation for space discretization to simulate this model. Moreover, the stability analysis for the time semi-discrete numerical scheme has been provided. We conduct two numerical experiments with the benchmark problems to verify the second-order convergence of this numerical scheme both in time and space. We investigate the pattern formation with various fractional orders, and observe the different effects of the super-diffusion ($1<\alpha <2$) and normal diffusion ($\alpha=2$). Furthermore, we use the spectral collocation method for space discretization to verify the accuracy of the patterns. The comparison of the RDFs further indicates the correctness of the numerical results. We analyze the scaling law 
for steady patterns from the RDFs in terms of the fractional order $1<\alpha \leq 2$.

%%%%%%%%%%%%%%%%%%%%%%%%%%%%%%%%%%%%%%%%%
\appendix
\section{The proof of the Theorem \ref{s3:t1}}\label{adx1}

\begin{proof}
We take the inner product of the first equation of $(\ref{s1:e1})$ with $u(t)$. From $(\ref{s3:ee2})$ and 
the identity $F\bigl((1-u),u\bigr)=\frac{F}{2}\bigl(|\Omega|-\|1-u\|^2-\|u\|^2\bigr)$, we have
\begin{equation}\label{s3:e1}
\frac{1}{2}\Bigl(\frac{d}{dt}\|u\|^{2}+2\mu_{u}\| D^{\frac{\alpha}{2}}u \|^{2}+F\|u\|^{2}\Bigr)
=-\frac{F}{2}\|1-u\|^{2}-\|uv\|^{2}+\frac{F}{2} |\Omega|  \leq \frac{F}{2} |\Omega|,
\end{equation}
which leads to 
\begin{equation}\label{s3:e2}
\frac{d}{dt}\|u\|^{2} +F\|u\|^{2} \leq F |\Omega| .
\end{equation}
By solving this ordinary differential inequality, we obtain $(\ref{tt:1})$. 

Next we add up the first and the second equation of $(\ref{s1:e1})$ to get the new equation $(\ref{tt:3})$ satisfied by $w(t)=u(t)+v(t)$. Taking the inner product of $(\ref{tt:3})$ with $w(t)$ and using $(\ref{s3:ee2})$, we get 
\begin{equation}\label{s3:e6}
\begin{aligned}
\frac{1}{2}\frac{d}{dt}\|w\|^{2}&+\mu_{v}\| D^{\frac{\alpha}{2}}w \|^{2}+(F+\kappa)\|w\|^{2}
=\int_{\Omega}[(\mu_{v}-\mu_{u})(-\Delta)^{\frac{\alpha}{2}}u+\kappa u+F]wdx. 
\end{aligned}
\end{equation}
By using the definition of the fraction Laplacian in (\ref{s2:ee1}), the first term on the right hand side of the above equation (\ref{s3:e6}) can be rewritten as follows 
\begin{equation}\label{s3:ee6}
\begin{aligned}
(&\mu_{v}-\mu_{u})((-\Delta)^{\frac{\alpha}{2}}u,~w) \\
&\leq \frac{|\mu_{v}-\mu_{u}|}{2|\cos(\frac{\pi \alpha}{2})|}
\Bigl(|({}_{a}D_{x}^{\alpha}u,~w)|+|({}_{x}D_{b}^{\alpha}u,~w)|+|({}_{c}D_{y}^{\alpha}u,~w)|+|({}_{y}D_{d}^{\alpha}u,~w)|\Bigr).
\end{aligned}
\end{equation}
Next, we bound the first term on the right hand side by
{\small\begin{equation}\label{s3:ee7}\begin{aligned}
\frac{|\mu_{v}-\mu_{u}|}{2|\cos(\frac{\pi \alpha}{2})|}|({}_{a}D_{x}^{\alpha}u,~w)|
&=\frac{|\mu_{v}-\mu_{u}|}{2|\cos(\frac{\pi \alpha}{2})|}|({}_{a}D_{x}^{\frac{\alpha}{2}}u,~{}_{x}D_{b}^{\frac{\alpha}{2}}w)|
\leq \frac{|\mu_{v}-\mu_{u}|}{2|\cos(\frac{\pi \alpha}{2})|}\|{}_{a}D_{x}^{\frac{\alpha}{2}}u\| \|{}_{x}D_{b}^{\frac{\alpha}{2}}w\|\\
&\leq \frac{\mu_{v}}{4}\|{}_{x}D_{b}^{\frac{\alpha}{2}}w\|^{2}+
\frac{|\mu_{v}-\mu_{u}|^2}{4\mu_{v}|\cos(\frac{\pi \alpha}{2})|^2}\|{}_{a}D_{x}^{\frac{\alpha}{2}}u\|^2.
\end{aligned}\end{equation}}
Similarly, the other three terms are bounded. Then we have the following inequality
\begin{equation}\label{s3:ee8}
(\mu_{v}-\mu_{u})((-\Delta)^{\frac{\alpha}{2}}u,~w)\leq \frac{\mu_{v}}{2}\|D^{\frac{\alpha}{2}}w\|^{2}+
\frac{|\mu_{v}-\mu_{u}|^2}{2\mu_{v}|\cos(\frac{\pi \alpha}{2})|^2}\|D^{\frac{\alpha}{2}}u\|^2.
\end{equation}
Substituting (\ref{s3:ee8}) into (\ref{s3:e6}), and canceling like terms, we get
\begin{equation}\label{s3:e7}
\begin{aligned}
\frac{d}{dt}\|w\|^{2}&+\mu_{v}\| D^{\frac{\alpha}{2}}w \|^{2}+(F+\kappa)\|w\|^{2}
\leq \frac{|\mu_{v}-\mu_{u}|^{2}}{\mu_{v}|\cos(\frac{\pi \alpha}{2})|^2} \| D^{\frac{\alpha}{2}}u\|^{2} +\kappa \|u\|^{2}+F|\Omega|\\
& \leq \frac{|\mu_{v}-\mu_{u}|^{2}}{\mu_{v}|\cos(\frac{\pi \alpha}{2})|^2} \| D^{\frac{\alpha}{2}}u\|^{2} +\kappa e^{-Ft}\|u_{0}\|^{2}+(F+\kappa)|\Omega|.
\end{aligned}
\end{equation}
Integrating the inequality $(\ref{s3:e7})$ in time, we have
\begin{equation}\label{s3:e8}
\|w(t)\|^{2}\leq \|u_{0}+v_{0}\|^{2}+\frac{|\mu_{v}-\mu_{u}|^{2}}{\mu_{v}|\cos(\frac{\pi \alpha}{2})|^2}\int_{0}^{t}\|D^{\frac{\alpha}{2}}u(s)\|^{2}ds
+\frac{\kappa}{F}\|u_{0}\|^{2}+(F+\kappa)|\Omega|t.
\end{equation}
Moreover, combining $(\ref{s3:e1})$ and $\eqref{s3:e8}$,
the estimate $(\ref{tt:2})$ holds.

Since $\|v(t)\|\leq \|w(t)\|+\|u(t)\|$, we conclude that $\|v(t)\|$ is bounded.
\end{proof}

\section{The proof of the Theorem \ref{s4:l1}}\label{adx2}
\begin{proof}
Firstly, from equation $(\ref{s3:e10})$, for $n=0$, we have
\begin{equation}\label{add1}
\frac{U^{1}-u_{0}}{\tau}=-\mu_{u}(-\Delta)^{\frac{\alpha}{2}}U^{\frac{1}{2}}-U^{\frac{1}{2}}v_{0}^2+F(1-U^{\frac{1}{2}}).
\end{equation}
Taking the inner product of $(\ref{add1})$ with $2\tau U^{\frac{1}{2}}$ and using the similar techniques in Appendix \ref{adx1}, we derive
\begin{equation}\label{add2}
\|U^{1}\|^{2}+2\tau \mu_{u}\|D^{\frac{\alpha}{2}}U^{\frac{1}{2}}\|^2+2\tau \|U^{\frac{1}{2}}v_{0}\|^{2}
\leq \|u_{0}\|^{2}+\tau F|\Omega|.
\end{equation}
Then we take the inner product of $(\ref{s3:e10})$ with $2\tau U^{n+\frac{1}{2}}$ for $1\leq n\leq M-1$ and using $(\ref{s3:ee2})$ similarly to get
\begin{equation}\label{s3:e12}
\begin{aligned}
\|U^{n+1}\|^{2}&+2\tau \mu_{u}\|D^{\frac{\alpha}{2}}U^{n+\frac{1}{2}}\|^{2}+
2\tau \|U^{n+\frac{1}{2}}V^{*,n+\frac{1}{2}}\|^{2}\\
&=\|U^{n}\|^{2}-\tau F\|1-U^{n+\frac{1}{2}}\|^{2}-\tau F\|U^{n+\frac{1}{2}}\|^{2}+\tau F|\Omega|\\
& \leq \|U^{n}\|^{2}+\tau F|\Omega|.
\end{aligned}\end{equation}
We sum the inequality $(\ref{s3:e12})$ from $k=1$ to $k=n$, $n\leq M-1$, to get
\begin{equation}\label{s3:e13}
\|U^{n+1}\|^{2}+2\tau \mu_{u}\sum_{k=1}^{n}\| D^{\frac{\alpha}{2}}U^{k+\frac{1}{2}}\|^{2}+
2\tau \sum_{k=1}^{n}\|U^{k+\frac{1}{2}}V^{*,k+\frac{1}{2}}\|^{2} \leq \|U^{1}\|^{2}+n\tau F|\Omega|.
\end{equation}
Thus combining the $(\ref{add2})$ with $(\ref{s3:e13})$, for $0\leq n \leq M-1$, we obtain
\begin{equation}\label{add3}
\|U^{n+1}\|^{2}+2\tau \mu_{u}\sum_{k=0}^{n}\| D^{\frac{\alpha}{2}}U^{k+\frac{1}{2}}\|^{2}
 \leq \|u_{0}\|^{2}+ FT|\Omega|,
\end{equation}
which leads to the estimate $(\ref{tt:4})$.

Next adding $(\ref{s3:e10})$ and $(\ref{s3:e11})$, we derive a new equation $(\ref{tt:6})$ satisfied by $W^{n}=U^{n}+V^{n}$. Taking the inner product of $(\ref{tt:6})$ with $2\tau W^{n+\frac{1}{2}}$ for $0\leq n\leq M-1$, we have
\begin{equation}\label{s3:e15}
\begin{aligned}
\|W^{n+1}\|^{2}-\|W^{n}\|^{2}=&-2\tau\mu_{v}\|D^{\frac{\alpha}{2}}W^{n+\frac{1}{2}}\|^2
-2\tau(F+\kappa)\|W^{n+\frac{1}{2}}\|^{2}\\
&+2\tau \int_{\Omega}\Big[(\mu_{v}-\mu_{u})(-\Delta)^{\frac{\alpha}{2}}U^{n+\frac{1}{2}}+\kappa U^{n+\frac{1}{2}}+F\Bigr]W^{n+\frac{1}{2}}dx.
\end{aligned}\end{equation}
using the same techniques from (\ref{s3:ee6}) to (\ref{s3:ee8}), we derive
\begin{equation*}\begin{aligned}
2\tau \int_{\Omega}\Big[(\mu_{v}-\mu_{u})&(-\Delta)^{\frac{\alpha}{2}}U^{n+\frac{1}{2}}+\kappa U^{n+\frac{1}{2}}+F\Bigr]W^{n+\frac{1}{2}}dx\\
&\leq \tau \mu_{v} \|D^{\frac{\alpha}{2}}W^{n+\frac{1}{2}}\|^2+\tau\frac{|\mu_{v}-\mu_{u}|^2}{\mu_{v}|\cos(\frac{\pi \alpha}{2})|^2}\|D^{\frac{\alpha}{2}}U^{n+\frac{1}{2}}\|^2\\
&\quad +\kappa \tau (\|W^{n+\frac{1}{2}}\|^{2}
+\|U^{n+\frac{1}{2}}\|^{2})+F\tau (|\Omega|+\|W^{n+\frac{1}{2}}\|^{2}).
\end{aligned}\end{equation*}
Hence, we get
\begin{equation}\label{s3:e16}
\begin{aligned}
\|W^{n+1}\|^{2}&+\tau\mu_{v}\|D^{\frac{\alpha}{2}}W^{n+\frac{1}{2}}\|^2
+\tau(F+\kappa)\|W^{n+\frac{1}{2}}\|^{2}\\
&\leq \|W^{n}\|^{2}+\tau \frac{|\mu_{v}-\mu_{u}|^{2}}{\mu_{v}|\cos(\frac{\pi \alpha}{2})|^2}\|D^{\frac{\alpha}{2}}U^{n+\frac{1}{2}}\|^{2}
+\kappa \tau \|U^{n+\frac{1}{2}}\|^{2}+F\tau|\Omega|.
\end{aligned}
\end{equation}
Since $\|U^{n+\frac{1}{2}}\|^2\leq \|u_{0}\|^2+FT|\Omega|$, we have
\begin{equation}\label{s3:e17}
\begin{aligned}
\|W^{n+1}\|^{2}&+\tau\mu_{v}\|D^{\frac{\alpha}{2}}W^{n+\frac{1}{2}}\|^2
+\tau(F+\kappa)\|W^{n+\frac{1}{2}}\|^{2}\\
&\leq \|W^{n}\|^{2}+\tau \frac{|\mu_{v}-\mu_{u}|^{2}}{\mu_{v}|\cos(\frac{\pi \alpha}{2})|^2}\|D^{\frac{\alpha}{2}}U^{n+\frac{1}{2}}\|^{2}
+\kappa \tau (\|u_{0}\|^{2}+FT|\Omega|)+F\tau|\Omega|.
\end{aligned}
\end{equation}
We sum the inequality $(\ref{s3:e17})$ from $k=0$ to $k=n$, $n\leq M-1$, and cancel like terms to obtain
\begin{equation}\label{s3:e18}
\|W^{n+1}\|^{2}\leq \|u_{0}+v_{0}\|^{2}+\tau \frac{|\mu_{v}-\mu_{u}|^{2}}{\mu_{v}|\cos(\frac{\pi \alpha}{2})|^2}\sum_{k=0}^{n}\|D^{\frac{\alpha}{2}}U^{k+\frac{1}{2}}\|^{2}+\kappa T (\|u_{0}\|^{2}+FT|\Omega|)+FT|\Omega|.
\end{equation}
Moreover, from $(\ref{add3})$ we see that
\begin{equation*}
\tau\sum_{k=0}^{n}\|D^{\frac{\alpha}{2}}U^{k+\frac{1}{2}}\|^{2}\leq
\frac{1}{2\mu_{u}}\|u_{0}\|^{2}+\frac{1}{2\mu_{u}}FT|\Omega|.
\end{equation*}
Hence, the estimate $(\ref{tt:5})$ holds.

Since $\|W^{n}\|$ and $\|U^{n}\|$ have been bounded, we use $\|V^{n}\|\leq \|W^{n}\|+\|U^{n}\|$ to conclude that $\|V^{n}\|$ is bounded, for $0\leq n\leq M$. We thus prove the estimates of the scheme.
\end{proof}

%%%%%%%%%%%%%%%%%%%%%%%%%%%%%%%%%%%%%%%%%

{\bf Acknowledgements.} This work was supported by the OSD/ARO/MURI on ``Fractional PDEs for Conservation Laws and Beyond: Theory, Numerics and Applications (W911NF-15-1-0562)" and the National Science Foundation under Grant DMS-1620194. The first author was supported by the China Scholarship Council under 201706220157. 

\bibliographystyle{abbrv}
\bibliography{elsarticle-template-num}
\end{document}